\documentclass{amsart}
\usepackage{amsmath}
\usepackage{amssymb}
\usepackage{extarrows}
\usepackage{mathabx}
\usepackage[all]{xy}
\usepackage{amsbsy}
\usepackage{graphicx}
\usepackage{appendix}
\usepackage{float}
\usepackage{subfigure}
\usepackage[numbers,sort&compress]{natbib}
\usepackage{graphicx}
\usepackage{amsthm}
\usepackage{geometry}
\usepackage{fancyhdr}
\usepackage{color} 
\usepackage{overpic}
\usepackage{mathabx}
\usepackage{tikz-cd}
\usepackage{enumerate}
\usepackage{mathrsfs}
\usepackage{colonequals}
\usepackage{microtype}

\newtheorem{thm}{Theorem}[section]
\newtheorem{cor}[thm]{Corollary}
\newtheorem{prop}[thm]{Proposition}
\newtheorem{lem}[thm]{Lemma}

\theoremstyle{definition}
\newtheorem{defn}[thm]{Definition}
\newtheorem{exmp}[thm]{Example}
\newtheorem{conj}[thm]{Conjecture}
\newtheorem*{convs}{Conventions}
\newtheorem*{conv}{Convention}

\newtheorem*{org}{Organization}
\newtheorem*{ack}{Acknowledgement}

\theoremstyle{remark}
\newtheorem{rem}[thm]{Remark}

\numberwithin{equation}{section}

\newcommand{\beq}{\begin{equation*}\begin{aligned}}
\newcommand{\eeq}{\end{aligned}\end{equation*}}
\newcommand{\bpf}{\begin{proof}}
\newcommand{\epf}{\end{proof}}
\newcommand{\bthm}{\begin{thm}}
\newcommand{\ethm}{\end{thm}}
\newcommand{\bprop}{\begin{prop}}
\newcommand{\eprop}{\end{prop}}
\newcommand{\bcor}{\begin{cor}}
\newcommand{\ecor}{\end{cor}}
\newcommand{\blem}{\begin{lem}}
\newcommand{\elem}{\end{lem}}
\newcommand{\bdefn}{\begin{defn}}
\newcommand{\edefn}{\end{defn}}
\newcommand{\bexmp}{\begin{exmp}}
\newcommand{\eexmp}{\end{exmp}}
\newcommand{\brem}{\begin{rem}}
\newcommand{\erem}{\end{rem}}
\newcommand{\benu}{\begin{enumerate}[(1)]}
\newcommand{\eenu}{\end{enumerate}}

\newcommand{\bdia}{\begin{displaymath}\xymatrix}
\newcommand{\edia}{\end{displaymath}}

\newcommand{\hfk}{\widehat{HFK}}
\newcommand{\hf}{\widehat{HF}}

\newcommand{\shib}{{\mathbf{SHI}}}
\newcommand{\shi}{\underline{\rm SHI}}
\newcommand{\shf}{\mathbf{SHF}}

\newcommand{\shg}{\mathbf{SH}}

\newcommand{\khi}{\underline{\rm KHI}^-}
\newcommand{\khg}{{\mathbf{KH}}}
\newcommand{\khib}{\mathbf{KHI}}
\newcommand{\sh}{\mathbf{SH}}
\newcommand{\hg}{\mathbf{H}}
\newcommand{\hhat}{\widetilde{\mathbf{H}}}

\newcommand{\deq}{\colonequals}
\newcommand{\sym}{{\rm Sym}}
\newcommand{\cfi}{{\bf CF}^\infty}
\newcommand{\hfi}{{\bf HF}^\infty}
\newcommand{\cfm}{{\bf CF}^-}
\newcommand{\hfm}{{\bf HF}^-}
\newcommand{\spin}{{\rm Spin}^c}

\newcommand{\al}{\alpha}
\newcommand{\be}{\beta}
\newcommand{\ga}{\gamma}
\newcommand{\Ga}{\Gamma}

\newcommand{\ot}{\otimes}

\newcommand{\p}{\prime}
\newcommand{\pp}{{\prime\prime}}
\newcommand{\bs}{\boldsymbol}
\newcommand{\gr}{\chi_{\rm gr}}

\newcommand{\aand}{~{\rm and}~}

\newcommand{\uw}{{U_{\bs{w}}}}
\newcommand{\uww}{{U_{\bs{w}}^{-1}}}

\newcommand{\mcr}{\mathcal{R}}

\newcommand{\mch}{\mathcal{H}}

\newcommand{\intg}{\mathbb{Z}}
\newcommand{\ft}{{\mathbb{F}_2}}

\newcommand{\ra}{\rightarrow}

\newcommand{\xra}{\xrightarrow}

\begin{document}

\title{Instanton Floer homology, sutures, and Euler characteristics}


\author{Zhenkun Li}
\address{Department of Mathematics, Stanford University}
\curraddr{}
\email{zhenkun@stanford.edu}
\thanks{}

\author{Fan Ye}
\address{Department of Pure Mathematics and Mathematical Statistics, University of Cambridge}
\curraddr{}
\email{fanye@math.harvard.edu}
\thanks{}

\date{}
\dedicatory{}
\begin{abstract}
This is a companion paper to an earlier work of the authors. In this paper, we provide an axiomatic definition of Floer homology for balanced sutured manifolds and prove that the graded Euler characteristic $\chi_{\rm gr}$ of this homology is fully determined by the axioms we proposed. As a result, we conclude that $\chi_{\rm gr}(SHI(M,\gamma))=\chi_{\rm gr}(SFH(M,\gamma))$ for any balanced sutured manifold $(M,\gamma)$. In particular, for any link $L$ in $S^3$, the Euler characteristic $\chi_{\rm gr}(KHI(S^3,L))$ recovers the multi-variable Alexander polynomial of $L$, which generalizes the knot case. Combined with the authors' earlier work, we provide more examples of $(1,1)$-knots in lens spaces whose $KHI$ and $\widehat{HFK}$ have the same dimension. Moreover, for a rationally null-homologous knot in a closed oriented 3-manifold $Y$, we construct canonical $\mathbb{Z}_2$-gradings on $KHI(Y,K)$, the decomposition of $I^\sharp(Y)$ discussed in the previous paper, and the minus version of instanton knot homology $\underline{\rm KHI}^-(Y,K)$ introduced by the first author.
\end{abstract}
\maketitle
\tableofcontents
\newpage

\section{Introduction}

Sutured manifold theory was introduced by Gabai \cite{gabai1983foliations}, and Floer theory was introduced by Floer \cite{floer1988instanton}. They are both powerful tools in the study of 3-manifolds and knots. The first combination of these theories, called sutured Floer homology, was introduced by Juh\'asz \cite{juhasz2006holomorphic} based on Heegaard Floer theory, with some pioneering work done by Ghiggini \cite{ghiggini2008knot} and Ni \cite{ni2007knot}. Later, Kronheimer and Mrowka made analogous constructions in monopole (Seiberg-Witten) theory and instanton (Donaldson-Floer) theory \cite{kronheimer2010knots}. Different versions of Floer theories have different merits. For example, Heegaard Floer theory is more computable, while instanton theory is closely related to representation varieties of fundamental groups. Hence it is important to understand the relationship between different versions of Floer theories. In this line, Lekili \cite{lekili2013heegaard} and Baldwin and Sivek \cite{Baldwin2020} proved that sutured (Heegaard) Floer homology is isomorphic to sutured monopole homology, though the relation to sutured instanton homology is still open.

\begin{conj}[{\cite[Conjecture 7.24]{kronheimer2010knots}}]\label{conj: I=HF}
    For a balanced sutured manifold $(M,\ga)$, we have $$SHI(M,\ga)\cong SFH(M,\ga)\otimes \mathbb{C}.$$In particular, for a knot $K$ in a closed oriented 3-manifold $Y$, there are isomorphisms
    $$I^{\sharp}(Y)\cong \widehat{HF}(Y)\otimes\mathbb{C}~{\rm and}~KHI(Y,K)\cong \widehat{HFK}(Y,K)\otimes \mathbb{C}.$$
Here $SHI$ is sutured instanton homology \cite{kronheimer2010knots}, $SFH$ is sutured (Heegaard) Floer homology \cite{juhasz2006holomorphic}, $I^{\sharp}$ is framed instanton Floer homology \cite{kronheimer2011khovanov}, $\widehat{HF}$ is the hat version of Heegaard Floer homology \cite{ozsvath2004holomorphic}, $KHI$ is instanton knot homology \cite{kronheimer2010knots}, and $\widehat{HFK}$ is the hat version of knot Floer homology \cite{ozsvath2004holomorphicknot,Rasmussen2003}.
\end{conj}

In this paper, instead of studying the full homologies, we study their graded Euler characteristics and obtain the following theorem.


\bthm\label{thm: main}
Suppose $(M,\ga)$ is a balanced sutured manifold and $S_1,\dots,S_n$ are properly embedded admissible surfaces (\textit{c.f.} Definition \ref{defn_2: admissible surfaces}) generating $H_2(M,\partial M)/{\rm Tors}$. Then there exist $\intg^n$-gradings on $SHI(M,\ga)$ and $SFH(M,\ga)$ induced by these surfaces, respectively. Equivalently, we have $$SHI(M,\ga)=\bigoplus_{(i_1,\dots,i_n)\in \mathbb{Z}^n}SHI(M,\ga,(S_1,\dots,S_n),(i_1,\dots,i_n))$$and a similar result holds for $SFH(M,\ga)$. Moreover, there exist relative $\mathbb{Z}_2$-gradings on $SHI(M,\ga)$ and $SFH(M,\ga)$, respecting the decompositions. Define 
\begin{equation}\label{eq: defn of chi}
    \chi_{\rm gr}(SHI(M,\ga))\deq \sum_{(i_1,\dots,i_n)\in\intg^n}\chi(SHI(M,\ga,(S_1,\dots,S_n),(i_1,\dots,i_n))\cdot t_1^{i_1}\cdots t_n^{i_n},
\end{equation}
and define $\chi_{\rm gr}(SFH(M,\ga))$ similarly. Then we have
$$\chi_{\rm gr}(SHI(M,\ga))\sim \chi_{\rm gr}(SFH(M,\ga)),$$
where $\sim$ means two polynomials equal up to multiplication by $\pm t_1^{j_1}\cdots t_n^{j_n}$ for some $(j_1,\dots,j_n)\in\mathbb{Z}^n$, 
\ethm
\brem
Suppose that $t_1,\dots,t_n$ represent generators of $$H=H_1(M;\mathbb{Z})/{\rm Tors}\cong H_2(M,\partial M;\mathbb{Z})/{\rm Tors}.$$Then $\sim$ means the equality holds for elements in $\mathbb{Z}[H]/\pm H$.
\erem

The graded Euler characteristic $\chi_{\rm gr}(SFH(M,\ga))$ was studied by Friedl, Juh{\'{a}}sz, and Rasmussen \cite{Friedl2009}. Applying their results, we can relate the graded Euler characteristics of links with classical invariants obtained from fundamental groups. 

Consider a finitely generated group $\pi=\langle x_1,\dots,x_n|r_1,\dots,r_k\rangle$. Let $H=H_1(\pi)/{\rm Tors}$ be the abelianization of $\pi$ modulo torsions. For a generator $x_i$ and a word $w$, let $\partial w/\partial x_i$ be the Fox derivative of $w$ with respect to $x_i$. Equivalently, it satisfies the following conditions.
\begin{enumerate}
    \item For any word $w=u\cdot v$, we have $\frac{\partial w}{\partial x_i}=\frac{\partial u}{\partial x}+u\cdot\frac{\partial v}{\partial x}$.
    \item $\frac{\partial x_i}{\partial x_i}=1$ and $\frac{\partial x_j}{\partial x_i}=0$ for any $j\neq i$.
\end{enumerate}

Consider $A=\{\partial r_j/\partial x_i\}_{i,j}$ as a matrix with entries in $\mathbb{Z}[H]$ by the projection map $\mathbb{Z}[\pi]\to \mathbb{Z}[H]$. Let $E(\pi)$ be the ideal generated by the minor determinants of $A$ of order $(n-1)$. Since $\mathbb{Z}[H]$ is a unique facterization domain, one can consider the greatest common divisor of the elements of $E(\pi)$, which is well-defined up to multiplication by a unit in $\pm H$. This is denoted by $\Delta(\pi)$ and called the \textbf{Alexander polynomial} of $\pi$ (\textit{c.f.} \cite{Turaev2002}). For a 3-manifold $M$, the Alexander polynomial of $M$ is defined by $\Delta(M)\deq \Delta (\pi_1(M))$. For an $n$-component link $L$ in $S^3$, we write $t_1,\dots,t_n$ for homology classes of meridians of components of $L$ and define $\Delta_L(t_1,\dots,t_n)\deq \Delta(S^3\backslash {\rm int}N(L))$ as the \textbf{multi-variable Alexander polynomial} of $L$. If $n=1$ and $L=K$ is a knot, we can fix the ambiguity of $\pm H$ by assuming $\Delta_K(t_1)=\Delta_K(t^{-1}_1)$ and $\Delta_K(1)=1$. In this case, we call it the \textbf{symmetrized Alexander polynomial} of $K$.

\bthm\label{thm: multi variable Alexander polynomial}Suppose $M$ is a compact manifold whose boundary consists of tori $T_1,\dots,T_n$ with $b_1(M)\ge 2$. Suppose $$\ga=\bigcup_{j=1}^nm_j\cup (-m_j)$$consists of two simple closed curves with opposite orientations on each torus. Suppose $H=H_1(M;\mathbb{Z})/{\rm Tors}$ and $[m_1],\dots,[m_n]$ are homology classes. Then we have
\begin{equation}\label{eq: turaev}
    \chi_{\rm gr}(SHI(M,\ga))=
\Delta(M)\cdot\prod_{j=1}^n([m_j]-1)\in \mathbb{Z}[H]/\pm H.
\end{equation}

In particular, suppose $L\subset S^3$ is an $n$-component link with $n\ge 2$. Define 
\begin{equation}\label{eq: khi(L)}
    KHI(L)\deq SHI(S^3\backslash {\rm int}N(L),\bigcup_{j=1}^nm_j\cup (-m_j)),
\end{equation}where $m_1,\dots,m_n$ are meridians of components of $L$. Let $(i_1,\dots,i_n)$ denote the $\intg^n$-grading on $KHI(L)$ induced by Seifert surfaces of components of $L$. Then we have
\begin{equation*}
    \chi_{\rm gr}(KHI(L))\deq\sum_{(i_1,\dots,i_n)\in\intg^n}\chi(KHI(L,(S_1,\dots,S_n),(i_1,\dots,i_n)))\cdot t_1^{i_1}\cdots t_n^{i_n}\sim
\Delta_L(t_1,\dots,t_n)\cdot\prod_{j=1}^n(t_j-1),
\end{equation*}
where $\sim$ means the equality holds for elements in $\mathbb{Z}[H]/\pm H$.
\ethm
\brem
A similar result to Theorem \ref{thm: multi variable Alexander polynomial} has been proved for link Floer homology in Heegaard Floer theory by Ozsv\'ath and Szab\'o \cite{ozsvath2008multivariable}. For instanton theory, the case of single-variable Alexander polynomial for links in $S^3$ was understood by Kronheimer and Mrowka \cite{kronheimer2010instanton} and independently by Lim \cite{Lim2009}, while the case of the multi-variable polynomial was unknown before.
\erem

For knots, the corresponding corollary is the following.
\bthm\label{thm: knot}
Suppose $K$ is a knot in a closed oriented 3-manifold $Y$. Suppose $Y(K)\deq Y\backslash {\rm int}N(K)$ is the knot complement and $b_1(Y(K))=1$. Let $[m]\in H=H_1(Y(K);\mathbb{Z})/{\rm Tors}\cong \mathbb{Z}\langle t\rangle$ be the homology class of the meridian of $K$. Define $KHI(Y,K)$ similarly to $KHI(L)$ as in (\ref{eq: khi(L)}). Then we have
\begin{equation*}\label{eq: turaev 2}
    \chi_{\rm gr}(KHI(Y,K))=\Delta(Y(K))\cdot \frac{[m]-1}{t-1}\in \mathbb{Z}[H]/\pm H.
\end{equation*}
\ethm

\brem
Analogous results of Theorem \ref{thm: knot} in Heegaard Floer theory can be found in {\cite[Proposition 2.1]{Rasmussen2017} and \cite[Proposition 3.1]{Rasmussen2007}}. Also, Theorem \ref{thm: knot} is a generalization of work of Kronheimer and Mrowka \cite{kronheimer2010instanton} and Lim \cite{Lim2009}, in which they proved the same results only for knots inside $S^3$.
\erem

\brem
Since the first version of this paper, the above theorems have many applications: \begin{enumerate}[(1)]
   \item in \cite{yixiesu2}, Xie and Zhang  used Theorem \ref{thm: multi variable Alexander polynomial} to show that the fundamental group of any link in $S^3$ that is neither the unknot nor the Hopf link admits an irreducible $SU(2)$-representation;
    \item in \cite{LY2021enhanced}, the authors used Theorem \ref{thm: main} to compute the instanton knot homology of some extra families of $(1,1)$-knots and classified all null-homologous instanton Floer minimal knots inside an instanton L-space;
    \item in \cite{BLSY2021su2}, Baldwin, Li, Sivek, and Ye classified all instanton Floer minimal knots inside the lens space $L(3,1)$, which finally contributed to the proof of the theorem that the fundamental group of the $3$-surgery of any nontrivial knot admits an irreducible $SU(2)$-representation.
\end{enumerate}
    
\erem

\brem
In Theorem \ref{thm: main}, we relate the Euler characteristics of $SHI$ and $SFH$. For this purpose, we deal with Heegaard Floer theory directly, and prove some elementary properties like the excision property in Section \ref{sec: hf}. An alternative approach could be to deal with the relation between $SHI$ and $SHM$ first, as all necessary preparations have already been made in \cite{kronheimer2010knots}, and then use the relation between $\widecheck{HM}_\bullet$ and $HF^+$ by work of Colin, Ghiggini, and Honda \cite{colin2011equivalence} and Taubes \cite{taubes2010ech1,taubes2010ech2,taubes2010ech3,taubes2010ech4,taubes2010ech5}, or independently Kutluhan, Lee, and Taubes \cite{kutluhan2010hf,kutluhan2010hf2,kutluhan2010hf3,kutluhan2011hf,kutluhan2012hf}, and the relation between $SHM$ and $SFH$ by work of Lekili \cite{lekili2013heegaard} or independently Baldwin and Sivek \cite{Baldwin2020}.
\erem

An application of Theorem \ref{thm: knot} is to compute the instanton knot homology of some special families of knots. In \cite{LY2020}, the authors proved the following.
\bthm[{\cite[Theorem 1.6]{LY2020}}]\label{thm: inequality for 1,1 knots}
Suppose $K\subset Y$ is a $(1,1)$-knot in a lens space (including $S^3$). Then we have
$${\rm dim}_{\mathbb{C}}KHI(Y,K)\leq\dim_\ft\widehat{HFK}(Y,K).$$
\ethm
Obviously, a lower bound of ${\rm dim}_{\mathbb{C}}KHI(Y,K)$ can be obtained from $\chi_{\rm gr}(KHI(Y,K))$. If this lower bound coincides with the upper bound from Theorem \ref{thm: inequality for 1,1 knots}, then we figure out the precise dimension of $KHI(Y,K)$. This trick applies to $(1,1)$-knots in $S^3$, which are either homologically thin knots or Heegaard Floer L-space knots. In \cite{LY2020}, the authors worked with knots in $S^3$ because prior to the current paper, the graded Euler characteristic of instanton knot homology was only understood in that case. On the other hand, in \cite{Ye2020}, the second author discovered a family of $(1,1)$-knots in general lens spaces whose $\dim_\ft\widehat{HFK}(Y,K)$ is determined by $\chi(\widehat{HFK}(Y,K))$. Hence, with Theorem \ref{thm: knot} and results in \cite{Ye2020}, we conclude the following.

\bcor\label{cor: 11 knot}
Suppose $Y$ is a lens space, and $K\subset Y$ is a $(1,1)$-knot such that .
\benu
\item either $K$ admits an $L$-space surgery (\textit{c.f.} \cite[Lemma 3.2]{Rasmussen2017} and \cite[Theorem 2.2]{Greene2018}), or $K$ is a constrained knot (\textit{c.f.} \cite[Section 4]{Ye2020}),
\item $H_1(Y(K);\mathbb{Z})\cong\mathbb{Z}$, where $Y(K)$ is the knot complement of $K$.
\eenu
Then, we have
$$\dim_{\mathbb{C}}KHI(Y,K)= \dim_{\mathbb{Z}}\widehat{HFK}(Y,K).$$
\ecor
\brem
Greene, Lewallen, and Vafaee \cite{Greene2018} provided a clear criterion to check if a $(1,1)$-knot admits an L-space surgery. 
\erem
\brem
The condition $H_1(Y(K);\mathbb{Z})\cong\mathbb{Z}$ is necessary since terms related to Euler characteristics of torsion spin$^c$ structures may cancel out when we consider the map between group rings induced by the projection $H_1(Y(K);\mathbb{Z})\to H_1(Y(K);\mathbb{Z})/{\rm Tors}$. In a subsequent paper \cite{LY2021enhanced}, we introduced an enhanced Euler characteristic of $SHI$ to deal with the torsion part of $H_1(Y(K);\mathbb{Z})$ and remove the second condition in Corollary \ref{cor: 11 knot}.
\erem

Now, we explain the rough idea to prove Theorem \ref{thm: main}. First, let us consider the case of a closed oriented 3-manifold $Y$. The Euler characteristic of the framed instanton Floer homology, $\chi(I^\sharp (Y))$, was understood by Scaduto \cite[Section 10]{scaduto2015instanton}. The strategy is to carry out an induction on the order of $H_1(Y;\mathbb{Z})$ using exact triangles. The grading behavior of $\chi(I^\sharp (Y))$ under a surgery exact triangle was fully described as in \cite[Section 42.3]{kronheimer2007monopoles} and it is known that $\chi(I^\sharp (Y^\p))=1$ for any integral homology sphere $Y^\p$. Hence we can prove $\chi(I^\sharp (Y))=|H_1(Y;\mathbb{Z})|$ inductively. 

However, the above argument requires more care when we take into account gradings associated to surfaces inside $3$-manifolds. Suppose $R\subset Y$ is a closed homologically essential surface. Then $R$ induces a $\intg$-grading on $I^\sharp (Y)$ by considering the generalized eigenspaces of the linear action $\mu(R)$ on $I^{\sharp}(Y)$ (\textit{c.f.} \cite[Section 7]{kronheimer2010knots}). When trying to understand the graded Euler characteristic in this case, the previous strategy does not apply directly. The reason is that, the surgery curves inducing the exact triangles may have nontrivial algebraic intersections with the surface $R$, so the maps in surgery exact triangles may not preserve the grading associated to $R$. We are faced with the same problem when proving Theorem \ref{thm: main}.

Our strategy is the following. Suppose $(M,\ga)$ is a balanced sutured manifold and suppose $S_1$,\dots, $S_n$ are properly embedded surfaces in $M$. Then $S_1,\dots,S_n$ induce a $\intg^n$-grading on $SHI(M,\ga)$. After attaching product 1-handles along $\partial M$, we can find a framed link in the interior of the resulting manifold such that the link is disjoint from all the surfaces. Moreover, surgeries along the link with all slopes chosen in $\{0,1\}$ produce only handlebodies. Since the surgery link is disjoint from the surfaces $S_1,\dots,S_n$ that induce the $\intg^n$-grading, the maps in surgery exact triangles preserve the grading. Hence, it suffices to understand the case of sutured handlebodies. In this case, we can further use bypass exact triangles to reduce any sutured handlebodies to product sutured manifolds. It is known that the Floer homology of any product sutured manifold is one-dimensional. Since the behavior of Euler characteristics under bypass exact triangles and surgery exact triangles are the same for both instanton theory and Heegaard Floer theory, we finally conclude that these two versions of Floer theories must have the same graded Euler characteristic. 

In the above argument, it is not necessary to treat instanton theory and Heegaard Floer theory separately. Instead, we only use some formal properties that are shared by both theories, and hence we can deal with them at the same time. This observation can be made more general. In Kronheimer and Mrowka's definition of sutured (monopole or instanton) Floer homology, they constructed a closed $3$-manifold, called a closure, out of a balanced sutured manifold in a topological manner, and defined the Floer homology for a balanced sutured manifold to be some direct summands of the Floer homology of its closure. Then they used the formal properties of monopole theory and instanton theory to show that the construction is independent of the choices of the closures. In the following series of work \cite{baldwin2015naturality,baldwin2016contact,baldwin2016instanton,li2019direct,li2019decomposition}, most arguments were also carried out based on topological constructions and hence only depend on the formal properties of Floer theories. 

In this paper, we summarize the necessary properties of Floer theory that are used to build a sutured homology for balanced sutured manifolds. $(3+1)$-TQFTs are functors from cobordism categories to categories of vector spaces. In Subsection \ref{subsec: axioms for SH}, we proposed three extra axioms for a $(3+1)$-TQFT called
\begin{itemize}
    \item (A1) the adjunction inequality axiom;
    \item (A2) the surgery exact triangle axiom;
    \item (A3) the canonical $\intg_2$ (mod 2) grading axiom.
\end{itemize}

A $(3+1)$-TQFT satisfying these axioms is called a \textbf{Floer-type} theory. For any Floer-type theory $\hg$ and any balanced sutured manifold $(M,\ga)$, we construct a vector space $\shg(M,\ga)$, called the \textbf{formal sutured homology} of $(M,\ga)$, over a suitable coefficient field. More precisely, we have the following theorem.

\bthm\label{thm: formal sutured Floer theory}
Suppose $\hg$ is a $(3+1)$-TQFT and suppose $(M,\ga)$ is a balanced sutured manifold. If $\hg$ satisfies Axioms (A1) and (A2), then there is a vector space $\shg(M,\ga)$ well-defined up to multiplication by a unit in the coefficient field $\mathbb{F}$. Suppose $S_1,\dots,S_n$ are properly embedded admissible surfaces inside $(M,\ga)$. Then there exists a $\intg^n$-grading on $\shg(M,\ga)$ induced by these surfaces, \textit{i.e.}
\begin{equation}\label{eq: zn grading}
    \shg(M,\ga)=\bigoplus_{(i_1,\dots,i_n)\in \mathbb{Z}^n}\shg(M,\ga,(S_1,\dots,S_n),(i_1,\dots,i_n)).
\end{equation}

Furthermore, if $\hg$ satisfies Axiom (A3), then there exists a relative $\mathbb{Z}_2$-grading $\shg(M,\ga)$, respecting the decomposition in (\ref{eq: zn grading}). Moreover, the graded Euler characteristic $\chi_{\rm gr}(\shg(M,\ga))$, defined similarly to (\ref{eq: defn of chi}) and determined up to multiplication by a unit in $\pm H_1(M)/{\rm Tors}$, is independent of the choice of the Floer-type theory.
\ethm
\brem
\textit{A priori}, the definition of formal sutured homology depends on a large and fixed integer $g$, which is the genus of the closure; see the Convention after Definition \ref{defn: conv}.
\erem
\brem
The construction of $\shg$ is essentially due to the work of Kronheimer and Mrowka \cite{kronheimer2010knots}. Note that instanton theory, monopole theory and Heegaard Floer theory all satisfy Axioms (A1), (A2), and (A3) with coefficients $\mathbb{C}$, $\mathbb{F}_2$ and $\mathbb{F}_2$, respectively, up to mild modifications (\textit{c.f.} Subsection \ref{subsec: axioms for SH}). Moreover, Axioms (A1), (A2), and (A3) are not limited by the scope of gauge-theoretic theories mentioned above and may hold for other more general $(3+1)$-TQFTs. 
\erem

There is one further step to prove Theorem \ref{thm: main} from Theorem \ref{thm: formal sutured Floer theory}. For Heegaard Floer theory, the construction coming from Theorem \ref{thm: formal sutured Floer theory} is different from the original version of sutured (Heegaard) Floer homology defined by Juh\'asz \cite{juhasz2006holomorphic}. It has been shown by Lekili \cite{lekili2013heegaard} and Baldwin and Sivek \cite{Baldwin2020} that these two constructions coincide with each other. Although not shown explicitly, their proofs also imply that the isomorphism between these two constructions respects gradings. Based on their work, we show the following proposition.

\bprop\label{prop: SFH=SHF}
Suppose $(M,\ga)$ is a balanced sutured manifold and suppose $H=H_1(M)/{\rm Tors}$. Suppose $\shf$ is the sutured homology for balanced sutured manifolds constructed in Theorem \ref{thm: formal sutured Floer theory} for Heegaard Floer theory. Then we have$$\chi_{\rm gr}(\shf(M,\ga))=\chi_{\rm gr}(SFH(M,\ga))\in \mathbb{Z}[H]/\pm H.$$
\eprop

Next, we discuss the $\intg_2$-grading on $\shg(M,\ga)$. Following \cite{kronheimer2010knots}, to construct $\shg(M,\ga)$, we first construct a closure $Y$ from $(M,\ga)$. From a fixed balanced sutured manifold $(M,\ga)$, we can construct infinitely many different closures (with the same genus), and the Floer homology of each closure has its own (absolute) $\intg_2$-grading. Although we can construct isomorphisms between the Floer homology of different closures, the maps do not necessarily respect the $\intg_2$-gradings. See \cite[Section 2.6]{kronheimer2010instanton} for a concrete example. Thus, we cannot obtain a canonical $\intg_2$-grading on $\shg(M,\ga)$ and the Euler characteristic can only be defined up to a sign (since we do know the maps between closures are homogenous with respect to the $\intg_2$-gradings). 

However, if we focus on balanced sutured manifolds whose underlying $3$-manifolds are knot complements of null-homologous knots and whose sutures have two components, it is possible to obtain a canonical $\intg_2$-grading. The idea is to compare closures of a general knot with closures of the unknot in $S^3$ and then fix the relative $\mathbb{Z}_2$-grading. When the suture on the boundary of the knot complement consists of two meridians, we recover the canonical $\intg_2$-grading on $KHI(Y,K)$ already known by Floer \cite{floer1990knot} and Kronheimer and Mrowka \cite{kronheimer2010instanton}. When the suture on the boundary of knot complement consists of two curves of slope $-n$, this canonical $\intg_2$-grading is also carried over to the following decomposition of $I^\sharp (Y)$ introduced by the authors.
\bthm[{\cite[Theorem 1.12]{LY2020}}]\label{thm: torsion spin c decomposition}
Suppose $Y$ is a closed 3-manifold, and $K\subset Y$ is a null-homologous knot. Suppose $\widehat{Y}$ is obtained from $Y$ by performing the $q/p$ surgery along $K$ with $q>0$. Then there is a decomposition up to isomorphism
$$I^{\sharp}(\widehat{Y})\cong \bigoplus_{i=0}^{q-1}I^{\sharp}(\widehat{Y},i),$$
associated to the knot $K$ and the slope $q/p$.
\ethm
\bprop\label{prop: Euler Char of I sharp}
Under the hypothesis and the statement of Theorem \ref{thm: torsion spin c decomposition}, there is a well-defined $\intg_2$-grading on $I^{\sharp}(\widehat{Y},i)$. Moreover, for $i=0,\dots,q-1$, we have
$$\chi(I^{\sharp}(\widehat{Y},i))=\chi(I^{\sharp}(Y)).$$
\eprop
\bcor\label{cor: L space surgery}
Suppose $K$ is a knot in an integral homology sphere $Y$. Suppose further that $r=q/p$ is a rational number with $q>0$. Then, the $3$-manifold $Y_{r}(K)$ is an instanton L-space (\textit{i.e.}, $\dim_\mathbb{C}I^\sharp(Y)=|H_1(Y;\mathbb{Z})|$) if and only if for $i=0,\dots,q-1$, we have
$$I^{\sharp}(Y_r(K),i)\cong\mathbb{C}.$$
\ecor
\begin{proof}
If for $i=0,\dots,q-1$, we have
$I^{\sharp}(Y_r(K),i))\cong \mathbb{C}$,
then it follows directly from Theorem \ref{thm: torsion spin c decomposition} that $Y_r(K)$ is an instanton L-space.

Now suppose $Y_r(K)$ is an instanton L-space. Applying Proposition \ref{prop: Euler Char of I sharp} to $Y$, we have \begin{equation}\label{S3 char}
    \dim_\mathbb{C}I^\sharp(Y_r(K))\ge|\chi(I^\sharp(Y_r(K))|=|\sum_{i=0}^{q-1}\chi(I^{\sharp}(Y_r(K),i))|=|\sum_{i=0}^{q-1}\chi(I^{\sharp}(Y))|=q.
\end{equation}By assumption, $Y_r(K)$ is an instanton $L$-space, \textit{i.e.}, $$\dim_\mathbb{C}I^\sharp(Y_r(K))=|H_1(Y_r(K))|=q.$$
Hence the inequality in (\ref{S3 char}) is sharp, which implies $\dim_\mathbb{C}I^{\sharp}(Y_r(K),i)=1$.
\end{proof}

The techniques to prove Proposition \ref{prop: Euler Char of I sharp} can also be applied to study the minus version of instanton knot homology $\khi$, which was introduced by the first author in \cite[Section 5]{li2019direct}.
\bprop\label{prop: Euler char of khi-minus}
Suppose $K\subset S^3$ is a knot and $\Delta_K(t)$ is the symmetrized Alexander polynomial of $K$. Then there is a canonical $\intg_2$-grading on $\khi(-S^3,K)$. Furthermore, we have
$$\sum_{i\in\intg}\chi(\khi(-S^3,K,i))\cdot t^i=-\Delta_K(t)\cdot \sum_{i=0}^{+\infty}t^{-i}.$$
\eprop
\brem
The analogous result of Proposition \ref{prop: Euler char of khi-minus} in Heegaard Floer theory had been known by the work of Ozsv\'ath and Szab\'o \cite{ozsvath2004holomorphicknot}.
\erem


\begin{convs}If it is not mentioned, homology groups and cohomology groups are with $\mathbb{Z}$ coefficients, \textit{i.e.}, we write $H_*(Y)$ for $H_*(Y;\mathbb{Z})$. A general field is denoted by $\mathbb{F}$, and the field with two elements is denoted by $\mathbb{F}_2$.

If it is not mentioned, all manifolds are smooth and oriented. Moreover, all manifolds are connected unless we indicate disconnected manifolds are also considered. This usually happens when discussing cobordism maps from a $(3+1)$-TQFT.

Suppose $M$ is an oriented manifold. Let $-M$ denote the same manifold with the reverse orientation, called the \textbf{mirror manifold} of $M$. If it is not mentioned, we do not consider orientations of knots. Suppose $K$ is a knot in a 3-manifold $M$. Then $(-M,K)$ is the \textbf{mirror knot} in the mirror manifold.

For a manifold $M$, let ${\rm int} (M)$ denote its interior. For a submanifold $A$ in a manifold $Y$, let $N(A)$ denote the tubular neighborhood. The knot complement of $K$ in $Y$ is denoted by $Y(K)=Y\backslash {\rm int} (N(K))$.

For a simple closed curve on a surface, we do not distinguish between its homology class and itself. The algebraic intersection number of two curves $\al$ and $\be$ on a surface is denoted by $\al\cdot\be$, while the number of intersection points between $\al$ and $\be$ is denoted by $|\al\cap \be|$. A basis $(m,l)$ of $H_1(T^2;\mathbb{Z})$ satisfies $m\cdot l=-1$. The \textbf{surgery} means the Dehn surgery and the slope $q/p$ in the basis $(m,l)$ corresponds to the curve $qm+pl$.

A knot $K\subset Y$ is called \textbf{null-homologous} if it represents the trivial homology class in $H_1(Y;\mathbb{Z})$, while it is called \textbf{rationally null-homologous} if it represents the trivial homology class in $H_1(Y;\mathbb{Q})$. We write $\mathbb{Z}_n$ for $\mathbb{Z}/n\mathbb{Z}$.


An argument holds for \textbf{large enough} or \textbf{sufficiently large} $n$ if there exists a fixed $N\in\mathbb{Z}$ so that the argument holds for any integer $n>N$. 
\end{convs}
\begin{org}
The paper is organized as follows. In Section \ref{sec: Axioms of sutured Floer homology}, we introduce three axioms to define formal sutured homology for balanced sutured manifolds and prove the first part of Theorem \ref{thm: formal sutured Floer theory}. Moreover, we state many useful properties for the proof of the second part of Theorem \ref{thm: formal sutured Floer theory}. In Section \ref{sec: hf},  we discuss the modification of Heegaard Floer theory to make it suitable to formal sutured homology and prove Proposition \ref{prop: SFH=SHF}. In Section \ref{sec: Equivalence of graded Euler characteristics}, we prove the second part of Theorem \ref{thm: formal sutured Floer theory}. In Section \ref{sec: mod 2}, we construct a canonical $\mathbb{Z}_2$-grading for balanced sutured manifolds obtained from knots in closed 3-manifolds and prove Proposition \ref{prop: Euler Char of I sharp} and Proposition \ref{prop: Euler char of khi-minus}.

\end{org}
\begin{ack}
 The authors thank Ian Zemke for pointing out a mistake in the previous statement of Corollary \ref{cor: free iso} and helping with the new proof of Floer's excision theorem for Heegaard Floer theory in Subsection \ref{subsec: Floer}. The authors also thank John A. Baldwin and Yi Xie for valuable discussions and anonymous referees for the helpful comments. The first author thanks Ciprian Manolescu and Joshua Wang for helpful conversations. The second author thanks his supervisor Jacob Rasmussen for patient guidance and helpful comments and thanks his parents for support and constant encouragement. The second author is also grateful to Yi Liu for inviting him to BICMR, Peking University, Linsheng Wang for valuable discussions on algebra, and Honghuai Fang, Zekun Chen, Fei Chen, and Shengyu Zou for the company when he was writing this paper.
\end{ack}

\section{Axioms and formal properties for sutured homology}
\label{sec: Axioms of sutured Floer homology}

In this section, we construct formal sutured homology and prove some basic properties.

\subsection{Axioms of a Floer-type theory for closed 3-manifolds}\label{subsec: axioms for SH}\quad

Let ${\bf Cob}^{3+1}$ be the cobordism category whose objects are closed oriented (possibly disconnected) 3-manifolds, and whose morphisms are oriented (possibly disconnected) 4-dimensional cobordisms between closed oriented 3-manifolds. The disjoint union gives a monoidal structure on ${\bf Cob}^{3+1}$ and reversing the orientation induces the dual of an object. Let ${\bf Vect}_{\mathbb{F}}$ be the category of $\mathbb{F}$-vector spaces, where $\mathbb{F}$ is a suitably chosen coefficient field. The monoidal and dual structure on ${\bf Vect}_{\mathbb{F}}$ is induced by tensor product and dual space. 

A \textbf{$(3+1)$ dimensional topological quantum field theory}, or in short \textbf{$(3+1)$-TQFT}, is a symmetric monoidal functor
$${\mathbf H}:{\bf Cob}^{3+1}\ra{\bf Vect}_{\mathbb{F}}$$preserving the dual structure. For a closed oriented 3-manifold $Y$, we write $\hg(Y)$ for the related vector space, called the \textbf{$\hg$-homology} of $Y$. For an oriented cobordism $W$, we write $\hg(W)$ for the induced map between $\hg$-homologies of boundaries, called the \textbf{$\hg$-cobordism map} associated to $W$. If $\hg$ is fixed, then we simply write \textbf{homology} and \textbf{cobordism map} for $\hg$-homology and $\hg$-cobordism map, respectively. Note that by the definition of the involved categories, we have $$\hg(Y_1\sqcup Y_2)=\hg(Y_1)\otimes_{\mathbb{F}}\hg(Y_2)\aand \hg(-Y)\cong {\rm Hom}_\mathbb{F}(\hg(Y),\mathbb{F}).$$

It is well-known that Floer theories are special cases of $(3+1)$-TQFTs. Summarized from known Floer theories, we propose the following definition.
\bdefn\label{defn: floer type theory}

A $(3+1)$-TQFT $\hg$ is called a \textbf{Floer-type theory} if it satisfies the following three Axioms (A1), (A2), and (A3).
\edefn

{\bf A1. Adjunction inequality}. For a closed oriented 3-manifold $Y$ and a second homology class $\al\in H_2(Y)$, there is a $\intg$-grading of $\hg(Y)$ associated to $\al$, \textit{i.e.}, we have
$$\hg(Y)=\bigoplus_{i\in\intg}\hg(Y,\al,i).$$
This grading satisfies the following properties.

\quad {\bf A1-1}. For any odd integer $i$, we have $\hg(Y,\al,i)=0$.

\quad {\bf A1-2}. For $i\in \intg\backslash\{0\}$, the summand $\hg(Y,\al,i)$ is a finite dimensional vector space over $\mathbb{F}$. 

\quad {\bf A1-3}. For $i\in\intg$, we have $\hg(Y,\al,i)\cong \hg(Y,\al,-i).$

\quad {\bf A1-4} (Adjunction inequality). Suppose $\Sigma$ is a connected closed oriented surface embedded in $Y$ with $g(\Sigma)\ge 1$. For $|i|>2g(\Sigma)-2$, we have $\hg(Y,[\Sigma],i)=0.$

\quad {\bf A1-5}. Suppose $\Sigma_g$ is a connected closed oriented surface of genus $g(\Sigma)\ge 2$. Suppose $Y=S^1\times \Sigma_g$ and $\Sigma=\{1\}\times \Sigma_g$. Then we have
$$\hg(Y,[\Sigma],2g(\Sigma)-2)\cong\mathbb{F}.$$


\quad {\bf A1-6}. The gradings coming from multiple homology classes are compatible with each other, \textit{i.e.}, if we have $\al_1,\dots,\al_n\in H_2(Y)$, then there is a $\intg^n$-grading on $\hg(Y)$, denoted by $$\hg(Y)=\bigoplus_{(i_1,\dots,i_n)\in\mathbb{Z}^n}\hg(Y,(\al_1,\dots,\al_n),(i_1,\dots,i_n)).$$Moreover, we have
\beq
&\hg(Y,\al_1+\dots+\al_n,i)=
\bigoplus_{\substack{(i_1,\dots,i_n)\in\mathbb{Z}^n\\i_1+\dots+i_n=i}}\hg(Y,(\al_1,\dots,\al_n),(i_1,\dots,i_n)).
\eeq

\quad {\bf A1-7}. Suppose $W$ is an oriented cobordism from $Y_1$ to $Y_2$. Suppose $\al_1,\dots,\al_n\in H_2(Y_1)$ and $\be_1,\dots,\be_n\in H_2(Y_2)$ are homology classes such that for $i=1,\dots,n$, we have
$$\al_i=\be_i\in H_2(W).$$
Then the cobordism map $\hg(W)$ respects the grading associated to those homology classes:
$$\hg(W): \hg(Y_1,(\al_1,\dots,\al_n),(i_1,\dots,i_n))\ra \hg(Y_1,(\be_1,\dots,\be_n),(i_1,\dots,i_n)).$$ 

{\bf A2. Surgery exact triangle}. Suppose $M$ is a connected compact oriented 3-manifold with toroidal boundary. Let $\ga_1$, $\ga_2$, $\ga_3$ be three connected oriented simple closed curves on $\partial M$ such that
$$\ga_1\cdot \ga_2=\ga_2\cdot\ga_3=\ga_3\cdot\ga_1=-1.$$
Let $Y_1$, $Y_2$, and $Y_3$ be the Dehn fillings of $M$ along curves $\ga_1$, $\ga_2$, and $\ga_3$, respectively. Then there is an exact triangle
\begin{equation}\label{eq_2: Floer's triangle}
\xymatrix{
\hg(Y_1)\ar[rr]&&\hg(Y_2)\ar[dl]\\
&\hg(Y_3)\ar[ul]&
}	
\end{equation}
Moreover, maps in the above triangle are induced by the natural cobordisms associated to different Dehn fillings.

\brem
It is worth mentioning that Axioms (A1) and (A2) are enough for defining formal sutured homology for balanced sutured manifolds. The following Axiom (A3) is only involved when considering Euler characteristics.
\erem

{\bf A3. $\intg_2$-grading}. For any closed oriented 3-manifold $Y$, there is a canonical $\intg_2$-grading on $\hg(Y)$, denoted by
$$\hg(Y)=\hg_0(Y)\oplus \hg_1(Y).$$
This grading satisfies the following properties.

\quad{\bf A3-1}. The $\intg_2$-grading is compatible with the grading in Axiom (A1). More precisely, if we have $\al_1,\dots,\al_n\in H_2(Y)$, then there is a $\intg_2\oplus\intg^n$-grading on $\hg(Y)$:
$$\hg(Y)=\bigoplus_{j\in\{0,1\}}\bigoplus_{(i_1,\dots,i_n)\in\intg^n}\hg_j(Y,(\al_1,\dots,\al_n),(i_1,\dots,i_n)).$$

\quad {\bf A3-2}. Suppose $\Sigma_g$ is a connected closed oriented surface of genus $g\geq 2$. Suppose $Y=S^1\times \Sigma_g$ and $\Sigma=\{1\}\times \Sigma_g$. Then we have
$$\hg(Y,[\Sigma],2g-2)=\hg_1(Y,[\Sigma],2g-2)\cong\mathbb{F}.$$

\quad {\bf A3-3}. Suppose $W$ is a cobordism from $Y_1$ to $Y_2$. Then $\hg(W)$ is homogeneous with respect to the canonical $\intg_2$-grading. Its degree can be calculated by the following degree formula
\begin{equation}\label{eq: degree formula}
    {\rm deg}(\hg(W))\equiv \frac{1}{2}(\chi(W)+\sigma(W)+b_1(Y_2)-b_1(Y_1)+b_0(Y_2)-b_0(Y_1)) \pmod 2.
\end{equation}

\brem
The canonical $\mathbb{Z}_2$-grading is essentially determined by Axioms (A3-2) and (A3-3) (\textit{c.f.} \cite[Section 25.4]{kronheimer2007monopoles}). The normalization of the $\mathbb{Z}_2$-grading for the generator of $\hg(Y,[\Sigma],2g-2)$ is not essential. Assuming $\hg(Y,[\Sigma],2g-2)=\hg_0(Y,[\Sigma],2g-2)$ shifts the canonical $\mathbb{Z}_2$-grading for all 3-manifolds. It is worth mentioning that two existing discussions on this $\intg_2$-grading in \cite{lidman2020framed,kronheimer2011knot},  adopted different normalizations.
\erem
The degrees of the maps in Axiom (A2) were described explicitly by Kronheimer and Mrowka \cite[Section 42.3]{kronheimer2007monopoles}. For the convenience of later usage, we present the discussion here.

\bprop[{\cite[Section 42.3]{kronheimer2007monopoles}}]\label{prop_2: parity of maps in an exact triangle}
Suppose $\delta$ is a unit in ${\rm ker}(i_*)$ for the map $$i_*:H_1(\partial M;\mathbb{Q})\ra H_1(M;\mathbb{Q}).$$In the surgery exact triangle (\ref{eq_2: Floer's triangle}), we can determine the parities of the maps $f_1$, $f_2$, and $f_3$ as follows. 
\begin{enumerate}[(1)]
    \item If there is an $i=1,2,3$ so that $\ga_i\cdot \delta=0$, then $f_{i-1}$ is odd and the other two are even. 
    \item If $\ga_i\cdot\delta\neq0$ for any $i=1,2,3$, then there is a unique $j\in\{1,2,3\}$ so that $\ga_j\cdot \delta$ and $\ga_{j+1}\cdot \delta$ are of the same sign. Then the map $f_j$ is odd and the other two are even.
\end{enumerate}
Here the indices are taken mod 3.
\eprop

With Proposition \ref{prop_2: parity of maps in an exact triangle}, the following lemma is straightforward.
\blem\label{lem_2: Floer's exact triangle after grading shift}
In the surgery exact triangle (\ref{eq_2: Floer's triangle}), after arbitrary shifts on the canonical $\intg_2$-gradings on $\hg(Y_i)$ for all $i=1,2,3$, exactly one of the following two cases happens. 
\begin{enumerate}[(1)]
    \item If all three maps $f_i$ are odd, then we have $$\chi(\hg(Y_1))+\chi(\hg(Y_2))+\chi(\hg(Y_3))=0.$$
   \item If there exists $i=1,2,3$ so that $f_i$ is odd and the other two are even, then
   $$\chi(\hg(Y_{i-1}))=\chi(\hg(Y_i))+\chi(\hg(Y_{i+1})).$$
Here the indices are taken mod 3.
\end{enumerate}
\elem

\brem
If there are no shifts, then case (2) in Lemma \ref{lem_2: Floer's exact triangle after grading shift} happens due to Proposition \ref{prop_2: parity of maps in an exact triangle}.
\erem

In this paper, we discuss three Floer theories, namely, instanton (Donaldson-Floer) theory, monopole (Seiberg-Witten) theory, and Heegaard Floer theory. However, for any of these theories, we need some modifications as follows. Suppose $Y$ is an object of ${\bf Cob}^{3+1}$ and $W$ is a morphism of ${\bf Cob}^{3+1}$.

{\bf Instanton theory}. We consider the decorated cobordism category ${\bf Cob}_\omega^{3+1}$ rather than ${\bf Cob}^{3+1}$. The objects are admissible pairs $(Y,\omega)$, where $\omega\subset Y$ is a 1-cycle such that any component of $Y$ contains at least one component of $\omega$. The admissible condition means that for any component $Y_0$ of $Y$, there exists a closed oriented surface $\Sigma\subset Y_0$ such that $g(\Sigma)\ge 1$ and the algebraic intersection number $\omega\cdot \Sigma$ is odd. Morphisms are pairs $(W,\nu)$, where $\nu$ is a 2-cycle restricting to the given 1-cycles on $\partial W$. The $\hg$-homology and the $\hg$-cobordism map are denoted by $I^\omega(Y)$ and $I(W,\nu)$ (\textit{c.f.} \cite{floer1990knot}), respectively. 

It is worth mentioning that for an admissible pair $(Y,\omega)$, we need $[\omega]\neq 0\in  H^1(Y;\intg)$ as a necessary condition (see \cite{Floer1990}), so $b_1(Y)>0$. Thus, strictly speaking, the objects of ${\bf Cob}_\omega^{3+1}$ do not involve all closed oriented $3$-manifolds.

We choose the coefficient field to be $\mathbb{F}=\mathbb{C}$. The decorations $\omega$ and $\nu$ do not influence Axiom (A1), where the $\mathbb{Z}$-grading is induced by the generalized eigenspaces of $(\mu(\al),\mu({\rm pt}))$ actions for $\al\in H_2(Y)$ (\textit{c.f.} \cite[Section 7]{kronheimer2010knots}). In particular, the Axiom (A1-5) follows from \cite[Proposition 7.4]{kronheimer2010knots}, where we choose $\omega=S^1\times\{{\rm pt}\}$.

In the original statement of the surgery exact triangle in \cite{floer1990knot}, different 3-manifolds in the surgery exact triangle may have different choices of $\omega$. However, by the argument in \cite[Section 2.2]{baldwin2020concordance}, we can assume that, in Axiom (A2), the 1-cycle $\omega$ is unchanged in all manifolds involved in the triangle.

The canonical $\mathbb{Z}_2$-grading for instanton theory was discussed by Kronheimer and Mrowka \cite[Section 2.6]{kronheimer2010instanton}. Indeed, the degree formula (\ref{eq: degree formula}) is from their discussion.

{\bf Monopole theory}. For connected 3-manifolds and cobordisms with connected incoming and outgoing ends, Kronheimer and Mrowka \cite{kronheimer2007monopoles} constructed well-defined vector spaces and linear maps in monopole theory ($\widehat{HM}_\bullet$ or $\widecheck{HM}_\bullet$). However, the cases for disconnect 3-manifolds or cobordisms with disconnected ends are subtle, and hence monopole theory may not satisfy all assumptions of $(3+1)$-TQFT. Fortunately, when restricting to nontorsion spin$^c$ structure, Kronheimer and Mrowka showed that two versions $\widehat{HM}_\bullet$ and $\widecheck{HM}_\bullet$ are canonically identified, and they induce well-defined vector spaces and linear maps that extend to disconnected cases (\textit{c.f.} \cite[Section 2.5-2.6]{kronheimer2007monopoles}), which are denoted by $HM(Y)$ and $HM(W)$ (the spin$^c$ structure is omitted), respectively. In such case, we can still check Axioms for Floer-type theory.

The $\mathbb{Z}$-grading in Axiom (A1) is induced by $\langle c_1(\mathfrak{s}),\al\rangle$ for $\mathfrak{s}\in \spin(Y)$ and $\al\in H_2(Y)$ (\textit{c.f.} \cite[Section 2.4]{kronheimer2010knots}). In particular, the Axiom (A1-5) follows from \cite[Lemma 2.2]{kronheimer2010knots}.

We choose the coefficient field to be $\mathbb{F}=\mathbb{F}_2$. This is because originally the surgery exact triangle is only proved in characteristic two, by work of Kronheimer, Mrowka, Ozsv\'ath, and Szab\'o \cite{kronheimer2007monopolesandlens}. It is worth mentioning that the surgery exact triangle with $\mathbb{Q}$ coefficients was proved by Lin, Ruberman, and Saveliev \cite[Section 4]{lin2018lefschetz}, and the one with $\mathbb{Z}$ coefficients was under working by Freeman \cite{freeman2021surgery}. So we can extend the discussion to $\mathbb{F}=\mathbb{Q}$ or $\mathbb{C}$ for monopole theory. It is also worth mentioning that in \cite{kronheimer2010knots}, when Kronheimer and Mrowka first introduced sutured monopole homology, they worked only with $\intg$ coefficients (and did not use the surgery exact triangle). The case of $\mathbb{F}_2$ coefficients was later verified by Sivek \cite[Section 2.2]{sivek2012monopole}. Sivek's modification can be applied to any field coefficient since all Tor terms vanish for vector spaces.

The canonical $\mathbb{Z}_2$-grading for monopole theory was discussed by Kronheimer and Mrowka \cite[Section 25.4]{kronheimer2007monopoles}. When considering cobordisms of connected 3-manifolds, the degree formula (\ref{eq: degree formula}) is the same as the formula in \cite[Definition 25.4.1]{kronheimer2007monopoles}.

{\bf Heegaard Floer theory}. Similar to monopole theory, the connected cases and disconnected cases need to be considered separately. For connected cases, Ozsv{\'{a}}th and Szab{\'{o}} \cite{ozsvath2004holomorphic,Ozsvath2004c,Ozsvath2006} constructed vector spaces and linear maps in Heegaard Floer theory ($HF^-$ and $HF^+$). However, there is some naturality problem, which was resolved by Juh{\'{a}}sz, Thurston, and Zemke \cite{Juhasz2012} and Zemke \cite{Zemke2020hat,Zemke2019}, at the cost of adding basepoints and paths to 3-manifolds and cobordisms, respectively. In Subsection \ref{subsec:Cobordism maps for restricted graph cobordisms}, we will propose a new transitive system to get rid of the dependence of basepoints and paths. The $\hg$-homology and the $\hg$-cobordism map in such cases are denoted by $HF(Y)$ and $HF(W)$.

Zemke's construction also extend to disconnected cases, at the cost of introducing graphs in cobordisms. The cobordism maps would be different when the cobordisms are the same but graphs are different (for example, see Corollary \ref{cor: free iso}, where $S_{w_0}^+$ and $S_{w_0}^-$ both correspond the product cobordism but the graphs are different).

Another subtlety is the duality. The duality for connected cases was proved in \cite[Theorem 3.5]{Ozsvath2006}, which is on the homology level. However, the duality for disconnected cases (with graphs) was proved only on the chain level \cite{Zemke2018}. Since the chain complexes constructed by Zemke are over $\ft[U]$ (which is not a field), the duality on the chain level does not imply the duality on the homology level (the example about $S_{w_0}^+$ and $S_{w_0}^-$ also indicates this subtlety).

When we construct the formal sutured homology, we only need two kinds of cobordisms, ones from Dehn surgeries and ones from Floer's excision theorem. The first kind is in connected cases and the second kind is possibly in disconnected cases. In Subsection \ref{subsec: Floer}, we construct cobordisms with specific graphs to prove Floer's excision theorem in Heegaard Floer theory. Hence, we can still check Axioms for Floer-type theory and use Heegaard Floer theory to build a formal sutured homology following the strategy in the next subsection.

We choose the coefficient field to be $\mathbb{F}=\mathbb{F}_2$. This is because we have to use the naturality results in \cite{Juhasz2012}, which works only for $\mathbb{F}_2$. 

For characteristic zero, the naturality results for closed 3-manifolds were obtained by Gartner in \cite{2019naturality}. However, the naturality results for cobordisms are still under working. Hence we choose to focus on characteristic two.

Similar to monopole theory, the $\mathbb{Z}$-grading in Axiom (A1) is induced by $\langle c_1(\mathfrak{s}),\al\rangle$ for $\mathfrak{s}\in \spin(Y)$ and $\al\in H_2(Y)$. In particular, Axiom (1-5) follows from \cite[Corollary 17]{lekili2013heegaard}. The surgery exact triangle in Axiom (A2) is proved by Ozsv{\'{a}}th and Szab{\'{o}} \cite{Ozsvath2004c}.

There are many ways to fix the $\mathbb{Z}_2$-grading for Heegaard Floer theory. See \cite[Section 10.4]{Ozsvath2004c} and \cite[Section 2.4]{Friedl2009}. However, we can arrange the canonical $\mathbb{Z}_2$-grading to be the same as those for instanton theory and monopole theory. This is possible because the degree formula (\ref{eq: degree formula}) only depends on algebraic-topological information of cobordisms and 3-manifolds.

\subsection{Formal sutured homology of balanced sutured manifolds}\label{sec: Floer theory for sutured manifold}\quad

In \cite{kronheimer2010knots}, Kronheimer and Mrowka constructed sutured monopole homology $SHM$ and sutured instanton homology $SHI$ by considering closures of balanced sutured manifolds. The discussion and construction in this subsection are based on \cite{kronheimer2010knots,baldwin2015naturality} except for the proof of Proposition \ref{prop: X_P induces an isomorphism}.

\bdefn[{\cite[Definition 2.2]{juhasz2006holomorphic}}]\label{defn_2: balanced sutured manifold}
A \textbf{balanced sutured manifold} $(M,\ga)$ consists of a compact 3-manifold $M$ with non-empty boundary together with a closed 1-submanifold $\ga$ on $\partial{M}$. Let $A(\ga)=[-1,1]\times\ga$ be an annular neighborhood of $\ga\subset \partial{M}$ and let $R(\ga)=\partial{M}\backslash{\rm int}(A(\ga))$, such that they satisfy the following properties.
\begin{enumerate}[(1)]
    \item Neither $M$ nor $R(\ga)$ has a closed component.
    \item If $\partial{A(\ga)}=-\partial{R(\ga)}$ is oriented in the same way as $\ga$, then we require this orientation of $\partial{R(\ga)}$ induces the orientation on $R(\ga)$, which is called the \textbf{canonical orientation}.
    \item Let $R_+(\ga)$ be the part of $R(\ga)$ for which the canonical orientation coincides with the induced orientation on $\partial{M}$ from $M$, and let $R_-(\ga)=R(\ga)\backslash R_+(\ga)$. We require that $\chi(R_+(\ga))=\chi(R_-(\ga))$. If $\ga$ is clear in the context, we simply write $R_\pm=R_\pm(\ga)$, respectively.
\end{enumerate}
\edefn
\bdefn[\cite{kronheimer2010knots}]\label{defn: closure}
Suppose $(M,\ga)$ is a balanced sutured manifold. Let $T$ be a connected compact oriented surface such that the numbers of components of $\partial T$ and $\ga$ are the same. Let the \textbf{preclosure} $\widetilde{M}$ of $(M,\ga)$ be
$$\widetilde{M}\deq M\mathop{\cup}_{\ga=-\partial T}[-1,1]\times T.$$
The boundary of $\widetilde{M}$ consists of two components $$\widetilde{R}_+=R_+(\ga)\cup \{1\}\times T\aand \widetilde{R}_-=R_-(\ga)\cup \{-1\}\times T.$$

Let $h:\widetilde{R}_+\xra{\cong}\widetilde{R}_-$ be a diffeomorphism which reverses the boundary orientations (\textit{i.e.} preserves the canonical orientations). Let $Y$ be the 3-manifold obtained from $\widetilde{M}$ by gluing $\widetilde{R}_+$ to $\widetilde{R}_-$ by $h$ and let $R$ be the image of $\widetilde{R}_+$ and $\widetilde{R}_-$ in $Y$. The pair $(Y,R)$ is called a {\bf closure} of $(M,\ga)$. The genus of $R$ is called the {\bf genus} of the closure $(Y,R)$. For a closure $(Y,R)$ with $g(R)\ge 2$ and a $(3+1)$-TQFT $\hg$ satisfying Axiom (A1), define
$$\hg(Y|R)\deq \hg(Y,[R],2g(R)-2).$$

\edefn
\brem
For instanton theory, we also choose a point $p$ on $T$ and choose a diffeomorphism $h$ such that $h(\{1\}\times p)=\{-1\}\times p$. The image of $[-1,1]\times p$ in $Y$ becomes a 1-cycle $\omega$ and we have $\omega\cdot R=1$. We use $(Y,R,\omega)$ for the definition of a closure. We do not mention this subtlety later since everything works well under this modification.
\erem

Suppose $(Y_1,R_1)$ and $(Y_2,R_2)$ are two closures of $(M,\ga)$ of the same genus. We now construct a canonical map
$$\Phi_{12}:\hg(Y_1|R_1)\ra \hg(Y_2|R_2).$$
Note that $Y_2$ can be obtained from $Y_1$ as follows. There exists an orientation preserving diffeomorphism $h_{12}:R_1\ra R_1$ so that if we cut $Y_1$ open along $R_1$ and reglue using $h_{12}$, then we obtain a new 3-manifold $Y_1^\p$ together with the surface $R_1^\p\subset Y_1^\p$. Furthermore, there exists a diffeomorphism $\phi: Y_1^\p\ra Y_2$ such that
$$\phi|_{M}={\rm id}_M~{\rm and~}\phi(R_1^\p)=R_2.$$
Let $X_\phi$ be a cobordism from $Y_1^\p$ to $Y_2$ induced by $\phi$. It is straightforward to check $$\hg (X_\phi):\hg(Y_1^\p|R_1^\p)\ra \hg(Y_2|R_2)$$is an isomorphism. We can regard $h_{12}$ as a composition of Dehn twists along curves on $R_1$:
$$h_{12}=D_{\al_1}^{e_1}\circ\cdots\circ D_{\al_n}^{e_n}.$$
Here $e_i\in\{\pm1\}$, where $e_1=1$ means a positive Dehn twist, and $e_1=-1$ means a negative Dehn twist. Suppose
$$N=\{i~|~e_i=-1\}~{\rm~and~}P=\{i~|~e_i=1\}.$$

Note that the resulting 3-manifold of cutting $Y_1$ open along $R_1$ and regluing by $D_{\al_i}^{e_i}$ is the same as the resulting 3-manifold of performing a $(-e_i)$-surgery along $\al_i\subset R_1\subset Y_1$. We take a neighborhood $N(R_1)$ of $R_1\subset Y_1$, and choose an identification $N(R_1)=[-1,1]\times R_1.$
Pick $$-1<t_1<\cdots<t_n<1$$so that $t_i\neq0$ for $i=1,\dots,n$, and isotope $\al_i$ to the level $\{t_i\}\times R_1\subset N(R_1)\subset Y_1$. Let $Y_P$ be the 3-manifold obtained from $Y_1$ by performing $(-1)$-surgeries along $\al_i$ for all $i\in P$. There is a natural cobordism $X_P$ from $Y_1$ to $Y_P$ by attaching framed $4$-dimensional 2-handles to the product $[0,1]\times Y_1$ along $\al_i\times\{1\}$. Furthermore, the manifold $Y_P$ can also be obtained from $Y_1^\p$ by performing $(-1)$-surgeries along $\al_i$ for all $i\in N$. Hence there is a similar cobordism $X_N$ from $Y_1^\p$ to $Y_P$. Since $t_i\neq 0$, the surface $R_1=\{0\}\times R_1$ survives in all surgeries. Let $R_P\subset Y_P$ be the corresponding surface.
\bdefn[\cite{baldwin2015naturality}]\label{defn: canonical maps}
Define
$$\Phi_{12}=\hg(X_\phi)\circ \hg(X_N)^{-1}\circ \hg(X_P):\hg(Y_1|R_1)\ra \hg(Y_2|R_2).$$
\edefn
\bprop\label{prop: X_P induces an isomorphism}
The maps
$$\hg(X_P):\hg(Y_1|R_1)\ra \hg(Y_P|R_N)~{\rm and}~ \hg(X_N):\hg(Y_1'|R_1')\ra \hg(Y_P|R_P)$$
are both isomorphisms.
\eprop

\brem
Proposition \ref{prop: X_P induces an isomorphism} restates \cite[Lemma 4.9]{baldwin2015naturality}. However, the proof in that paper involves a non-vanishing result for minimal Lefschetz fibrations. See \cite[Proposition B.1]{baldwin2015naturality}. Yet this non-vanishing result is not covered by Axioms (A1), (A2), and (A3), so we present an alternative proof of Proposition \ref{prop: X_P induces an isomorphism} based on surgery exact triangles from Axiom (A2). Also, it is worth mentioning that Baldwin and Sivek worked with $\intg$ coefficients for monopole theory in \cite{baldwin2015naturality}, while we work with $\intg_2$ coefficients. The choice of coefficients matters since the existing proof of the surgery exact triangle in monopole theory is only carried out in 
$\mathbb{F}_2$ and $\mathbb{Q}$.
\erem

\bpf[Proof of Proposition \ref{prop: X_P induces an isomorphism}]
The cobordisms $X_P$ and $X_N$ are constructed similarly, so we only prove $X_P$ is an isomorphism. Furthermore, we can assume that $P$ has only one element $\al_1$. If it has more elements, then $X_P$ is simply the composition of cobordisms associated to single Dehn surgeries. With this assumption, the manifold $Y_P$ is obtained from $Y_1$ by performing a $(-1)$-surgery along $\al_1$. Let $Y_0$ be obtained from $Y_1$ by performing a 0-surgery along $\al_1$, and $R_1$ survives to become $R_0\subset Y_0$. Then we have an exact triangle by Axioms (A1-7) and (A2):
\begin{equation}\label{eq: exact triangle, appendix}
	\xymatrix{
	\hg(Y_1|R_1)\ar[rr]^{\hg(X_P)}&&\hg(Y_P|R_P)\ar[dl]\\
	&\hg(Y_0|R_0)\ar[ul]&
	}
\end{equation}
To show that $\hg(X_P)$ is an isomorphism, it suffices to show that $\hg(Y_0|R_0)=0$. Indeed, since $Y_0$ is obtained from a 0-surgery along $\al_0$, and $\al_0$ can be isotoped to be a simple closed curve on $R_1$, the surface $R_0\subset Y_0$ is compressible. Hence $\hg(Y_0|R_0)=0$ by the adjunction inequality in Axiom (A1-4).
\epf

With Proposition \ref{prop: X_P induces an isomorphism} settled down, the rest of the argument in \cite{baldwin2015naturality} can be applied to our setup verbatim, and we have the following theorem.

\bthm[\cite{baldwin2015naturality}]\label{thm: canonical maps}
Suppose $(M,\ga)$ is a balanced sutured manifold and $(Y_1,R_1)$ and $(Y_2,R_2)$ are two closures of the same genus. Then the isomorphism
$$\Phi_{12}:\hg(Y_1|R_1)\ra \hg(Y_2|R_2)$$
defined in Definition \ref{defn: canonical maps} satisfies the following properties.
\begin{enumerate}
\item The map $\Phi_{12}$ is well-defined up to multiplication by a unit in $\mathbb{F}$.
\item If $(Y_1,R_1)=(Y_2,R_2)$, then
$$\Phi_{12}\doteq {\rm id},$$
where $\doteq$ means the equation holds up to multiplication by a unit in $\mathbb{F}$.
\item If there is a third closure $(Y_3,R_3)$ of the same genus, then we have
$$\Phi_{13}\doteq\Phi_{23}\circ\Phi_{12}.$$
\end{enumerate}
\ethm

\brem\label{rem: grading is well-defined}
In Baldwin and Sivek's original work, the requirement that the two closures have the same genus could be dropped, at the cost of involving local coefficient systems. However, the naturality of Heegaard Floer theory only holds when we restrict to one spin$^c$ structure at a time and consider the transition map projectively (\textit{c.f.} \cite[Remark 12.1]{Zemke2018}, see also \cite[Section 3]{juhasz1804twisted}), which makes arguments more subtle. Since it is enough to work with closures of a large and fixed closure in the current paper, we choose not to discuss the local coefficients.
\erem
\bdefn[{\cite{Juhasz2012,baldwin2015naturality}}]\label{defn: transitive system}
A \textbf{projectively transitive system} of vector spaces over a field $\mathbb{F}$ consists of
\benu
    \item a set $A$ and collection of vector spaces $\{V_\al\}_{\al\in A}$ over $\mathbb{F}$,
    \item a collection of linear maps $\{g_{\be}^\al\}_{\al,\be\in A}$ well-defined up to multiplication by a unit in $\mathbb{F}$ such that
    \begin{enumerate}
        \item $g_\be^\al$ is an isomorphism from $V_\al$ to $V_\be$ for any $\al,\be\in A$, called a \textbf{canonical map},
    \item $g_\al^\al\doteq {\rm id}_{V_\al}$ for any $\al\in A$,
    \item $g^\be_\ga\circ g^\al_\be\doteq g^\al_\ga$ for any $\al,\be,\ga\in A$.
    \end{enumerate}
\eenu
A morphism of projectively transitive systems of vector spaces over a field $\mathbb{F}$ from $(A, \{V_\al\}, \{g_\be^\al\})$ to $(B, \{U_\ga\}, \{h^\ga_\delta\})$
is a collection of maps $\{f_\ga^\al\}_{\al\in A,\ga\in B}$ such that
\benu
\item $f_\ga^\al$ is a linear map from $V_\al$ to $U_\ga$ well-defined up to multiplication by a unit in $\mathbb{F}$ for any $\al\in A$ and $\ga\in B$,
\item $f^\be_\delta\circ g^\al_\be\doteq h^\ga_\delta\circ f^\al_\ga$ for any $\al,\be\in A$ and $\ga,\delta\in B$.
\eenu

A \textbf{transitive system} of vector spaces over a field $\mathbb{F}$ if it is a projectively transitive system and all equations with $\doteq$ are replaced by ones with $=$. A morphism of transitive systems of vector spaces over a field $\mathbb{F}$ is defined similarly.

We can replace vector spaces with groups or chain complexes of vector spaces and define the projectively transitive system and the transitive system similarly.
\edefn
\brem
A transitive system of vector spaces $(A, \{V_\al\}, \{g_\be^\al\})$ over a field $\mathbb{F}$ canonically defines an actual vector space over $\mathbb{F}$
$$V\deq \coprod _{\al\in A} V_\al/\sim,$$
where $v_\al\sim v_\be$ if and only if $g^\al_\be(v_\al)=v_\be$ for any $v_\al\in V_\al$ and $v_\be\in V_\be$. A morphism of transitive systems of vector spaces canonically defines an linear map between corresponding actual vector spaces.
\erem
\begin{conv}
If $\mathbb{F}=\mathbb{F}_2$, a projectively transitive system over $\mathbb{F}$ is simply a transitive system since $\mathbb{F}_2$ has only one unit. In this case, we do not distinguish the projectively transitive system, the transitive system and the corresponding actual vector space. For a general field $\mathbb{F}$, the morphisms between projectively transitive systems are also called maps.
\end{conv}
\bdefn\label{defn: conv}
Suppose $\hg$ is a $(3+1)$-TQFT satisfying Axioms (A1) and (A2), and $(M,\ga)$ is a balanced sutured manifold, the {\bf formal sutured homology} $\shg^g(M,\ga)$ is the projectively transitive system consisting of 
\benu
\item the $\hg$-homology $\hg(Y|R)$ for closures $(Y,R)$ of $(M,\ga)$ with a fixed and large enough genus $g$.
\item the canonical maps $\Phi$ between $\hg$-homologies as in Definition \ref{defn: canonical maps}.
\eenu
\edefn
\begin{conv}
Throughout the paper, when discussing formal sutured homology, we will pre-fix a large enough genus. So we omit it from the notation and write simply $\shg(M,\ga)$ instead of $\shg^g(M,\ga)$.
\end{conv}

\brem\label{rem: z2 grading}
When $\hg$ also satisfies Axiom (A3), since $\Phi$ is constructed by cobordism maps and their inverses, it is homogeneous with respect to the $\mathbb{Z}_2$-grading from Axiom (A3). Then there exists an induced relative $\mathbb{Z}_2$-grading on $\sh(M,\ga)$.
\erem

In \cite{baldwin2016contact,baldwin2018khovanov}, Baldwin and Sivek proved the bypass exact triangle for sutured monopole homology and sutured instanton homology. Their proof can be exported to our setup.

\begin{figure}[ht]
\centering
\begin{overpic}[width=0.5\textwidth]{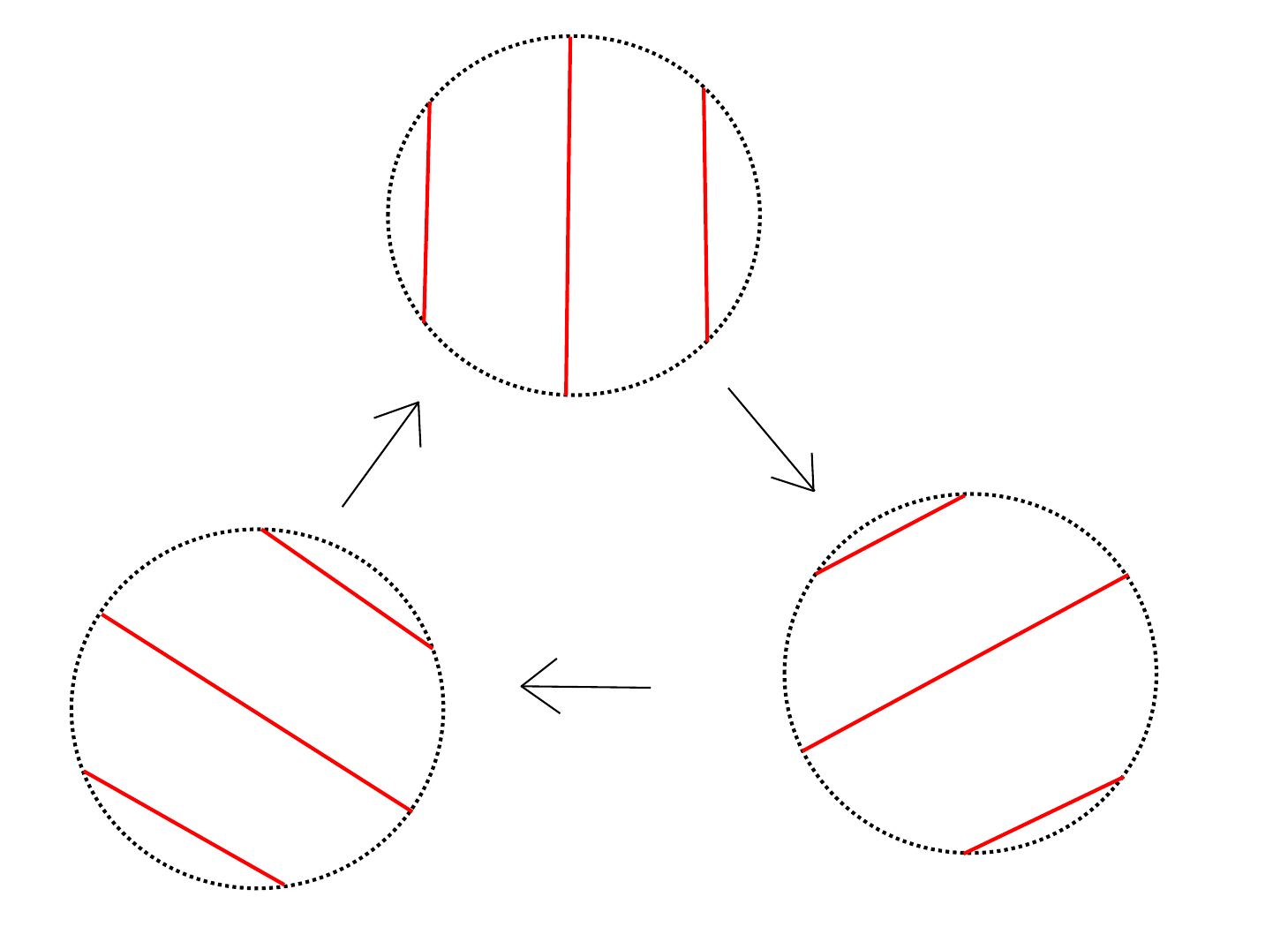}
    \put(16,22){$\ga_3$}
    \put(38,58){$\ga_1$}
    \put(73,28){$\ga_2$}
\end{overpic}
\vspace{-0.05in}
\caption{The bypass triangle.}\label{fig: the bypass triangle}
\end{figure}

\bthm[{\cite[Theorem 5.2]{baldwin2016contact} and \cite[Theorem 1.21]{baldwin2018khovanov}}]\label{thm_2: bypass exact triangle on general sutured manifold}
Suppose $(M,\ga_1)$, $(M,\ga_2)$, $(M,\ga_3)$ are three balanced sutured manifold such that the underlying 3-manifold is the same, and the sutures $\ga_1$, $\ga_2$, and $\ga_3$ only differ in a disk as depicted in Figure \ref{fig: the bypass triangle}. Then there exists an exact triangle
\begin{equation}\label{eq: bypass exact triangle}
\xymatrix@R=6ex{
\shg(-M,-\ga_1)\ar[rr]^{\psi_1}&&\shg(-M,-\ga_2)\ar[dl]^{\psi_2}\\
&\shg(-M,-\ga_3)\ar[lu]^{\psi_3}&
}    
\end{equation}

Moreover, the maps $\psi_i$ are induced by cobordisms and hence are homogeneous with respect to the relative $\mathbb{Z}_2$-grading on $\sh(M,\ga_i)$.
\ethm

We un-package the proof of Theorem \ref{thm_2: bypass exact triangle on general sutured manifold} for later convenience. 
\bprop[{\cite[Section 5]{baldwin2016contact}} and {\cite[Section 4]{baldwin2018khovanov}}]\label{prop_2: description of bypass maps}
Consider $(M,\ga_i)$ for $i=1,2,3$ in Theorem \ref{thm_2: bypass exact triangle on general sutured manifold}, there is a closure $(Y_1,R)$ of $(-M,-\ga_1)$ with the following significance.
\begin{enumerate}[(1)]
    \item The genus $g(R)$ is large enough.
    \item There are pairwise disjoint curves $\zeta_1,\zeta_2,\zeta_3\subset Y_1$ so that the following is true.
\begin{enumerate}[(a)]
    \item For $i=1,2,3$, we have $\zeta_i\cap{\rm int}(M)=\emptyset$ and $\zeta_i$ can be isotoped to be disjoint from $R$.
    \item If we perform a suitable Dehn surgery along $\zeta_1$, then we obtain a closure $(Y_2,R)$ of $(-M,-\ga_2)$. If we perform a suitable Dehn surgery in $Y_2$ along $\zeta_2$, then we obtain a closure $(Y_3,R)$ of $(-M,-\ga_3)$. If we perform a suitable Dehn surgery in $Y_3$ along $\zeta_3$, then we obtain the closure $(Y_1,R)$ of $(-M,-\ga_1)$ again.
    \item The maps $\psi_1$, $\psi_2$, and $\psi_3$ are induced the cobordism associated to Dehn surgeries along $\zeta_1$, $\zeta_2$, and $\zeta_3$, respectively.
\end{enumerate}
\item There are two curves $\eta_1$ and $\eta_2$ on $R$, so that if we perform $(-1)$-surgeries on both of them, with respect to the surface framings from $R$, then the surgeries along $\zeta_1$, $\zeta_2$, and $\zeta_3$ as stated in (b) lead to an exact triangle as in Axiom (A2).
\end{enumerate}
\eprop

\subsection{Gradings on formal sutured homology}\label{subsec: gradings on SH}\quad

Suppose $(M,\ga)$ is a balanced sutured manifold and $S\subset (M,\ga)$ is a properly embedded surface in $M$. If $S$ satisfies some admissible conditions, the first author \cite{li2019direct} constructed a $\mathbb{Z}$-grading on $SHM(M,\ga)$ and $SHI(M,\ga)$. In this subsection, we adapt the construction to formal sutured homology $\shg(M,\ga)$. 

\bdefn[{\cite{li2019decomposition}}]\label{defn_2: admissible surfaces}
Suppose $(M,\ga)$ is a balanced sutured manifold and $S\subset M$ is a properly embedded surface. The surface $S$ is called an \textbf{admissible surface} if the following holds.
\begin{enumerate}[(1)]
    \item Every boundary component of $S$ intersects $\ga$ transversely and nontrivially.
    \item $\frac{1}{2}|S\cap \ga|-\chi(S)$ is an even integer.
\end{enumerate}
\edefn
Recall the construction of a closure of $(M,\ga)$ in Definition \ref{defn: closure}. Let $T$ be a connected compact oriented surface of large enough genus and $\partial T\cong -\ga.$
Then we take
$$\widetilde{M}=M\cup[-1,1]\times T,\text{ with}~\partial \widetilde{M}=\widetilde{R}_+\sqcup(-\widetilde{R}_-).$$

Suppose $n=\frac{1}{2}|\partial S\cap \ga|$ and $\partial S\cap \ga=\{p_1,\dots,p_{2n}\}.$

\bdefn[\cite{li2019direct}]\label{defn: pairing}
A {\bf pairing} $\mathcal{P}$ of size $n$ is a collection of $n$ couples
$$\mathcal{P}=\{(i_1,j_1),\dots,(i_n,j_n)\}$$
such that the following holds.
\begin{enumerate}
\item $\{1,\dots,2n\}=\{i_1,j_1,\dots,i_n,j_n\}$.
\item For any $k\in\{1,\dots,n\}$, the points $p_{i_k}$ and $p_{j_k}$ are positive and negative, respectively, as intersection points of oriented curves $\partial S$ and $\ga$ on $\partial M$.
\end{enumerate}
\edefn

Given a pairing $\mathcal{P}$ of size $n$, and assuming that $g(T)$ is large enough, we can extend $S$ to a properly embedded surface in $\widetilde{M}$ as follows. Let $\al_1$,\dots,$\al_n$ be pairwise disjoint properly embedded arcs on $T$ such that the following holds.
\begin{enumerate}
	\item The arcs $\al_1,\dots,\al_n$ represent linearly independent homology classes in $H_1(T,\partial T)$.
	\item For any $k\in\{1,\dots,n\}$, we have $\partial \al_i=\{p_{i_k},p_{j_k}\}$.
\end{enumerate}
Given such $\al_1,\dots,\al_n$, take
$$\widetilde{S}_{\mathcal{P}}\deq S\cup[-1,1]\times(\al_1\cup\dots\cup\al_n).$$
Then $\widetilde{S}_{\mathcal{P}}$ is a properly embedded surface inside $\widetilde{M}$.

\bdefn\label{defn: balanced pairing}
A pairing $\mathcal{P}$ is called {\bf balanced} if $\widetilde{S}_{\mathcal{P}}\cap \widetilde{R}_+$ and $\widetilde{S}_{\mathcal{P}}\cap \widetilde{R}_-$ have the same number of components.
\edefn

For any balanced pairing $\mathcal{P}$, we can pick an orientation preserving diffeomorphism $h:\widetilde{R}_+\xra{\cong} \widetilde{R}_-$ so that
$$h(\widetilde{S}_{\mathcal{P}}\cap \widetilde{R}_+)=\widetilde{S}_{\mathcal{P}}\cap \widetilde{R}_-.$$
Thus, we obtain a closed oriented surface $\bar{S}_{\mathcal{P}}\subset Y$ in the closure $(Y,R)$ induced by $h$. Define
$$\shg(M,\ga,S,i)\deq \hg(Y,([R],[\bar{S}_{\mathcal{P}}]),(2g(R)-2,2i)).$$
\bthm\label{thm: grading is well defined, appendix}
Given an admissible surface $S$ in a balanced sutured manifold $(M,\ga)$, the decomposition$$\shg(M,\ga)=\bigoplus_{i\in \mathbb{Z}}\shg(M,\ga,S,i)$$is independent of all the choices made in the construction and hence is well-defined.
\ethm

\brem
As mentioned in the convention after Definition \ref{defn: conv}, when writing $\shg$, we actually mean $\shg^g$ for some large and fixed integer $g$. This means that all closures involved have the same genus $g$.
\erem

\bpf[Proof of Theorem \ref{thm: grading is well defined, appendix}]
The decomposition follows from Axioms (A1-1) and (A1-7). This gives a $\mathbb{Z}$-grading on $\shg(M,\ga)$. To show that this grading is well-defined, we need to show that it is independent of the following three types of choices:
\begin{enumerate}
\item the choice of the balanced pairing $\mathcal{P}$,
\item the choice of arcs $\al_1$,\dots,$\al_n$ with fixed endpoints,
\item the choice of the diffeomorphism $h$.	
\end{enumerate}

In \cite[Section 3.1]{li2019direct}, the grading has been shown to be independent of the choices of type (2) and (3). The proof involves only Axioms (A1) and (A2) and hence can be applied to our current setup. However, the original argument for choices of type (1) in \cite[Section 3.3]{li2019direct} involves closures of different genera, which is beyond the scope of our current paper as mentioned in Remark \ref{rem: grading is well-defined}. Hence, we provide an alternative proof here. For the moment, let us write the grading as
$$\shg(M,\ga,S,\mathcal{P},i)$$
to emphasize that the grading \textit{a priori} depends on the choice of the balanced pairing. Theorem \ref{thm: grading is well defined, appendix} then follows from the following proposition.
\epf

\bprop\label{prop: independence of balanced pairings}
Suppose $\mathcal{P}$ and $\mathcal{P}^\p$ are two balanced pairings, then for any $i\in\intg$, we have
$$\shg(M,\ga,S,\mathcal{P},i)=\shg(M,\ga,S,\mathcal{P}^\p,i).$$
\eprop

To relate two different pairings, in \cite{li2019direct}, the author introduced the following operation.
\bdefn
Suppose $\mathcal{P}$ is a pairing of size $n$ and $\al_1,\dots,\al_n$ are related arcs. Suppose $k,l\in\{1,\dots,n\}$ are two indices so that the following holds.
\begin{enumerate}
\item The arcs $\{1\}\times\al_{i_k}$ and $\{1\}\times\al_{i_l}$ belong to different components of $\widetilde{S}_{\mathcal{P}}\cap R_+$.
\item The arcs $\{-1\}\times\al_{i_k}$ and $\{-1\}\times\al_{i_l}$ belong to different components of $\widetilde{S}_{\mathcal{P}}\cap R_-$.
\end{enumerate}
Then we can construct another pairing
$$\mathcal{P}^\p=(\mathcal{P}\backslash\{(i_k,j_k),(i_l,j_l)\})\cup\{(i_k,j_l),(i_l,j_k)\}.$$

The operation of replacing $\mathcal{P}$ by $\mathcal{P}^\p$ is called a {\bf cut and glue operation}.
\edefn

\bthm[\cite{kavi2019cutting}]\label{thm: cut and glue}
Balanced pairings always exist. Moreover, any two balanced pairings are related by a finite sequence of cut and glue operations and their inverses.
\ethm

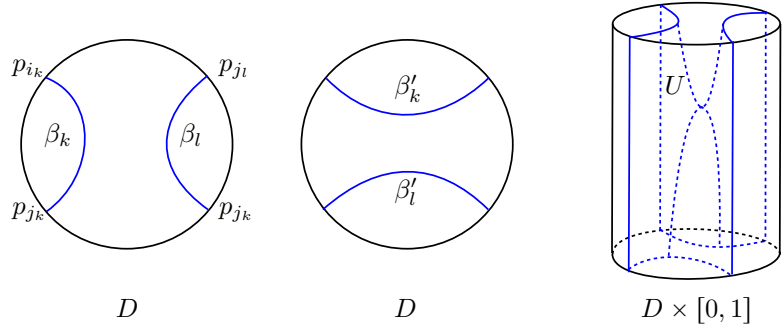
\begin{figure}[ht]
\centering
\begin{overpic}[width=0.7\textwidth]{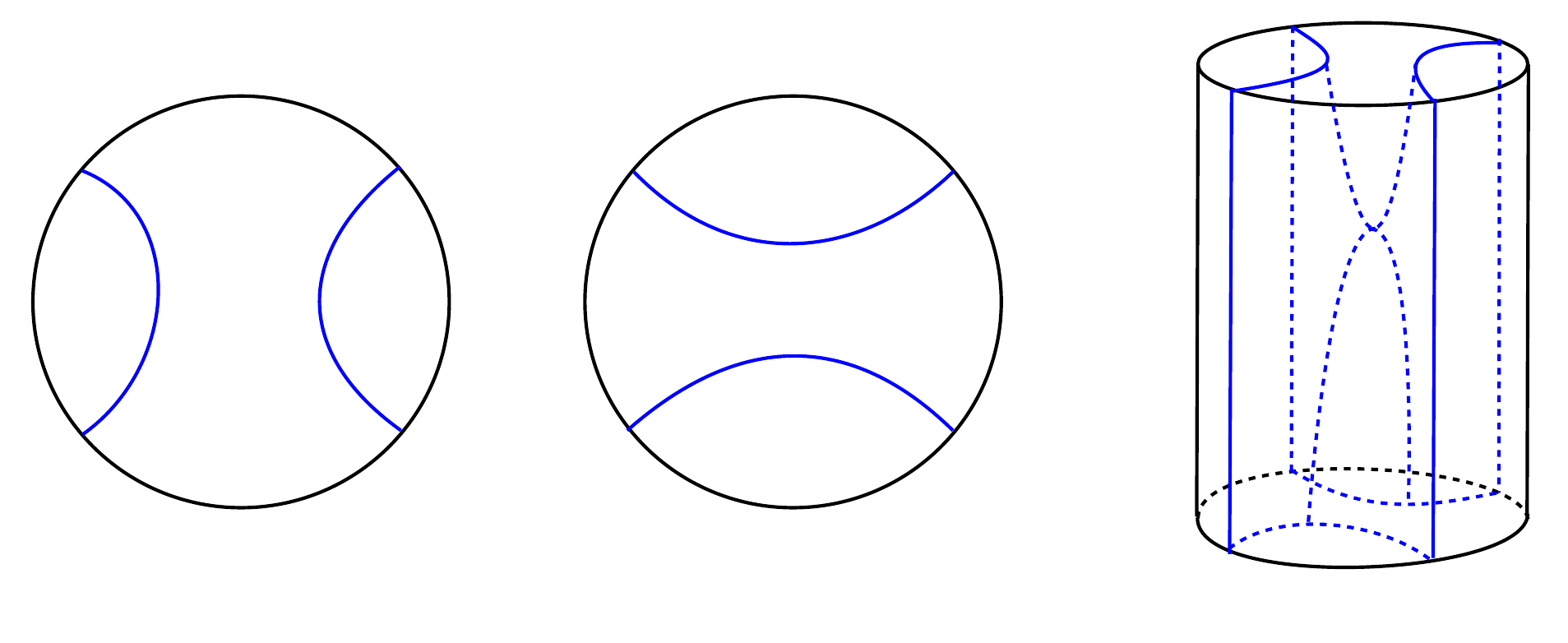}
	\put(1,31){$p_{i_k}$}
	\put(1,13){$p_{j_k}$}
	\put(27,31){$p_{j_l}$}
	\put(27,13){$p_{j_k}$}
    \put(22,22){$\be_l$}
    \put(5,22){$\be_k$}
    \put(49,28){$\be_k'$}
    \put(49,15){$\be_l'$}
    \put(83,28){$U$}
    \put(14,0){$D$}
    \put(49,0){$D$}
    \put(80,0){$D\times[0,1]$}
\end{overpic}
\vspace{-0.05in}
\caption{The disk $D$, the arcs $\be_k,\be_l,\be_k',\be_l'$, and the surface $U$.}\label{fig: cut and glue}
\end{figure}

\blem\label{lem: cut and glue}
Suppose $\mathcal{P}$ and $\mathcal{P}^\p$ are two balanced pairings that are related by a cut and glue operation, then for any $i\in\intg$, we have
$$\shg(M,\ga,S,\mathcal{P},i)=\shg(M,\ga,S,\mathcal{P}^\p,i).$$
\elem
\bpf
Suppose $k$ and $l$ are the indices involved in the operation. From the first part of the proof of Theorem \ref{thm: grading is well defined, appendix}, we can freely make choices of type (2) and (3). Hence we can assume that there is a disk $D\subset {\rm int}(T)$ so that $\al_k$ and $\al_l$ intersects $D$ in two arcs as depicted in Figure \ref{fig: cut and glue}. Suppose
$$D_+=\{1\}\times D\subset R_+{~\rm and~} D_-=\{-1\}\times D\subset R_-.$$
We can choose an orientation preserving diffeomorphism $h:\widetilde{R}_+\xra{\cong} \widetilde{R}_-$ such that 
$$h(\widetilde{S}_{\mathcal{P}}\cap \widetilde{R}_+)=\widetilde{S}_{\mathcal{P}}\cap \widetilde{R}_-{\rm~and~} h(D_+)=D_-.$$
Let $(Y,R)$ be the corresponding closure of $(M,\ga)$ and $\widebar{S}_{\mathcal{P}}$ be the closed surface defining the grading $\shg(M,\ga,S,\mathcal{P})$. Let $\be_k=\al_k\cap D$ and $\be_l=\al_l\cap D$. It is straightforward to check that if we remove the two arcs $\be_k$ and $\be_l$ from $D\subset T$, and glue back two new arcs $\be_k'$ and $\be_l'$ as shown in the middle subfigure of Figure \ref{fig: cut and glue}, then we obtain two new properly embedded arcs $\al_k'$ and $\al_l'$ on $T$ so that
$$\partial \al_k'=\{p_{i_k},p_{j_l}\}~{\rm and~}\partial \al_l'=\{p_{i_l},p_{j_k}\}.$$
Hence we change from $\mathcal{P}$ to $\mathcal{P}^\p$. Inside $Y$, if we remove $S^1\times (\be_k\cup\be_l)\subset S^1\times D\subset Y$ and glue back $S^1\times (\be_k'\cup\be_l')\subset S^1\times D$, then we obtain the surface $\widebar{S}_{\mathcal{P}^\p}\subset Y$ that gives rise to the grading $\shg(M,\ga,S,\mathcal{P}^\p)$. The lemma then follows from the fact
\begin{equation}\label{eq: P and P'}
	[\widebar{S}_{\mathcal{P}}]=[\widebar{S}_{\mathcal{P}^\p}]\in H_2(Y).
\end{equation}

The equality (\ref{eq: P and P'}) can be proved by constructing an explicit cobordism in $Y\times[0,1]$ from $\widebar{S}_{\mathcal{P}}\subset Y\times\{0\}$ to $\widebar{S}_{\mathcal{P'}}\subset Y\times\{1\}$: in the product $(Y\times[0,1], \widebar{S}_{\mathcal{P}}\times [0,1])$, we can remove $$S^1\times(\be_1\cup\be_2)\times[0,1]\subset S^1\times D\times [0,1]\subset Y\times[0,1]$$ and glue back $S^1\times U\subset D\times[0,1]$, where $U\subset D\times[0,1]$ is the surface shown in the right subfigure of Figure \ref{fig: cut and glue}.
\epf

\bpf[Proof of Proposition \ref{prop: independence of balanced pairings}]
It follows immediately from Theorem \ref{thm: cut and glue} and Lemma \ref{lem: cut and glue}.
\epf

Having constructed the grading, the rest of the arguments in \cite[Section 3.3]{li2019direct} can be applied to our current setup verbatim. Hence we have the following.

\bprop[\cite{li2019direct,li2018gluing}]\label{thm_2: grading in SHG}
Suppose $(M,\ga)$ is a balanced sutured manifold and $S\subset (M,\ga)$ is an admissible surface. Then there is a $\intg$-grading on $\shg(M,\ga)$ induced by $S$ , which we write as
$$\shg(M,\ga)=\bigoplus_{i\in\intg}\shg(M,\ga,S,i).$$
This decomposition satisfies the following properties.

\begin{enumerate}[(1)]
    \item Suppose $n=\frac{1}{2}|\partial S\cap\ga|$. If $|i|>\frac{1}{2}(n-\chi(S))$, then $\shg(M,\ga,S,i)=0.$
    \item Suppose $S$ is a \textbf{product disk}, \textit{i.e.} $S\cong D^2$ and $|\partial S\cap \ga|=2$. Let $(M',\ga')$ be the balanced sutured manifold under the sutured manifold decomposition $(M,\ga)\stackrel{S}{\leadsto}(M',\ga')$ in the sense of Gabai \cite{gabai1983foliations} (this is called \textbf{product decomposition} in \cite[Definition 9.11]{juhasz2006holomorphic}). Then we have
$$\shg(M,\ga,S,0)\cong \shg(M',\ga').$$
    \item For any $i\in\intg$, we have
$$\shg(M,\ga,S,i)=\shg(M,\ga,-S,-i).$$
    \item For any $i\in\intg$, we have
$$\shg(M,-\ga,S,i)\cong \shg(M,\ga,S,-i).$$
    \item For any $i\in\intg$, we have
$$\shg(-M,\ga,S,i)\cong {\rm Hom}_\mathbb{F}(\shg(M,\ga,S,-i),\mathbb{F}).$$
\end{enumerate}
\eprop

\bpf
Term (1) comes from the adjunction inequality in (A1-4). Term (2) follows from discussion in \cite[Section 4.2]{baldwin2016contact} (see also \cite[Section 6]{kronheimer2010knots}). Note that decomposing along such a disk is the inverse operation of attaching a product 1-handle, and the disk is precisely the co-core of the product 1-handle. Term (3) is straighforward from the definition. Term (4) is from Axiom (A1-3). Term (5) is from the pairing (\textit{c.f.} \cite{li2018gluing}):$$\langle\cdot,\cdot\rangle:\shg(M,\ga)\times \shg(-M,\ga)\to \mathbb{F}.$$
\epf
\brem\label{rem: general decomposition theorem}
We do not state a general surface decomposition theorem in term (2) of Proposition \ref{thm_2: grading in SHG} because the proof of the general theorem in \cite[Section 6]{kronheimer2010knots} involves closures of different genera. Indeed, proofs in Section \ref{sec: Equivalence of graded Euler characteristics} of the current paper only involve the decomposition theorem for a product disk.

Though not used in this paper, it is worth mentioning here that the general decomposition theorem is still valid when we consider admissible surfaces and also take the genus of the closures into account. Suppose $S\subset (M,\ga)$ is an admissible surface and let $\shg^g$ denote the formal sutured Floer homology obtaining from closures of genus $g$. Suppose $(M',\gamma')$ is obtained from $(M,\ga)$ by decomposing along $S$. Then for large enough $g$, we have 
$$\shg^g(M,\ga,S,\frac{1}{4}|S\cap\ga|-\frac{1}{2}\chi(S))\cong \shg^{g'}(M',\ga'),$$
where $g'=g+\frac{1}{4}|S\cap\ga|-\frac{1}{2}\chi(S)$.
The argument in \cite[Section 6]{kronheimer2010knots} applies verbatim to prove the above isomorphism.
\erem

\bdefn[\cite{juhasz2006holomorphic}]
A sutured manifold $(M,\ga)$ is called \textbf{taut} if $M$ is irreducible and $R_{+}(\ga)$ and $R_{-}(\ga)$
are both incompressible and Thurston norm-minimizing in the homology class that they represent in $H_2(M,\ga)$.
\edefn

\bprop\label{prop: SHI detects tautness}
Suppose $(M,\ga)$ is a balanced sutured manifold so that $M$ is irreducible. If $(M,\ga)$ is not taut, then $\shg(M,\ga)=0$.
\eprop
\bpf
It follows from the construction and the adjunction inequality.
\epf

\brem
Though not needed in this paper, using the arguments in the proof of \cite[Theorem 1.4]{juhasz2008floer} and Remark \ref{rem: general decomposition theorem}, we could also show that an irreducible balanced sutured manifold $(M,\ga)$ is taut if and only if $\shg(M,\ga)\neq0$.
\erem

\bdefn
Suppose $(M,\ga)$ is a balanced sutured manifold. It is called a \textbf{homology product} if $H_1(M, R_+(\ga)) = 0$ and $H_1(M, R_-(\ga)) = 0$. It is called a \textbf{product sutured manifold} if $$(M,\ga)\cong ([-1,1]\times \Sigma,\{0\}\times\partial \Sigma),$$where $\Sigma$ is a compact surface with boundary.
\edefn
\bdefn[\cite{juhasz2006holomorphic}]
A balanced sutured manifold $(M,\ga)$ is called a \textbf{product sutured manifold} if $$(M,\ga)\cong ([-1,1]\times \Sigma,\{0\}\times\partial \Sigma),$$where $\Sigma$ is a compact surface with boundary.
\edefn
\bprop\label{prop: product}
Suppose $(M,\ga)$ is a product sutured manifold. Then $\shg(M,\ga)\cong \mathbb{F}$.
\eprop
\bpf
It follows from the construction and Axiom (A1-5).
\epf


If $S\subset(M,\ga)$ is not admissible, then we can perform an isotopy on $S$ to make it admissible. 

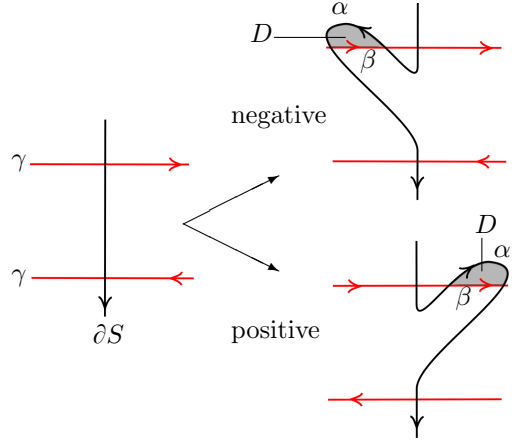
\begin{figure}[ht]
\centering
\begin{overpic}[width=2.5in]{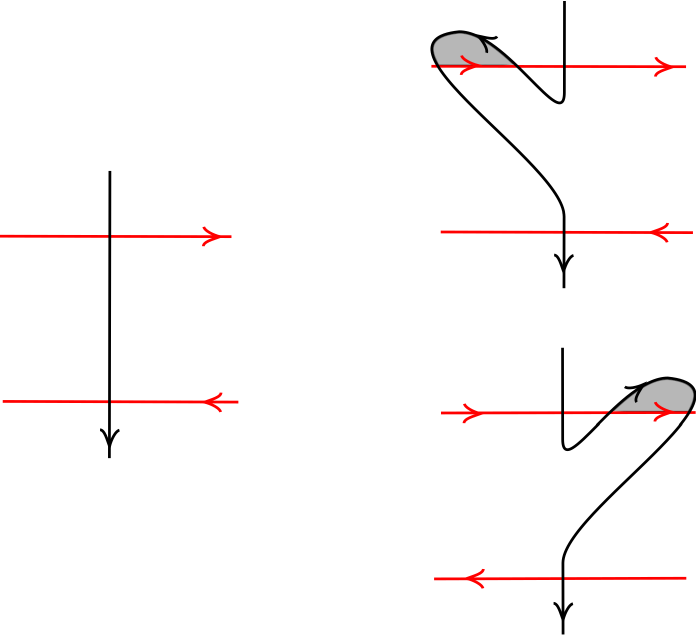}
    \put(13,20){$\partial{S}$}
    \put(-4,57){$\ga$}
    \put(-4,33){$\ga$}
    \put(32,45){\vector(2,1){20}}
    \put(32,45){\vector(2,-1){20}}
    \put(42,21){positive}
    \put(42,66){negative}
    \put(46,83){$D$}
    \put(51,84){\line(1,0){15}}
    \put(63,89){$\al$}
    \put(69,77.5){$\be$}
    \put(93,43){$D$}
    \put(94.5,42){\line(0,-1){7}}
    \put(97,38){$\al$}
    \put(89,28){$\be$}
\end{overpic}
\vspace{0.05in}
\caption{The positive and negative stabilizations of $S$.}\label{fig: pm_stabilization_of_surfaces}
\end{figure}

\bdefn\label{defn_2: stabilization of surfaces}
Suppose $(M,\ga)$ is a balanced sutured manifold, and $S$ is a properly embedded surface. A \textbf{stabilization} of $S$ is a surface $S^\p$ obtained from $S$ by isotopy in the following sense. This isotopy creates a new pair of intersection points:
$$\partial S'\cap\ga=(\partial{S}\cap\ga)\cup \{p_+,p_-\}.$$
We require that there are arcs $\al\subset \partial{S'}$ and $\be\subset \ga$, oriented in the same way as $\partial{S'}$ and $\ga$, respectively, and the following holds.
\begin{enumerate}[(1)]
    \item $\partial{\al}=\partial{\be}=\{p_+,p_-\}$.
    \item $\al$ and $\be$ cobound a disk $D$ with ${\rm int}(D)\cap (\ga\cup \partial{S}')=\emptyset$.
\end{enumerate}
The stabilization is called \textbf{negative} if $\partial{D}$ is the union of $\al$ and $\be$ as an oriented curve. It is called \textbf{positive} if $\partial{D}=(-\al)\cup\be$. See Figure \ref{fig: pm_stabilization_of_surfaces}. We denote by $S^{\pm k}$ the surface obtained from $S$ by performing $k$ positive or negative stabilizations, repsectively.
\edefn

\brem\label{rem_2: pm switches according to orientations of the suture}
The definition of stabilizations of a surface depends on the orientations of the suture and the surface. If we reverse the orientation of the suture or the surface, then positive and negative stabilizations switch between each other.
\erem

One can also relate the gradings associated to different stabilizations of a fixed surface. The proof for $SHM$ and $SHI$ in \cite{li2019direct,wang2020cosmetic} can be adapted to our setup as well.

\bthm[{\cite[Proposition 4.9]{li2019direct}} and {\cite[Proposition 4.17]{wang2020cosmetic}}]\label{thm_2: grading shifting property}
Suppose $(M,\ga)$ is a balanced sutured manifold and $S$ is a properly embedded surface in $M$ that intersects $\ga$ transversely. Suppose all the stabilizations mentioned below are performed on a distinguished boundary component of $S$. Then, for any $p,k,l\in \intg$ such that the stabilized surfaces $S^{p}$ and $S^{p+2k}$ are both admissible, we have
$$\shg(M,\ga,S^{p},l)=\shg(M,\ga,S^{p+2k},l+k).$$
Note that $S^p$ is a stabilization of $S$ as introduced in Definition \ref{defn_2: stabilization of surfaces}, and, in particular, $S^0=S$.
\ethm
\brem
The original form of Theorem \ref{thm_2: grading shifting property} in \cite{li2019direct} was stated for a Seifert surface in the case of a knot complement. However, it is straightforward to generalize the proof to the case of a general admissible surface in a general balanced sutured manifold, given the condition that the decompositions along $S$ and $-S$ are both taut. This extra condition on taut decompositions was then dropped due to the work in \cite{wang2020cosmetic}.
\erem

If we have multiple admissible surfaces, then they together induce a multi-grading. This is proved for $SHM$ and $SHI$ by Ghosh and the first author \cite{li2019decomposition}. The proof can be adapted to our case without essential changes.

\bthm[{\cite[Proposition 1.14]{li2019decomposition}}]\label{thm: zn grading}
Suppose $(M,\ga)$ is a balanced sutured manifold and $S_1,\dots,S_n$ are admissible surfaces in $(M,\ga)$. Then there exists a $\intg^n$-grading on $\shg(M,\ga)$ induced by $S_1,\dots,S_n$, which we write as
$$\shg(M,\ga)=\bigoplus_{(i_1,\dots,i_n)\in\intg^n}\shg(M,\ga,(S_1,\dots,S_n),(i_1,\dots,i_n)).$$
\ethm

\bthm[{\cite[Theorem 1.12]{li2019decomposition}}]\label{thm: shift}
Suppose $(M,\ga)$ is a balanced sutured manifold and $\al\in H_2(M,\partial M)$ is a nontrivial homology class. Suppose $S_1$ and $S_2$ are two admissible surfaces in $(M,\ga)$ such that$$[S_1]=[S_2]=\al\in H_2(M,\partial M).$$Then, there exists a constant $C$ so that $$\shg(M,\ga,S_1,l)=\shg(M,\ga,S_2,l+C).$$
\ethm

Based on the relative $\mathbb{Z}_2$-grading from Remark \ref{rem: z2 grading} and the $\mathbb{Z}^n$-grading from Theorem \ref{thm: zn grading}, we can define graded Euler characteristic of formal sutured homology.


\bdefn\label{defn: chi(shg)}
Suppose $(M,\ga)$ is a balanced sutured manifold and $S_1,\dots,S_n$ are admissible surfaces in $(M,\ga)$ such that $[S_1],\dots,[S_n]$ generate $H_2(M,\partial M)$. For $i=1,\dots,n$, let $\rho_i\in H=H_1(M)/{\rm Tors}$ be the class satisfying $\rho_i\cdot S_j=\delta_{i,j}$. Define the \textbf{graded Euler characteristic} of $\shg(M,\ga)$ to be $$\chi_{\rm gr}(\shg(M,\ga)) \deq\sum_{(i_1,\dots,i_n)\in\intg^n}\chi(\shg(M,\ga,(S_1,\dots,S_n),(i_1,\dots,i_n)))\cdot (\rho_1^{i_1}\cdots\rho_n^{i_n})\in \mathbb{Z}[H]/\pm H.$$
\edefn
\brem\label{rem: fix closure}
It can be shown by Theorem \ref{thm: shift} that the definition of graded Euler characteristic is independent of the choices of $S_1,\dots,S_n$ if we regard it as an element in $\mathbb{Z}[H]/\pm H$. If the admissible surfaces $S_1,\dots,S_n$ and a particular closure of $(M,\ga)$ is fixed, then the ambiguity of $\pm H$ can be removed.
\erem

From Theorem \ref{thm_2: bypass exact triangle on general sutured manifold}, Proposition \ref{prop_2: description of bypass maps}, and Axiom (A1-7), the following proposition is straightforward.
\bprop\label{prop: bypass maps preserves gradings}
Suppose $(M,\ga)$ is a balanced sutured manifold and $S\subset (M,\ga)$ is an admissible surface. Suppose the disk as in Figure \ref{fig: the bypass triangle}, where we perform the bypass change, is disjoint from $\partial S$. Let $\ga_2$ and $\ga_3$ be the resulting two sutures. Then all the maps in the bypass exact triangle (\ref{eq: bypass exact triangle}) are grading preserving, \textit{i.e.}, for any $i\in\intg$, we have an exact triangle
\begin{equation*}
\xymatrix@R=6ex{
\shg(-M,-\ga_1,S,i)\ar[rr]^{\psi_{1,i}}&&\shg(-M,-\ga_2,S,i)\ar[dl]^{\psi_{2,i}}\\
&\shg(-M,-\ga_3,S,i)\ar[lu]^{\psi_{3,i}}&
}    
\end{equation*}
where $\psi_{k,i}$ are the restriction of $\psi_k$ in (\ref{eq: bypass exact triangle}).
\eprop


\section{Heegaard Floer homology and the graph TQFT}\label{sec: hf}

In this section, we discuss the modification of Heegaard Floer theory to make it suitable to formal sutured homology.

\subsection{Heegaard Floer homology for multi-pointed 3-manifolds}\quad

In this subsection and the next subsection, we provide an overview of the graph TQFT for Heegaard Floer theory, constructed by Zemke \cite{Zemke2019} (see also \cite{Hendricks2018,Zemke2018}), and list some properties which are relevant to proofs in the third subsection about Floer's excision theorem.

\bdefn
A \textbf{multi-pointed 3-manifold} is a pair $(Y,\bf{w})$ consisting of a closed, oriented 3-manifold $Y$ (not necessarily connected), together with a finite collection of basepoints $\bs{w}\subset Y$, such that each component of $Y$ contains at least one basepoint.

Given two multi-pointed 3-manifolds $(Y_1,\bs{w}_1)$ and $(Y_2,\bs{w}_2)$, a \textbf{ribbon graph cobordism} from $(Y_1,\bs{w}_1)$ to $(Y_2,\bs{w}_2)$ is a pair $(W,\Ga)$ satisfying the following conditions.
\benu
\item $W$ is a cobordism from $Y_1$ to $Y_2$.
\item $\Ga$ is an embedded graph in $W$ such that $\Ga \cap Y_i = \bs{w}_i$ for $i=1,2$. Furthermore, each point of $\bs{w}_i$ has valence 1 in $\Ga$.
\item $\Ga$ has finitely many edges and vertices, and no vertices of valence 0.
\item The embedding of $\Ga$ is smooth on each edge.
\item $\Ga$ is decorated with a formal ribbon structure, \textit{i.e.}, a formal choice of cyclic ordering of the edges adjacent to each vertex.
\eenu
\edefn
\bdefn\label{defn: restrict}
A ribbon graph cobordism $(W,\Ga)$ from $(Y_1,\bs{w}_1)$ to $(Y_2,\bs{w}_2)$ is called a \textbf{restricted graph cobordism} if $W$ is obtained from $Y_1\times I$ by attaching 4-dimensional 1-, 2-, and 3-handles away from all basepoints and $\Ga=\bs{w}_1\times I$ is the induced graph in $W$ (so the cyclic ordering is unique and $|\bs{w}_1|=|\bs{w}_2|$).
\edefn

\bdefn[{\cite[Definition 4.1]{Zemke2019}}]
Suppose $(Y,\bs{w})$ is a connected multi-pointed 3-manifold. A \textbf{multi-pointed Heegaard diagram} $\mch = (\Sigma,\alpha,\beta,\bs{w})$ for $(Y,\bs{w})$ is a tuple satisfying the following conditions.
\benu
\item $\Sigma$ is a closed, oriented surface, embedded in $Y$, such that $\bs{w}\subset \Sigma\backslash(\al\cup \be)$. Furthermore, $\Sigma$ splits $Y$ into two handlebodies $U_\al$ and $U_\be$, oriented so that $\Sigma = \partial U_\al=-U_\be$.
\item $\al=\{\al_1,\dots,\al_n\}$ is a collection of $n = g(\Sigma)+|\bs{w}|-1$ pairwise disjoint simple closed curves on $\Sigma$, bounding pairwise disjoint compressing disks in $U_\al$. Each component of $\Sigma\backslash\al$ is planar
and contains a single basepoint.
\item $\be = \{\be_1,\dots,\be_n\}$ is a collection of pairwise disjoint, simple, closed curves on $\Sigma$ bounding pairwise disjoint compressing disks in $U_\be$. Each component of $\Sigma\backslash\be$ is planar and contains a single basepoint.
\eenu
\edefn
Suppose $\bs{w}=\{w_1,\dots,w_m\}$. Let the polynomial ring associated to $\bs{w}$ be $$\mathbb{F}_2[U_{\bs{w}}]\deq \ft[U_{w_1},\dots,U_{w_m}].$$Let $\ft[U_{\bs{w}},U_{\bs{w}}^{-1}]$ be the ring obtained by formally inverting each of the variables. 

If $\bs{k}=(k_1,\dots,k_m)$ is an $m$-tuple, let $$U_{\bs{w}}^{\bs{k}}\deq U_{w_1}^{k_1}\cdots U_{w_m}^{k_m}.$$For simplicity, we will also write $U_i$ for $U_{w_i}$.

Suppose $\mch=(\Sigma,\al,\be,\bs{w})$ is a multi-pointed Heegaard diagram of a connected multi-pointed 3-manifold $(Y,\bs{w})$. Suppose $n = g(\Sigma)+|\bs{w}|-1$. Consider two tori$$\mathbb{T}_\al\deq \al_1\times\cdots\times\al_n~{\rm and}~\mathbb{T}_\be\deq \be_1\times\cdots\times\be_n$$in the symmetric product $${\rm Sym}^n\Sigma\deq (\prod_{i=1}^n\Sigma)/S_n.$$

The chain complex $CF^-(\mch)$ is a free $\ft[U_{\bs{w}}]$-module generated by intersection points $\bs{x}\in\mathbb{T}_\al\cap \mathbb{T}_\be$. Define $$CF^\infty(\mch)\deq CF^-(\mch)\otimes_{\ft[U_{\bs{w}}]}\ft[\uw,\uww]\aand CF^+(\mch)\deq CF^\infty(\mch)/CF^-(\mch).$$

To construct a differential on $CF^-(\mch)$, suppose $\mch$ satisfies some extra admissibility conditions if $b_1(Y)>0$ (\textit{c.f.} \cite[Section 4.7]{Zemke2019}). Let $(J_s)_{s\in[0,1]}$ be an auxiliary path of almost complex structures on $\sym^n\Sigma$ and let $\pi_2(\bs{x},\bs{y})$ be the set of homology classes of Whitney disks connecting intersection points $\bs{x}$ and $\bs{y}$ (\textit{c.f.} \cite[Section 3.4]{ozsvath2008multivariable}). For $\phi\in \pi_2(\bs{x},\bs{y})$, let $\mathcal{M}_{J_s}(\phi)$ be the moduli space of $J_s$-holomorphic maps $u:[0,1]\times \mathbb{R}\to \sym^n\Sigma$ which represent $\phi$. The moduli space $\mathcal{M}_{J_s}(\phi)$ has a natural action of $\mathbb{R}$, corresponding to reparametrization of the source. We write
$$\widehat{\mathcal{M}}_{J_s}(\phi)\deq \mathcal{M}_{J_s}(\psi)/\mathbb{R}.$$

For $\phi\in \pi_2(\bs{x},\bs{y})$, let $\mu(\phi)$ be the expected dimension of $\mathcal{M}_{J_s}(\phi)$ for generic $J_s$ and let $n_{w_i}(\phi)$ be the algebraic intersection number of $\{w_i\}\times \sym^{n-1}\Sigma$ and any representative of $\phi$. Define $$n_{\bs{w}}(\phi)\deq (n_{w_1}(\phi),\dots,n_{w_m}(\phi)).$$

For a generic path $J_s$, define the differential on $CF^-(\mch)$ by$$\partial_{J_s}(\bs{x})=\sum_{\bs{y}\in \mathbb{T}_\al\cap \mathbb{T}_\be}\sum_{\substack{\phi\in\pi_2(\bs{x},\bs{y})\\\mu(\phi)=1}}\#\widehat{\mathcal{M}}_{J_s}(\phi)U_{\bs{w}}^{n_{\bs{w}}(\phi)}\cdot \bs{y},$$extended linearly over $\ft[\uw].$ The differential $\partial_{J_s}$ can be extended on $CF^\infty(\mch)$ and $CF^+(\mch)$ by tensoring with the identity map.
\blem[{\cite[Lemma 4.3]{ozsvath2008multivariable}}]
For a generic path $J_s$, the map $\partial _{J_s}$ on $CF^\circ(\mch)$, where $\circ\in\{\infty,+,-\}$, satisfies$$\partial _{J_s}\circ \partial _{J_s}=0.$$
\elem
For a disconnected multi-pointed 3-manifold $(Y,\bs{w})=(Y_1,\bs{w_1})\sqcup(Y_2,\bs{w_2})$, where $Y_i$ is connected for $i=1,2$, suppose $\mch_i$ is an admissible multi-pointed Heegaard diagram of $Y_i$ and suppose $J_{s_i}$ are corresponding generic paths of almost complex structures. For $\circ\in\{\infty,+,-\}$, let the chain complex associated to $(Y,\bs{w})$ be\begin{equation}\label{eq: circ}
    (CF^\circ(\mch_1\sqcup\mch_2),\partial_{J_s})\deq (CF^\circ(\mch_1),\partial_{J_{s_1}})\ot_\ft (CF^\circ(\mch_2),\partial_{J_{s_2}}).
\end{equation}
\brem\label{rem: color}
In Zemke's original construction \cite[Section 4.3]{Zemke2019}, one should choose colors for basepoints and graphs to achieve the functoriality of the TQFT. For basepoints with the same color, the corresponding $U$-variables should be the same. In above notations, we implicitly choose different colors for all basepoints so that the $U$-variable for each basepoint is different. This is to obtain the following relation on the homology level \begin{equation}\label{eq: homology}H(CF^\circ(\mch_1\sqcup\mch_2),\partial_{J_s})=H(CF^\circ(\mch_1),\partial_{J_{s_1}})\ot_\ft H(CF^\circ(\mch_1),\partial_{J_{s_2}}).\end{equation}Note that in the construction of \cite{Hendricks2018,Zemke2018}, the colors of all basepoints are the same and all $U$-variables are identified as $U$, so (\ref{eq: circ}) should be a tensor product over $\ft[U]$ rather than $\ft$ and (\ref{eq: homology}) does not hold in general. 
\erem
\brem\label{rem: color2}
Given a finite set of multi-pointed 3-manifolds and ribbon graph cobordisms, the chain complex $CF^-(\emptyset)$ is set to be $\ft[U_{\bs{w}}]$, where $U_{\bs{w}}$ contains all $U$-variables associated to basepoints in the set. For any multi-pointed 3-manifold $(Y,\bs{w}^\p)$ with $\bs{w}^\p\subset\bs{w}$ that is in the given set, the actual chain complex in the TQFT should be $$CF^-(Y,\bs{w}^\p)\otimes_{\ft}\ft[U_{\bs{w}\backslash\bs{w}^\p}].$$In the statements of results in this paper, we always have $\bs{w}^\p=\bs{w}$ for any multi-pointed 3-manifold $(Y,\bs{w}^\p)$. However, in the proof of those results (\textit{e.g.} Lemma \ref{lem: two cobordism equal} and Theorem \ref{thm: floer excision}), we may have multi-pointed 3-manifold $(Y,\bs{w}^\p)$ such that $\bs{w}^\p\neq \bs{w}$; see Remark \ref{rem: mcr}. Also, in the proof, the colors of basepoints may be different.
\erem
The chain homotopy type of $(CF^\circ(\mch),\partial_{J_s})$ is independent of the choices of the admissible diagram $\mch$ and the generic path $J_s$. Indeed, we have the following theorem about naturality.
\bthm[{\cite[Proposition 4.6]{Zemke2019}, see also \cite{ozsvath2004holomorphic,Juhasz2012}}]\label{thm: naturality hf}
Suppose that $(Y,\bs{w})$ is a multi-pointed 3-manifold. To each (admissible) pairs $(\mch,J)$ and $(\mch^\p,J^\p)$, there is a well-defined map
$$\Psi_{(\mch,J)\to (\mch^\p,J^\p)}:(CF^-(\mch),\partial_J)\to CF^-(\mch^\p),\partial_{J^\p}),$$
which is well-defined up to $\ft[U_{\bs{w}}]$-equivariant chain homotopy. Furthermore, the following holds.
\benu
\item If $(\mch, J)$, $(\mch^\p,J^\p)$ and $(\mch^\pp ,J^\pp)$ are three pairs, then there is a chain homotopy equivalence
$$\Psi_{(\mch,J)\to(\mch^\pp,J^\pp)}\simeq \Psi_{(\mch^\p,J^\p)\to(\mch^\pp,J^\pp)}\circ\Psi_{(\mch,J)\to(\mch^\p,J^\p)}.$$
\item $\Psi_{(\mch,J)\to(\mch,J)}\simeq {\rm id}_{(CF^-(\mch),\partial_{J})}.$
\eenu
Moreover, similar results hold for $CF^\infty$ and $CF^+$.
\ethm
\begin{conv}
If it is not mentioned, chain homotopy means $\ft[U_{\bs{w}}]$-equivariant chain homotopy.
\end{conv}
Since all chain complexes discussed above can be decomposed into spin$^c$ structures (\textit{c.f.} \cite[Section 2.6]{ozsvath2004holomorphic}), we have the following definition. 
\bdefn
Suppose $(Y,\bs{w})$ is a multi-pointed 3-manifold and $\mathfrak{s}\in{\rm Spin}^c(Y)$. For $\circ\in\{\infty,+,-\}$, define $CF^\circ(Y,\bs{w},\mathfrak{s})$ to be the transitive system of chain complexes with canonical maps from Theorem \ref{thm: naturality hf}, with respect to $\mathfrak{s}$, and define $HF^\circ(Y,\bs{w},\mathfrak{s})$ to be the induced transitive system of homology groups.
\edefn

For later use, we also define the completions of the chain complexes.
\bdefn
Let $\ft[[\uw]]$ be the ring of formal power series of $\uw$. For $\circ\in\{\infty,+,-\}$, define$${\bf CF}^\circ(Y,\bs{w},\mathfrak{s})\deq  CF^\circ(Y,\bs{w},\mathfrak{s})\ot_{\ft[\uw]}\ft[[\uw]].$$Let ${\bf HF}^\circ(Y,\bs{w},\mathfrak{s})$ be the induced homology groups.
\edefn
\begin{conv}
When omitting the module structure, we have ${\bf CF}^+(Y,\bs{w},\mathfrak{s})=CF^+(Y,\bs{w},\mathfrak{s})$. Hence we do not distinguish them.
\end{conv}
The advantage of the completions is that we have the following proposition.
\bprop[{\cite[Section 2]{Manolescu2017}, see also \cite[Lemma 2.3]{ozsvath2004symplectic}}]
If $(Y,\bs{w})$ is a multi-pointed 3-manifold and $\mathfrak{s}\in{\rm Spin}^c(Y)$ on each component is nontorsion, then $\hfi(Y,\bs{w},\mathfrak{s})=0.$
\eprop
Then the boundary map in the following long exact sequence induces a canonical isomorphism between $\hfm(Y,\bs{w},\mathfrak{s})$ and $HF^+(Y,\bs{w},\mathfrak{s})$ for any nontorsion spin$^c$ structure $\mathfrak{s}$.
\bprop\label{prop: long exact sequence}
From the short exact sequence$$0\to \cfm(Y,\bs{w},\mathfrak{s})\to\cfi(Y,\bs{w},\mathfrak{s})\to CF^+(Y,\bs{w},\mathfrak{s})\to 0,$$we have a long exact sequence
\begin{equation*}\label{eq: long exact sequence}
    \cdots\to\hfm(Y,\bs{w},\mathfrak{s})\to \hfi(Y,\bs{w},\mathfrak{s})\to HF^+(Y,\bs{w},\mathfrak{s})\to\cdots
\end{equation*}We also have a long exact sequence for $HF^-$, $HF^\infty$, and $HF^+$.
\eprop
\bdefn
Suppose $(Y,\bs{w})$ is a multi-pointed 3-manifold and $\mathfrak{s}\in{\rm Spin}^c(Y)$ is a nontorsion spin$^c$ structure. We write $$HF(Y,\bs{w},\mathfrak{s})= HF_{\rm red}(Y,\bs{w},\mathfrak{s})\deq HF^+(Y,\bs{w},\mathfrak{s})\cong \hfm(Y,\bs{w},\mathfrak{s}).$$
\edefn

\subsection{Cobordism maps for restricted graph cobordisms}\label{subsec:Cobordism maps for restricted graph cobordisms}\quad

\bthm[{\cite[Theorem A]{Zemke2019}}]
Suppose $(W,\Ga): (Y_0,\bs{w}_0) \to (Y_1, \bs{w}_1)$ is a ribbon graph cobordism and $\mathfrak{s}\in{\rm Spin}^c(W)$. Then there are two chain maps
$$F^A_{W,\Ga,\mathfrak{s}},F^B_{W,\Ga,\mathfrak{s}}:CF^-(Y_0,\bs{w}_0,\mathfrak{s}|_{Y_0})\to CF^-(Y_1,\bs{w}_1,\mathfrak{s}|_{Y_1}),$$
which are diffeomorphism invariants of $(W,\Ga)$, up to $\ft[\uw]$-equivariant chain homotopy.
\ethm
\bprop[{\cite[Theorem C]{Zemke2019}}]\label{prop: composition}
Suppose that $(W, \Ga)$ is a ribbon graph cobordism which decomposes as a composition
$(W,\Ga) = (W_2,\Ga_2) \cup (W_1,\Ga_1)$. If $\mathfrak{s}_1$ and $\mathfrak{s}_2$ are spin$^c$ structures on $W_1$ and $W_2$, respectively, then
$$F^A_{W_2,\Ga_2,\mathfrak{s}_2}\circ F^A_{W_1,\Ga_1,\mathfrak{s}_1}\simeq\sum_{\substack{\mathfrak{s}\in {\rm Spin}^c(W)\\\mathfrak{s}|_{W_2}=\mathfrak{s}_2\\\mathfrak{s}|_{W_1}=\mathfrak{s}_1}}F^A_{W,\Ga,\mathfrak{s}}.$$
A similar relation holds for $F^B_{W,\Ga,\mathfrak{s}}$.
\eprop

Since we will only consider restricted graph cobordisms, the map $F^A_{W,\Ga,\mathfrak{s}}$ is chain homotopic to $F^B_{W,\Ga,\mathfrak{s}}$. Hence we write $CF^-(W,\Ga,\mathfrak{s})$ for the chain map and $HF^-(W,\Ga,\mathfrak{s})$ for the induced map on the homology group. If $\Ga$ and $\mathfrak{s}$ are specified, we write $CF^-(W)$ and $HF^-(W)$ for simplicity, respectively. The chain maps on $CF^\infty,CF^+,\cfm,\cfi$ are obtained by tensoring with the identity maps, respectively. We use similar notations for these chain maps and the induced maps on homology groups. All maps are called \textbf{cobordism maps}.

\begin{figure}[ht]
\centering
\includegraphics[width=0.25\textwidth]{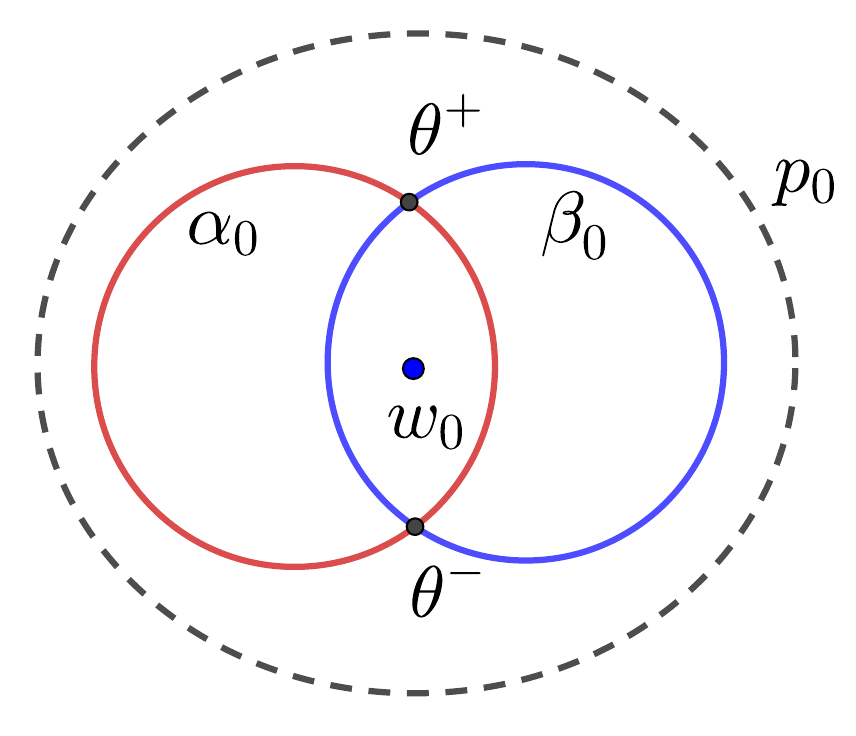}
\caption{Free-stabilization in a small disk $D$.}
\label{freestab}
\end{figure}

For $CF^-$, the cobordism map is defined by the composition of the following maps.
\begin{itemize}
    \item For 4-dimensional 1-, 2-, and 3-handle attachments away from the basepoints, we use the maps defined by Ozsv{\'{a}}th and Szab{\'{o}} \cite{Ozsvath2006}.
\item For 4-dimensional 0- and 4-handle attachments, or equivalently adding and removing a copy of $S^3$ with a single basepoint, respectively, we use the maps defined by the canonical isomorphism from the tensor product with $CF^-(S^3,w_0)\cong \ft[U_0]$.
\item For a ribbon graph cobordism $(Y\times[0,1],\Ga)$, we project the graph into $Y$ and use the \textbf{graph action map} defined in \cite[Section 7]{Zemke2019}.
\end{itemize}
\brem
For 4-dimensional 1-, 2-, and 3-handle attachments, Ozsv{\'{a}}th and Szab{\'{o}}'s original construction was for connected cobordisms between connected 3-manifolds. Zemke \cite[Section 8]{Zemke2019} extended the construction to cobordisms between possibly disconnected 3-manifolds. For 4-dimensional 0- and 4-handle attachments, the isomorphism is indeed $$CF^-(Y\sqcup S^3,\bs{w}\cup \{w_0\})\cong CF^-(Y,\bs{w})\otimes_{\ft} CF^-(S^3,w_0)\cong CF^-(Y,\bs{w})\otimes_{\ft} \ft[U_0].$$
\erem
The graph action map is obtained by the composition of maps associated to elementary graphs. The construction involves \textbf{free-stabilization maps} $S^{\pm}_w$ \cite[Section 6]{Zemke2019} and \textbf{relative homology maps} $A_\lambda$ \cite[Section 5]{Zemke2019}, where $S^{\pm}_w$ correspond to adding or removing a basepoint $w$ and $A_\lambda$ correspond to a path $\lambda$ between two basepoints. When considering restricted graph cobordisms, we only need maps associated to 1-, 2-, 3-handle attachments.
\bdefn\label{defn: freestab}
Suppose $\mch=(\Sigma,\al,\be,\bs{w})$ is a multi-pointed Heegaard diagram for a multi-pointed 3-manifold $(Y,\bs{w})$. Let $D\subset \Sigma\backslash(\al\cup\be)$ be a small disk containing a new basepoint $w_0\in\Sigma\backslash(\al\cup\be)$. Let $\al_0$ and $\be_0$ be two simple closed curves on $\Sigma$ bounding a disk containing $w_0$ and $|\al_0\cap \be_0|=2$. Suppose $\theta^+$ and $\theta^-$ are the higher and the lower graded intersection points, respectively. See Figure \ref{freestab}. Consider the Heegaard diagram $\mch^\p=(\Sigma,\al\cup\{\al_0\},\be\cup\{\be_0\},\bs{w}\cup\{w_0\})$, where $\al_0$ and $\be_0$ are in the region of a basepoint $z\in \bs{w}$.

\begin{figure}[ht]
\centering
\includegraphics[width=0.7\textwidth]{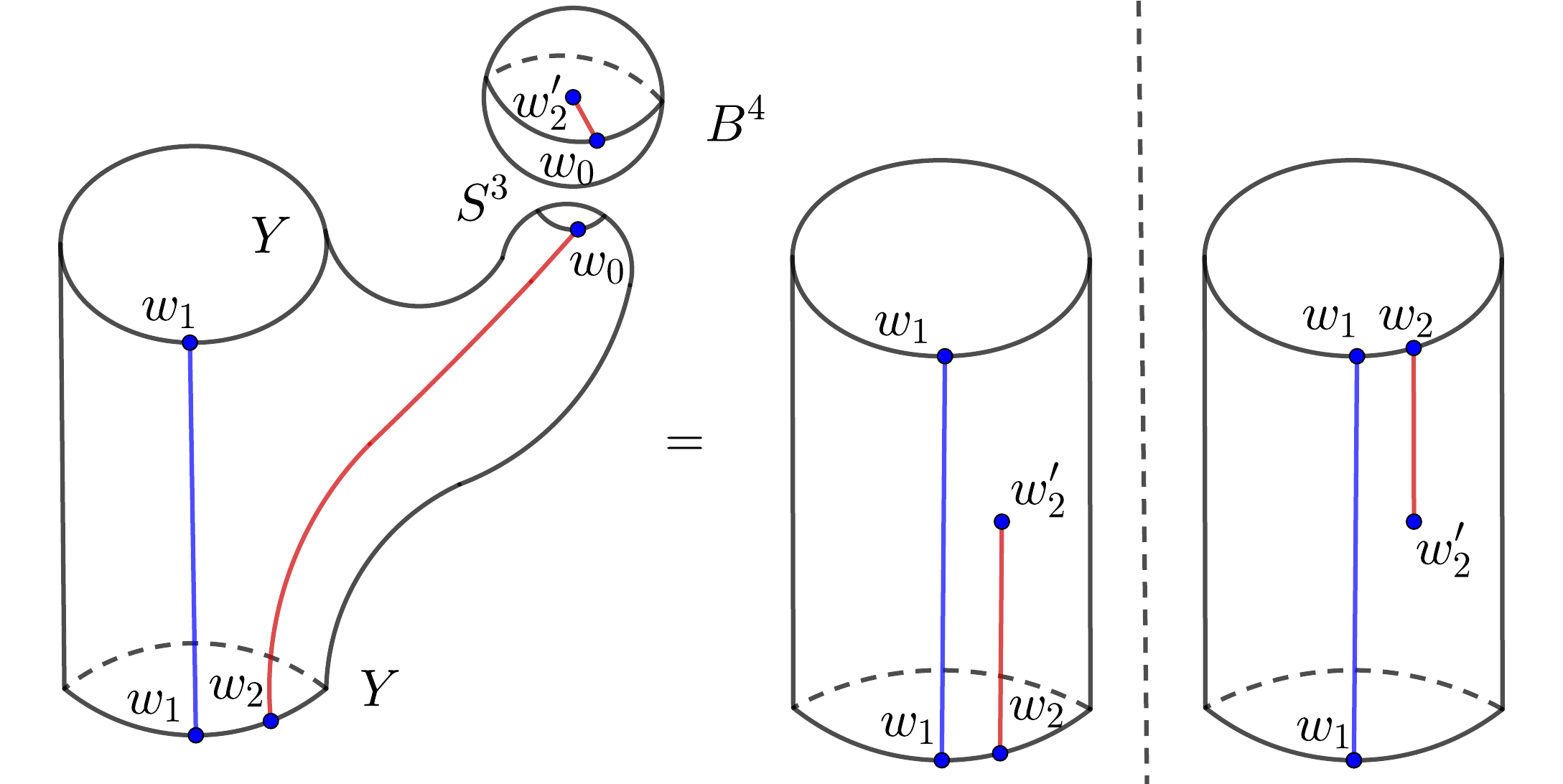}
\caption{Ribbon graph cobordisms related to free-stabilization maps.}
\label{excision31}
\end{figure}

For appropriately chosen almost complex structures, define the \textbf{free-stabilization} maps $S^{\pm}_{w_0}$ by$$S^+_{w_0}(\bs{x})=\bs{x}\times \theta^+,$$
$$S^-_{w_0}(\bs{x}\times \theta^-)=\bs{x}\aand S_{w_0}^-(\bs{x}\times \theta^+)=0.$$
\edefn
\brem\label{rem: connected sum}
If we collapse $\partial D$ to a point $p_0$, we obtain a doubly-pointed diagram on $S^2$ with two curves. Hence $\mch^\p$ can be considered as the connected sum of $\mch$ and $(S^2,\al_0,\be_0,\{w_0,p_0\})$ at the basepoint $z$ in $\mch$ and the basepoint $p_0$ (\textit{c.f.} \cite[Section 6.1]{ozsvath2008multivariable}).
\erem
\bprop[{\cite[Section 6 and Lemma 8.13]{Zemke2019}}]\label{prop: stabcommute}
The maps $S^{\pm}_{w_0}$ in Definition \ref{defn: freestab} determine well-defined chain maps on the level of transitive systems of chain complexes$$S_{w_0}^+:CF^-(Y,\bs{w})\to CF^-(Y,\bs{w}\cup\{w_0\}),$$$$S_{w_0}^-:CF^-(Y,\bs{w}\cup\{w_0\})\to CF^-(Y,\bs{w}).$$Moreover, they have the following properites.
\benu
\item The maps $S^{\pm}_{w_0}$ commute with maps associated to 1-, 2-, and 3-handle attachments.
\item For $\circ_1,\circ_2\in\{+,-\},$ we have $S^{\circ_1}_{w_1}S^{\circ_2}_{w_2}\simeq S^{\circ_2}_{w_2}S^{\circ_1}_{w_1}.$
\eenu
\eprop
\brem\label{rem: freestab}
The free-stabilization maps can be regarded as ribbon graph cobordisms with $W=Y\times [0,1]$. The graphs are shown in Figure \ref{excision31}. Alternatively, we can regard them as compositions of maps associated to handle attachments. The map $S_{w_2}^+$ is obtained by first attaching a 0-handle with an arc whose one endpoint is on the boundary, and the other is in the interior, and then attaching a 1-handle away from basepoints; see the left of Figure \ref{excision31}. The map $S_{w_2}^-$ is obtained by first attaching a 3-handle and then a 4-handle with an arc similarly; see the right of Figure \ref{excision31}.
\erem
\begin{conv}
All illustrations of cobordisms are from top to bottom.
\end{conv}
We can calculate the effect of free-stabilization maps on the homology explicitly.
\bprop[{\cite[Proposition 6.5]{ozsvath2008multivariable}}]\label{prop: freestab}Consider the construction in Definition \ref{defn: freestab}. For suitable choices of almost complex structures, the chain complex $CF^-(\mch^\p)$ is identified with the mapping cone of the following map
$$CF^-(\mch)\otimes_{\ft} \ft[U_0]\langle\theta^-\rangle\xlongrightarrow{U_0-U_1} CF^-(\mch)\otimes_{\ft} \ft[U_0]\langle\theta^+\rangle,$$where $U_1$ corresponds to the basepoint in the original diagram $\mch$ for the connected sum construction in Remark \ref{rem: connected sum}.
\eprop
\bcor\label{cor: free iso}
If $U_0\neq U_1$ in Proposition \ref{prop: freestab}, \textit{i.e.} the colors of corresponding basepoints are different (\textit{c.f.} Remark \ref{rem: color}), then the map $S^+_{w_0}$ induces isomorphisms on $HF^\circ$ and ${\bf HF}^\circ$ for $\circ\in\{\infty,+,-\}$, and the map $S^-_{w_0}$ induces zero maps on all versions of Heegaard Floer homology.
\ecor
\bpf
The arguments for $\circ\in\{\infty,-\}$ follows directly from Definition \ref{defn: freestab}, Proposition \ref{prop: freestab}, and definitions of Heegaard Floer homology groups. For $\circ=+$, note that the free-stabilization maps are compatible with the long exact sequence in Proposition \ref{prop: long exact sequence}. Hence the behaviors of maps for $\circ\in\{\infty,-\}$ imply the behavior for $\circ=+$.
\epf
The following proposition implies the choice of the basepoints is not important.
\bprop[{\cite[{Corollary 14.19 and Corollary F}]{Zemke2019}}]\label{prop: indep of path}
Suppose $(Y,\bs{w})$ is a multi-pointed 3-manifold and $w_1\in\bs{w}$. Then the $\pi_1(Y,w_1)$ action on $HF^-(Y,\bs{w})$ is always the identity map.

Suppose $(Y_1,\bs{w}_1)$ and $(Y_2,\bs{w}_2)$ are two multi-pointed 3-manifolds with $|\bs{w}_1|=|\bs{w}_2|$. Suppose $W$ is a cobordism from $Y_1$ to $Y_2$ such that the boundary of each component of $W$ consists one component of $-Y_1$ and one component of $Y_2$. Suppose $\Ga\subset W$ is a collection of paths connecting $\bs{w}_1$ and $\bs{w}_2$. Then the cobordism map $HF^-(W,\Ga)$ is independent of the choice of $\Ga$. Moreover, if $W=Y\times I$, then $HF^-(W,\Ga)$ is an isomorphism.

Similar results also hold for $HF^\infty,HF^+,\hfm,\hfi$.
\eprop
From Corollary \ref{cor: free iso} and Proposition \ref{prop: indep of path}, we can define a transitive system of groups based on different choices of basepoints.
\bdefn\label{defn: naturality of basepoints}
Suppose $Y$ is a closed, oriented 3-manifold and $\bs{w}_1,\bs{w}_2\subset Y$ are two collections of basepoints in $Y$. Let $\bs{w}_1^\p=\bs{w}_1\backslash \bs{w}_2$ and $\bs{w}_2^\p=\bs{w}_2\backslash \bs{w}_1$. For $\circ\in\{\infty,+,-\}$, define \textbf{transition maps} associated to $(\bs{w}_1,\bs{w}_2)$ as$$
\Psi_{\bs{w}_1\to\bs{w}_2}^\circ\deq \prod_{w\in\bs{w}_1^\p}(S^{+}_w)^{-1}\circ\prod_{w\in\bs{w}_2^\p}S^{+}_w\quad\text{on }HF^\circ\aand{\bf HF}^\circ$$where the products mean compositions. The order of maps is not important by the following lemma.
\edefn
\blem\label{lem: transitive 1}
Suppose $Y$ is a closed, oriented 3-manifold and $\bs{w}_1,\bs{w}_2,\bs{w}_3\subset Y$ are three collections of basepoints in $Y$. Suppose $w$ is a basepoint in $Y$ that is not in $\bs{w}_i$ for $i=1,2$. Then the following holds for transition maps.
\benu
\item $\Psi^\circ_{\bs{w}_i\to\bs{w}_j}$ is well-defined for $i,j\in\{1,2,3\}$, \textit{i.e.}, the composition is independent of the order of maps.
\item $\Psi^\circ_{\bs{w}_i\to\bs{w}_j}$ is an isomorphism for $i,j\in\{1,2,3\}$.
\item $\Psi^\circ_{\bs{w}_i\to\bs{w}_i}={\rm id}$ for $i=1,2,3.$
\item $\Psi^\circ_{\bs{w}_2\to\bs{w}_3}\circ \Psi^\circ_{\bs{w}_1\to\bs{w}_2}= \Psi^\circ_{\bs{w}_1\to\bs{w}_3}$.
\item $\Psi^\circ_{\bs{w}_1\cup \{w\}\to \bs{w}_2\cup \{w\}}\circ S^+_w=S^+_w\circ \Psi^\circ _{\bs{w}_1\to \bs{w}_2}$.
\item $\Psi^\circ_{\bs{w}_1\to \bs{w}_2}\circ S^-_w=S^-_w\circ \Psi^\circ _{\bs{w}_1\cup \{w\}\to \bs{w}_2\cup \{w\}}$.
\eenu
\elem
\bpf
Terms (1), (4), (5) and (6) follow from term (2) of Proposition \ref{prop: stabcommute}. Note that maps in terms (5) are both isomorphisms and the maps in term (6) are both zero maps. Term (3) is trivial from the definition. Term (2) follows from Corollary \ref{cor: free iso}.
\epf
\blem\label{lem: transitive 2}
Suppose $Y_1$ and $Y_2$ are closed, oriented 3-manifolds and $\bs{w}_1,\bs{w}_2\subset Y_1,\bs{w}_3,\bs{w}_4\subset Y_2$ are collections of basepoints. Suppose $W$ is a cobordism from $Y_1$ to $Y_2$ that is obtained from $Y_1\times I$ by attaching 4-dimensional 1-, 2-, 3-handles away from all basepoints. Let $\Ga_1=\bs{w}_1\times I$ be the induced graph in $W$ and suppose $\bs{w}_3$ is the image of $\bs{w}_1\times \{1\}$. The cobordism $W$ can also be obtained from $-Y_2\times I$ by attaching handles away from basepoints and let $\Ga_2=\bs{w}_4\times I$. Suppose the image of $\bs{w}_4$ is $\bs{w}_2$. Then we have a commutative diagram
\begin{equation*}
\xymatrix@R=6ex{
HF^-(Y_1,\bs{w}_1)\ar[rr]^{HF^-(W,\Ga_1)}\ar[d]^{\Psi^-_{\bs{w}_1\to\bs{w}_2}}&&HF^-(Y_2,\bs{w}_3)\ar[d]^{\Psi^-_{\bs{w}_3\to\bs{w}_4}}\\
HF^-(Y_1,\bs{w}_2)\ar[rr]^{HF^-(W,\Ga_2)}&&HF^-(Y_2,\bs{w}_4)
}    
\end{equation*}
Similar commutative diagrams hold for $\hfm$ and $HF^+$.
\elem
\bpf
This follows from term (1) of Proposition \ref{prop: stabcommute}.
\epf
\bthm\label{thm: transitive system}
Suppose $Y$ is a closed, oriented 3-manifold. Then groups $HF^-(Y,\bs{w})$ for all $\bs{w}\subset Y$ and transition maps $\Psi_{\bs{w}_1\to\bs{w}_2}^-$ for all $\bs{w}_1,\bs{w}_2\subset Y$ form a transitive system, which is denoted by $HF^-(Y)$. Moreover, suppose $(W,\Ga)$ is a restricted graph cobordism from $(Y_1,\bs{w}_1)$ to $(Y_2,\bs{w}_2)$. Then $HF^-(W,\Ga)$ induces a well-defined map from $HF^-(Y_1)$ to $HF^-(Y_2)$, which is independent of the choice of the restricted graph $\Ga$ and denoted by $HF^-(W)$.

Similar arguments hold for infinity and plus versions of Heegaard Floer homology groups.
\ethm
\bpf
The well-definedness of $HF^-(Y)$ and $HF^-(W,\Ga)$ follows from Lemma \ref{lem: transitive 1} and Lemma \ref{lem: transitive 2}. Note that the restricted graph cobordism is a composition of maps associated to 1-, 2-, 3-handle attachments. Then the independence of $\Ga$ follows from the functoriality of the map associated to a ribbon graph cobordism. The proofs for infinity and plus versions of Heegaard Floer homology groups are similar.
\epf
\brem\label{rem: hfred}
Groups and maps in Theorem \ref{thm: transitive system} also split into spin$^c$ structures. Suppose $\mathfrak{s}\in{\rm Spin}^c(W)$ is a nontorsion spin$^c$ structure which restricts to nontorsion spin$^c$ structure $\mathfrak{s}_i$ on $Y_i$ for $i=1,2$. Then $\hfm(Y_i,\mathfrak{s}_i)$ and $HF^+(Y_i,\mathfrak{s}_i)$ are canonically identified by the boundary map in Proposition \ref{prop: long exact sequence}. Moreover, the maps $\hfm(W,\mathfrak{s})$ and $HF^+(W,\mathfrak{s})$ are the same under this identification. We write the map as $HF(W,\mathfrak{s}).$
\erem

\subsection{Floer's excision theorem}\label{subsec: Floer}\quad

Note that the proofs of Theorem \ref{thm: canonical maps} and Theorem \ref{thm: grading is well defined, appendix} (\textit{c.f.} \cite{baldwin2015naturality,li2019direct}) both involve Floer's excision theorem in an essential way. In this subsection, we follow Kronheimer and Mrowka's idea in \cite[Section 3]{kronheimer2010knots} to prove an excision theorem for Heegaard Floer theory. The proof in \cite[Section 3]{kronheimer2010knots} depends essentially on the TQFT properties and Axiom (A1), so it works for a general TQFT satisfying Axiom (A1). Though for Heegaard Floer theory, we need to modify the proof to fit the settings of multi-basepoints 3-manifolds and ribbon graph cobordisms.

Let $Y$ be a closed, oriented 3-manifold, of either one or two components. In the latter case, let $Y_1$ and $Y_2$ be two components of $Y$. Let $\Sigma_1$ and $\Sigma_2$ be two closed, connected, oriented surfaces in $Y$ with $g(\Sigma_1)=g(\Sigma_2)$. If $Y$ has two components, suppose $\Sigma_i$ is a non-separating surface in $Y_i$ for $i=1,2$. If $Y$ is connected, suppose $\Sigma_1$ and $\Sigma_2$ represent independent homology classes. In either
case, let $F=\Sigma_1\cup\Sigma_2$. Let $h$ be an orientation-preserving diffeomorphism from $\Sigma_1$ to $\Sigma_2$. 

\begin{figure}[ht]
\centering
\includegraphics[width=0.4\textwidth]{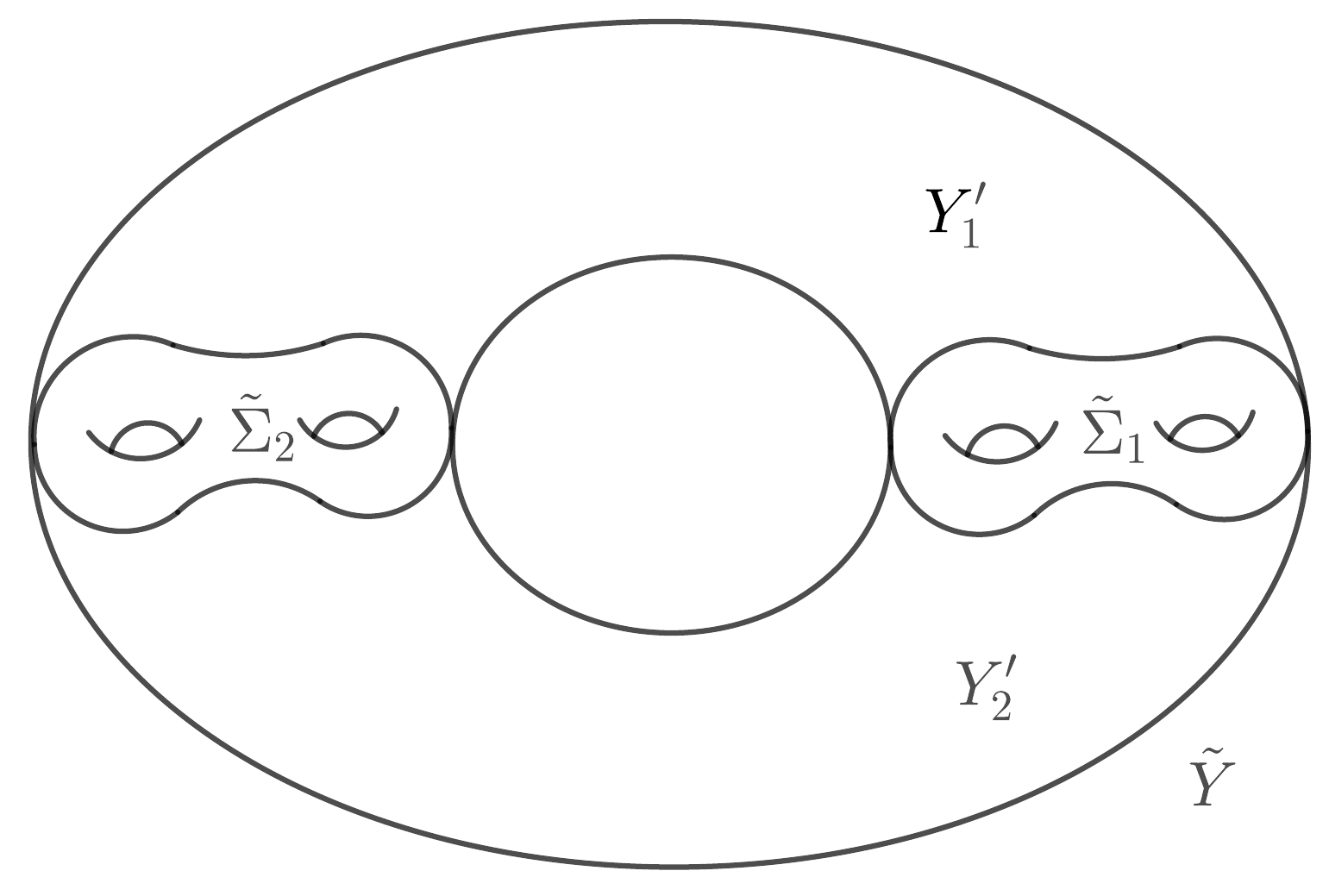}
\caption{Construction of $\widetilde{Y}$.}
\label{excision0}
\end{figure}

We construct a new manifold $\widetilde{Y}$ as follows. Let $Y^\p$ be obtained from $Y$ by cutting along $\Sigma$. Then$$\partial Y^\p=\Sigma_1\cup(-\Sigma_1)\cup\Sigma_2\cup(-\Sigma_2).$$
If $Y$ has two components, then we have $Y^\p=Y_1^\p\cup Y_2^\p$, where $Y_i^\p$ is obtained from $Y_i$ by cutting along $\Sigma_i$ for $i=1,2$. Let $\widetilde{Y}$ be obtained from $Y^\p$ by gluing the boundary component $\Sigma_1$ to the boundary component
$-\Sigma_2$ and gluing $\Sigma_2$ to $-\Sigma_1$, using the diffeomorphism of $h$ in both cases; see Figure \ref{excision0} for the case that $Y$ has two components. 

In either case, $\widetilde{Y}$ is connected. Let $\widetilde{\Sigma}_i$ be the image of $\Sigma_i$ in $\widetilde{Y}$ for $i=1,2$ and let $\widetilde{F}=\widetilde{\Sigma}_1\cup\widetilde{\Sigma}_2$.
\bdefn
Suppose $Y$ is a closed, oriented 3-manifold and $F\subset Y$ is a closed, oriented surface. Let $F_i$ for $i=1,\dots,m$ be the components of $F$. Suppose further that $g(F_i)\ge 2$ and any component of $Y$ contains at least one component of $F$. Let ${\rm Spin}^c(Y|F)$ denote the set of spin$^c$ structures $\mathfrak{s}\in{\rm Spin}^c(Y)$ satisfying \begin{equation}\label{eq: top class}
    \langle c_1(\mathfrak{s}),F_i\rangle=2g(F_i)-2 \text{ for any }F_i.
\end{equation}Define $$HF(Y|F)\deq \bigoplus_{\mathfrak{s}\in{\rm Spin}^c(Y|F)}HF(Y,\mathfrak{s}).$$

Suppose $(W,\Ga)$ is a restricted graph cobordism and $G\subset W$ is a closed, oriented surface. Let $G_i$ for $i=1,\dots,n$ be components of $G$. Suppose further that $g(G_i)\ge 2$ and any component of $W$ contains at least one component of $G$. Let ${\rm Spin}^c(W|G)$ denote the set of  spin$^c$ structures $\mathfrak{s}\in{\rm Spin}^c(W)$ satisfying similar conditions in (\ref{eq: top class}) by replacing $F_i$ by $G_i$. Define$$HF^-(W,\Ga|G)\deq \sum_{\mathfrak{s}\in{\rm Spin}^c(W|G)}HF^-(W,\Ga,\mathfrak{s}).$$Let $HF^+(W,\Ga|G)$, $\hfm(W,\Ga|G)$ and $HF(W,\Ga|G)$ be defined similarly. We also denote the corresponding map on the chain level by replacing $HF$ by $CF$.
\edefn
\brem
All spin$^c$ structures in $\spin(Y|F)$ are nontorsion, so $HF(Y,\mathfrak{s})$ is well-defined.
\erem
The following is the main theorem of this subsection.
\bthm[Floer's excision theorem]\label{thm: floer excision}
Consider $Y$ and $\widetilde{Y}$ constructed as above. If $g(\Sigma_1)=g(\Sigma_2)\ge 2$, then there is an isomorphism$$HF(Y|F)\cong HF(\widetilde{Y}|\widetilde{F}).$$Moreover, this isomorphism and its inverse are induced by restricted graph cobordisms.
\ethm
Before proving the main theorem, we introduce some lemmas analogous to results in monopole theory (\textit{c.f.} \cite[Lemma 2.2, Proposition 2.5 and Lemma 4.7]{kronheimer2010knots})
\blem[{\cite[Theorem 16 and Corollary 17]{lekili2013heegaard}, see also \cite[Theorem 5.2]{ozsvath2004symplectic}}]\label{lem: HF1d}
Let $Y \to S^1$ be a fibred 3-manifold whose fibre $F$ is a closed, connected, oriented surface with $g=g(F)\ge 2$. Then $\cfm(Y|F)$ is chain homotopic to the chain complex \begin{equation}\label{eq: 1d chain}
    0\to\ft[[U_0]]\langle x\rangle \xra{U_0}\ft[[U_0]]\langle y\rangle\to 0.
\end{equation}Moreover, there is a unique $\mathfrak{s}_0\in\spin(Y|R)$ so that $HF(Y,\mathfrak{s}_0)\neq 0$ and we have$$HF(Y|F)=HF(Y,\mathfrak{s}_0)\cong \ft.$$
\elem
\brem
Indeed, for $Y$ in Lemma \ref{lem: HF1d}, we can construct a \textit{weakly} admissible Heegaard diagram $\mch$ for the singly-pointed 3-manifold $(Y,w)$ so that $\cfm(\mch,\mathfrak{s}_0)$ is generated by $8g$ generators $\bs{x}_1,\dots,\bs{x}_{4g},\bs{y}_1,\dots,\bs{y}_{4g}$ and $$\partial \bs{x}_1=U_0\bs{y}_1,\partial \bs{x}_j=\bs{y}_j,\aand\partial \bs{y}_k=0 \text{ for }j>1,k\ge 1.$$The reason to use $\cfm$ rather than $CF^-$ is because the computation of $CF^-$ is based on \textit{strongly} admissible Heegaard diagram.
\erem
\blem\label{lem: nontrivial map}
Suppose $Y=\Sigma\times S^1$ such that $\Sigma=\Sigma\times\{1\}\subset Y$ is a closed, connected, oriented surface with $g(\Sigma)\ge 2$. Suppose $w_0\in S^3$ and $w\in Y$ are basepoints. Let $W$ be obtained from $\Sigma\times D^2$ by removing a 4-ball, considered as a cobordism from $S^3$ to $Y$. Let $\Ga\subset W$ be any path connecting $w_0$ to $w$. Then the map \begin{equation}\label{eq: nontrivial}
    \hfm(W,\Ga|\Sigma):\ft[[U_0]]\cong \hfm(S^3,w_0)\to HF(Y,w|\Sigma)\cong\ft
\end{equation}is nonzero. 

\elem
\begin{figure}[ht]
\centering
\includegraphics[width=0.3\textwidth]{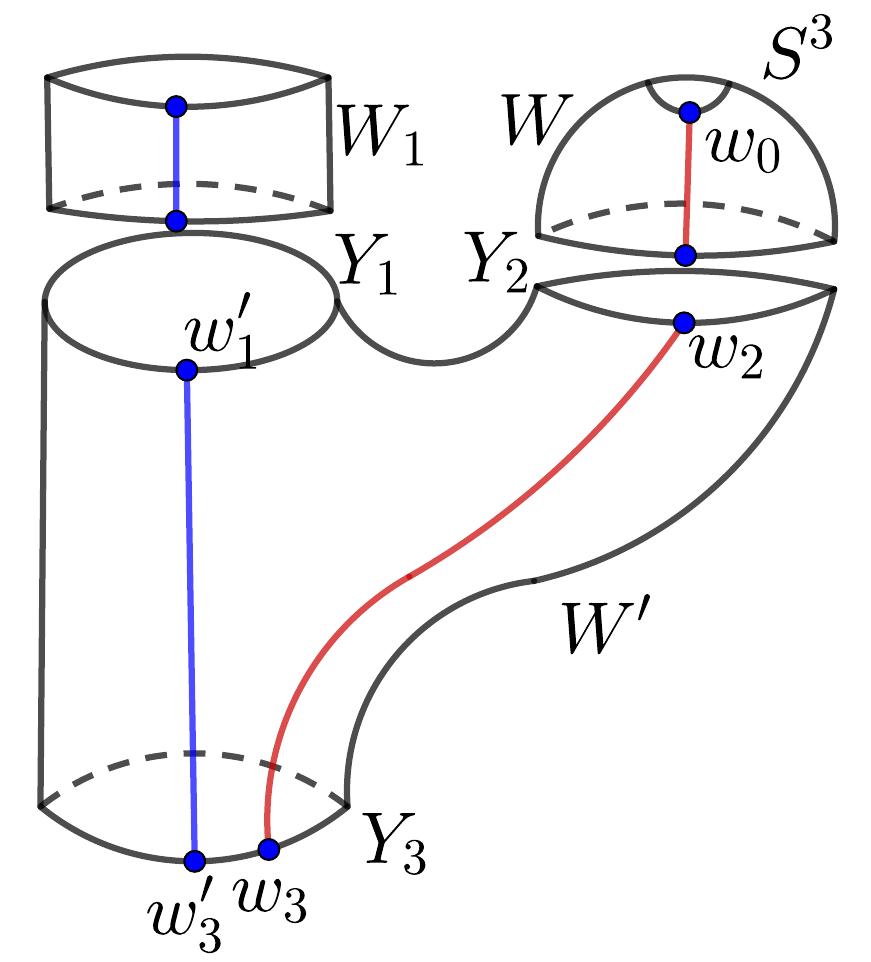}
\caption{Nontrivial cobordism map from composition.}
\label{excision21}
\end{figure}
\bpf
Suppose $P$ is 2-dimensional pair of pants as shown in Figure \ref{excision21}. Consider $W^\p=\Sigma\times P$ as a cobordism from $Y_1\sqcup Y_2$ to $Y_3$, where $Y_i\cong Y$ for $i=1,2,3$. Suppose $w^\p$ is another basepoint in $Y$. Let $w_i$ and $w_i^\p$ be the images of $w$ and $w^\p$ in $Y_i$ for $i=1,2,3$. Let $\Ga^\p\subset W^\p$ be a collection of two paths $\ga_1$ and $\ga_2$, where $\ga_1$ connects $w^\p_1$ to $w^\p_3$ and $\ga_2$ connects $w_2$ to $w_3$.

Let $(W_1,\Ga_1)=(Y_1\times I,w^\p\times I)$ be the product cobordism. Suppose $\Sigma_i\subset Y_i$ is the image of $\Sigma\subset Y$ for $i=1,2,3$. Consider the composition of the cobordism maps
$$\begin{aligned}
\hfm(W^\p,\Ga^\p|\Sigma_1\cup \Sigma_2\cup \Sigma_3)&\circ \hfm(W_1\sqcup W,\Ga_1\sqcup\Ga|\Sigma_1\cup \Sigma_2):\\HF(Y_1,w_1^{\prime}|\Sigma_1)&\ot_\ft \hfm( S^3,w_0)\to HF(Y_3,\{w_3,w_3^\p\}|\Sigma_3).
\end{aligned}$$After filling the $S^3$ component by a 4-ball, or equivalently composing it with the map associated to a 0-handle attachment, we obtain the free-stabilization map $S_w^+$ (\textit{c.f.} Remark \ref{rem: freestab}). By Corollary \ref{cor: free iso}, the resulting map is an isomorphism$$HF(Y_1,w_1^{\prime}|\Sigma_1)\cong HF(Y_3,\{w_3,w_3^\p\}|\Sigma_3).$$Since $$\hfm(W_1\sqcup W,\Ga_1\sqcup\Ga|\Sigma_1\cup \Sigma_2)=\hfm(W_1,\Ga_1|\Sigma_1)\ot_\ft \hfm(W,\Ga|\Sigma_2),$$and $\hfm(W_1|\Sigma_1)$ is the identity map, we know $\hfm(W|\Sigma_2)$ is nonzero.

\begin{figure}[ht]
\centering
\includegraphics[width=0.95\textwidth]{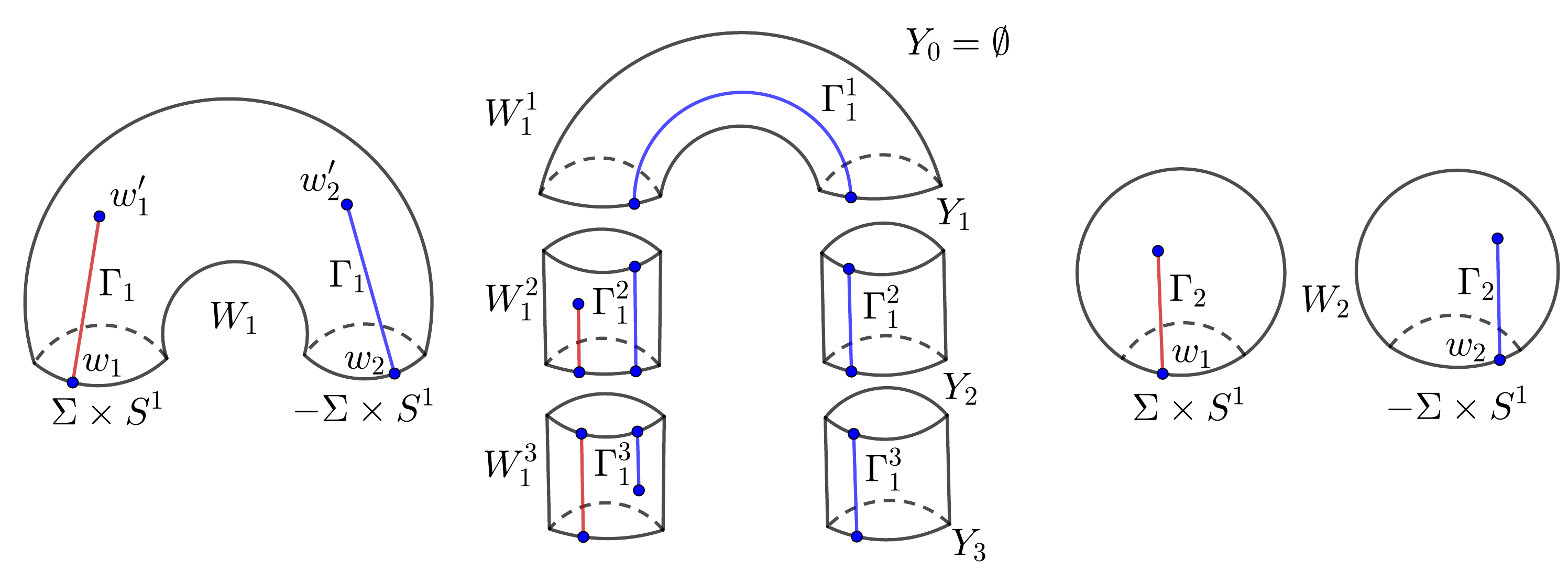}
\caption{Ribbon graph cobordisms $(W_1,\Ga_1)$ and $(W_2,\Ga_2)$.}
\label{excision11}
\end{figure}

\epf
\bcor[]\label{cor: nontrivial map}
On the chain level of (\ref{eq: nontrivial}), the cobordism map $\cfm(W,\Ga|\Sigma)$ sends the generator of $\cfm(S^3,w_0)\cong \ft[[U_0]]$ to the generator of second copy of $\ft[[U_0]]$ in (\ref{eq: 1d chain}).
\ecor
\bpf
The map in the statement is the only $\ft[U_0]$-equivariant chain map that induces a nonzero map on the homology.
\epf
The proof of the following lemma is due to Ian Zemke.
\blem\label{lem: two cobordism equal}
Let $Y=\Sigma\times S^1$ and let $W_1\cong Y\times I$ be a cobordism from $\emptyset$ to $Y\sqcup (-Y)$. Let $w_1\in Y,w_2\in(-Y),w_1^{\prime},w_2^{\prime}\in W_1$ and let $\Ga_1\subset W_1$ consist of two paths whose enpoints are $w_i$ and $w_i^{\prime}$ for $i=1,2$, as shown in the left subfigure of Figure \ref{excision11}. Let $W_2\cong \Sigma\times D^2\sqcup (-\Sigma\times D^2)$ be another cobordism from $\emptyset$ to $Y\sqcup (-Y)$ and let $\Ga_2\subset W_2$ be obtained from two copies of the cobordism in Lemma \ref{lem: nontrivial map} associated to $\Sigma$ and $-\Sigma$ by filling the $S^3$ components by 4-balls (\textit{c.f.} Remark \ref{rem: freestab}), as shown in the right subfigure of Figure \ref{excision11}. Then we have\begin{equation}\label{eq: same cf}
    \cfm(W_1,\Ga_1|\Sigma\sqcup(-\Sigma))\simeq \cfm(W_2,\Ga_2|\Sigma\sqcup(-\Sigma)):\cfm(\emptyset)\to \cfm(Y\sqcup(-Y),\{w_1,w_2\}|\Sigma\sqcup(-\Sigma)).
\end{equation}
\elem
\bpf
Set $\mcr=\ft[[U_1,U_2]]$. By Remark \ref{rem: color}, we implicitly choose $w_1$ and $w_2$ to have different colors and then $$\cfm(Y\sqcup(-Y),\{w_1,w_2\}|\Sigma\sqcup(-\Sigma))\deq \cfm(Y|\Sigma)\otimes_\ft \cfm(-Y|-\Sigma).$$By Remark \ref{rem: color2}, we have $\cfm(\emptyset)=\mathcal{R}$. By TQFT property in \cite{Zemke2019}, we have a canonical chain isomorphism$$\cfm(-Y,w_2|-\Sigma)\cong \cfm(Y,w_2|\Sigma)^\vee\deq {\rm Hom}_{\mathcal{R}}(\cfm(Y,w_2|\Sigma),\mcr).$$Then by Lemma \ref{lem: HF1d}, we have \begin{equation}\label{eq: cf1}
    \begin{aligned}
    \cfm(Y\sqcup(-Y),\{w_1,w_2\}|\Sigma\sqcup(-\Sigma))\simeq \xymatrix{
\mcr\langle x\otimes y^\vee\rangle\ar[r]^{U_1}\ar[d]^{U_2}&\mcr\langle x\otimes x^\vee\rangle\ar[d]^{U_2}&\\
\mcr\langle y\otimes y^\vee\rangle\ar[r]^{U_1}&\mcr\langle y\otimes x^\vee\rangle&
}	
    \end{aligned}
\end{equation}
where $x^\vee$ and $y^\vee$ are duals of $x$ and $y$, respectively. By Corollary \ref{cor: nontrivial map}, we know $\cfm(W_2,\Ga_2|\Sigma\sqcup(-\Sigma))$ sends the generator of $\cfm(\emptyset)$ to $y\otimes x^\vee$ in (\ref{eq: cf1}). 

By Proposition \ref{prop: composition}, we compute $\cfm(W_1,\Ga_1|\Sigma\sqcup(-\Sigma))$ by decomposing $(W_1,\Ga_1)$ into three parts $(W_1^i,\Ga_1^i):(Y_{i-1},\bs{w}_{i-1})\to (Y_i,\bs{w}_i)$ for $i=1,2,3$ as shown in the middle subfigure of Figure \ref{excision11}. Note that $(Y_0,\bs{w}_0)=\emptyset$. Let $F$ be the images of $\Sigma\sqcup (-\Sigma)$.

First, we compute $\cfm(W_1^1,\Ga_1^1|F)$. Since the two basepoints in $\bs{w}_1$ have the same color (also the same as $w_2$), we have \begin{equation}\label{eq: cf2}
    \cfm(Y_1,\bs{w}_1|F)\simeq \xymatrix{
\mcr\langle x\otimes y^\vee\rangle\ar[r]^{U_2}\ar[d]^{U_2}&\mcr\langle x\otimes x^\vee\rangle\ar[d]^{U_2}\\
\mcr\langle y\otimes y^\vee\rangle\ar[r]^{U_2}&\mcr\langle y\otimes x^\vee\rangle
}
\end{equation}
From Zemke's calculation \cite[Theorem 1.7]{Zemke2018}, the cobordism map $\cfm(W_1^1,\Ga_1^1|F)$ is the canonical cotrace map, \textit{i.e.}, it sends the generator of $\cfm(\emptyset)=\mcr$ to $x\otimes x^\vee+y\otimes y^\vee$. Note that the original calculation is for $CF^-$ but it is easy to extend the result to $\cfm$.
\brem\label{rem: mcr}
Though we only have one color in $\bs{w}_1$, we use $\mcr$ rather than $\ft[[U_2]]$ in (\ref{eq: cf2}) to achieve the functoriality (\textit{c.f.} Remark \ref{rem: color2}). Thus, when applying Proposition \ref{prop: freestab} in the following computation, we do not need to add another $U$-variable.
\erem

Second, we compute $\cfm(W_1^2,\Ga_1^2|F)$. Note that the left component of $(W_1^2,\Ga_1^2)$ corresponds to the free-stabilization map $S^+_{w_1}$ and the right component is just the identity map. By Proposition \ref{prop: freestab}, the chain complex $\cfm(Y_2,\bs{w}_2|F)$ is chain homotopic to the mapping cone of\begin{equation}\label{eq: cf3}
    \left(\vcenter{\xymatrix{
\mcr\langle x\otimes y^\vee\otimes \theta^-\rangle\ar[r]^{U_2}\ar[d]^{U_2}&\mcr\langle x\otimes x^\vee\otimes \theta^-\rangle\ar[d]^{U_2}\\
\mcr\langle y\otimes y^\vee\otimes \theta^-\rangle\ar[r]^{U_2}&\mcr\langle y\otimes x^\vee\otimes \theta^-\rangle\\
}} \right)\xra{U_1-U_2}\left(\vcenter{\xymatrix{
\mcr\langle x\otimes y^\vee\otimes \theta^+\rangle\ar[r]^{U_2}\ar[d]^{U_2}&\mcr\langle x\otimes x^\vee\otimes \theta^+\rangle\ar[d]^{U_2}\\
\mcr\langle y\otimes y^\vee\otimes \theta^+\rangle\ar[r]^{U_2}&\mcr\langle y\otimes x^\vee\otimes \theta^+\rangle\\
}}\right)
\end{equation}where $u\otimes v\otimes \theta^{\pm}$ for $u\in\{x,x^\vee\},y\in\{y,y^\vee\}$ represents $(u\times \theta^{\pm})\otimes v$. Then $\cfm(W_1^2,\Ga_1^2|F)$ sends any generator $u\otimes v$ to $u\otimes v \otimes \theta^+$ in (\ref{eq: cf3}).

Third, we compute $\cfm(W_1^3,\Ga_1^3|F)$. Note that the left component of $(W_1^3,\Ga_1^3)$ corresponds to the free-stabilization map $S^-_{w_2}$ and the right component is just the identity map. Also by Proposition \ref{prop: freestab}, the chain complex $\cfm(Y_2,\bs{w}_2|F)$ is chain homotopic to the mapping cone of\begin{equation}\label{eq: cf4}
    \left(\vcenter{\xymatrix{
\mcr\langle x\otimes y^\vee\otimes \theta^-\rangle\ar[r]^{U_1}\ar[d]^{U_2}&\mcr\langle x\otimes x^\vee\otimes \theta^-\rangle\ar[d]^{U_2}\\
\mcr\langle y\otimes y^\vee\otimes \theta^-\rangle\ar[r]^{U_1}&\mcr\langle y\otimes x^\vee\otimes \theta^-\rangle\\
}} \right)\xra{U_2-U_1}\left(\vcenter{\xymatrix{
\mcr\langle x\otimes y^\vee\otimes \theta^+\rangle\ar[r]^{U_1}\ar[d]^{U_2}&\mcr\langle x\otimes x^\vee\otimes \theta^+\rangle\ar[d]^{U_2}\\
\mcr\langle y\otimes y^\vee\otimes \theta^+\rangle\ar[r]^{U_1}&\mcr\langle y\otimes x^\vee\otimes \theta^+\rangle\\
}}\right)
\end{equation}Then $\cfm(W_1^3,\Ga_1^3|F)$ sends $u\otimes v\otimes \theta^-$ to $u\otimes v$ in (\ref{eq: cf1}) and sends $u\otimes v\otimes \theta^+$ to $0$ for $u\in\{x,x^\vee\},y\in\{y,y^\vee\}$. 

To compute the composition, we need to find the explicit chain homotopy between the above two mapping cones (\ref{eq: cf3}) and (\ref{eq: cf4}), which is calculated by Zemke \cite[Theorem 14.1]{Zemke2019}. Since we only care about the image of $\cfm(\emptyset)$, we only need to calculate the image of $\ast$ map in \cite[(14.3)]{Zemke2019} (from the target in (\ref{eq: cf3}) to the source in (\ref{eq: cf4}))\begin{equation}\label{eq: star map}
    (\Psi_{\al\to \al^\p}^{\be^\p})_{U_w}^{U_z\to U_{w^\p}}\circ \left( \sum_{i,j\ge 0}U_w^iU_{w^\p}^j(\partial_{i+j+1})_{U_w,U_{w^\p}}\right)\circ(\Psi_{\al}^{\be\to\be^\p})_{U_{w^\p}}^{U_z\to U_{w}}
\end{equation}for the element \begin{equation}\label{eq: element}
    x\otimes x^\vee\otimes \theta^++y\otimes y^\vee\otimes \theta^+
\end{equation}in (\ref{eq: cf3}). In (\ref{eq: star map}), we have $z\in Y_1$ for the connected sum construction in Remark 
\ref{rem: connected sum}, $w=w_2,w^\p=w_1,U_w=U_2,U_{w^\p}=U_1$ and $\al^\p,\be^\p$ being small isotopies of $\al,\be$, respectively. The differential $\partial_k$ comes from \begin{equation}\label{eq: partial}
    \partial=\sum_{k\in\mathbb{N}}U_{z}^k\partial_k,
\end{equation}where $\partial$ is the differential in \begin{equation}\label{eq: cf5}
    \begin{aligned}
    \cfm(Y_1,\{z,w_2\}|\Sigma\sqcup(-\Sigma))\simeq \xymatrix{
\mcr\langle x\otimes y^\vee\rangle\ar[r]^{U_z}\ar[d]^{U_2}&\mcr\langle x\otimes x^\vee\rangle\ar[d]^{U_2}&\\
\mcr\langle y\otimes y^\vee\rangle\ar[r]^{U_z}&\mcr\langle y\otimes x^\vee\rangle&
}	
    \end{aligned}
\end{equation} For a map $f$, the notation $(f)^{U_z\to U_{w}}$ means we replace $U_z$ by $U_w$ in the image of $f$ and the notation $(f)_{U_w}$ means tensoring $f$ with the identity map in $\ft[U_w]$.

Since the element (\ref{eq: element}) has no $U$-power, the transition maps $(\Psi_{\al\to \al^\p}^{\be^\p})_{U_w}^{U_z\to U_{w^\p}}$ and $(\Psi_{\al}^{\be\to\be^\p})_{U_{w^\p}}^{U_z\to U_{w}}$ can be regarded as identity maps. By (\ref{eq: partial}) and (\ref{eq: cf5}), we know $\partial_k=0$ for $k\ge 1$ and $\partial _1$ sends $x\otimes x^\vee$ to $0$ and sends $y\otimes y^\vee$ to $y\otimes x^\vee$. Hence the $\ast$ map (\ref{eq: star map}) sends the element (\ref{eq: element}) to $y\otimes x^\vee\otimes \theta^-$ in (\ref{eq: cf4}).

Thus, by composing three cobordism maps and up to chain homotopy, we show that $\cfm(W_1,\Ga_1|\Sigma\sqcup(-\Sigma))$ also sends the generator of $\cfm(\emptyset)=\mcr$ to $y\otimes x^\vee$ in (\ref{eq: cf1}).
\epf
Now we start to prove the main theorem of this subsection. The basic idea is from Kronheimer and Mrowka \cite[Section 3.2]{kronheimer2010knots}, which originally came from Floer's work \cite{floer1990knot}, where he dealt with the excision theorem in instanton theory for the genus one case.
\bpf[Proof of Theorem \ref{thm: floer excision}]
{\bf Step 1}. We construct a cobordism $W$ from $\widetilde{Y}$ to $Y$ and a cobordism $\bar{W}$ from $Y$ to $\widetilde{Y}$.

Recall that $Y^\p$ is obtained from $Y$ by cutting along $\Sigma_1$ and $\Sigma_2$ and we have $$\partial Y^\p=\Sigma_1\cup(-\Sigma_1)\cup\Sigma_2\cup(-\Sigma_2).$$

\begin{figure}[ht]
\centering
\includegraphics[width=0.8\textwidth]{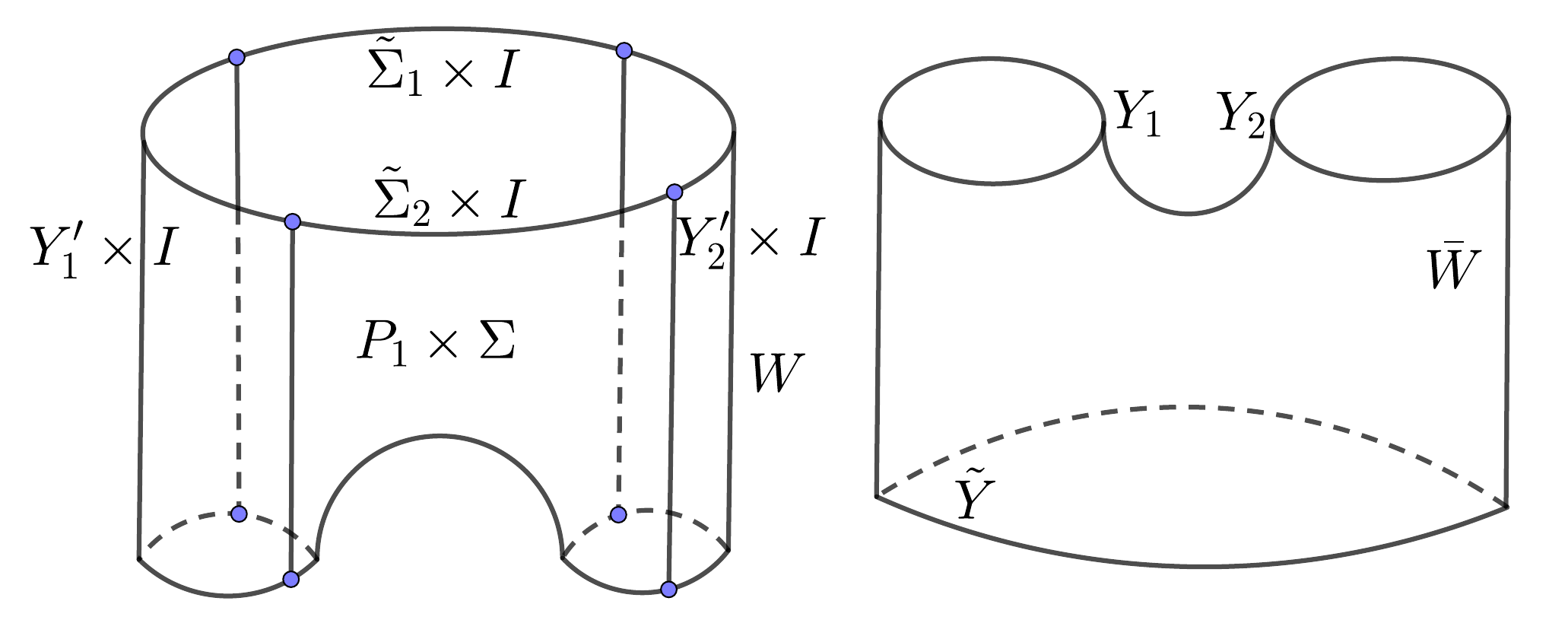}
\caption{Cobordisms $W$ and $\bar{W}$.}
\label{excision41}
\end{figure}

Suppose $P_1$ is a saddle surface, which can be regarded as a submanifold of a pair of pants with one boundary component on the top and two boundary components at the bottom; see the left subfigure of Figure \ref{excision41}. Suppose $$\partial P_1=\lambda_1\cup\lambda_2\cup\mu_1\cup\mu_2\cup\eta_{1,1}\cup\eta_{1,2}\cup\eta_{2,1}\cup\eta_{2,2},$$where $\lambda_1$ and $\lambda_2$ are two arcs in the top boundary component of the pair of pants, $\mu_1$ and $\mu_2$ are two arcs in the bottom boundary components of the pair of pants, and $\eta_{i,j}$ is the arc connecting $\lambda_i$ and $\mu_j$ for $i,j\in\{1,2\}$.

Suppose $\Sigma\cong\Sigma_1\cong\Sigma_2$. Note that we have fixed a diffeomorphism $h$ from $\Sigma_1$ to $\Sigma_2$. Suppose $h^\p$ is an orientation-preserving diffeomorphism from $\Sigma$ to $\Sigma_1$. Let $W$ be the union$$P_1\times \Sigma\cup Y^\p\times I,$$where $\eta_{1,1}\times \Sigma$ is glued to $\Sigma_1\times I$, $\eta_{2,1}\times \Sigma$ is glued to $-\Sigma_1\times I$, $\eta_{2,2}\times \Sigma$ is glued to $\Sigma_2\times I$, and $\eta_{1,2}\times \Sigma$ is glued to $-\Sigma_2\times I$, using $h^\p$ and $h\circ h^\p$, respectively. Figure \ref{excision41} illustrates the case that $Y^\p$ has two components $Y_1^\p$ and $Y_2^\p$. By the construction of $\widetilde{Y}$, the resulting manifold $W$ is a cobordism from $\widetilde{Y}$ to $Y$.

The cobordism $\bar{W}$ is constructed similarly. Let $P_2$ be another saddle surface and let $\bar{W}$ be obtained by gluing $P_2\times \Sigma$ and $Y^\p\times I$ as shown in the right subfigure of Figure \ref{excision41}.

{\bf Step 2}. For some restricted graph $\Ga_A$ and some surface $G_A$ in $W_A=\bar{W}\cup_{\widetilde{Y}}W$, we show the cobordism map $$HF(W_A,\Ga_A|G_A)\deq HF^+(W_A,\Ga_A|G_A)=\hfm(W_A,\Ga_A|G_A)$$induces the identity map on $$HF(Y|F)\deq HF^+(Y|F)\cong \hfm(Y|F).$$

\begin{figure}[ht]
\centering
\includegraphics[width=0.9\textwidth]{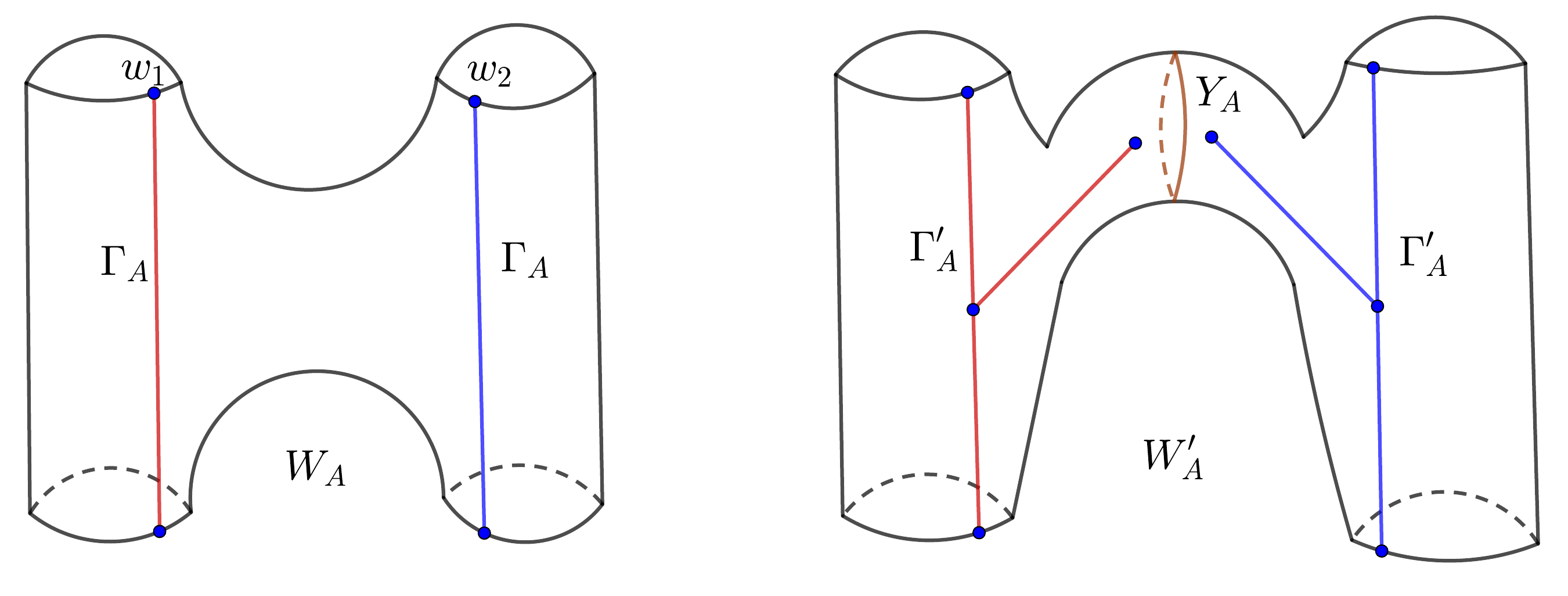}
\caption{Ribbon graph cobordisms $(W_A,\Ga_A)$ and $(W_A^\p,\Ga_A^\p)$.}
\label{excision51}
\end{figure}

We prove this for the case that $Y$ has two components $Y_1$ and $Y_2$. The proof for the case that $Y$ is connected is similar. For $i=1,2$, let $w_i\in Y_i$ be basepoints and let $\Ga_A\subset W_A$ consist of paths connecting basepoints $w_i$ in different ends of $W_A$; see the left subfigure of Figure \ref{excision51}. Suppose $W_A^\p$ is diffeomorphic to $W_A$ but drawn in a different position and suppose $\Ga_A^\p\subset W_A^\p$ is obtained from $\Ga_A$ by adding an arc to each path and choosing any ordering for the vertex with valence $3$; see the middle subfigure of Figure \ref{excision51}. By \cite[Section 11.2]{Zemke2019}, the ribbon graph cobordisms $(W_A,\Ga_A)$ and $(W_A^\p,\Ga_A^\p)$ induce the same cobordism map. Suppose $Y_A\cong \Sigma\times S^1\subset W_A$ is the manifold in the neck of $W_A^\p$. We know a neighborhood $N(Y_A)$ is diffeomorphic to $Y_0\times I$. Let $G_A$ consist of the images of $\Sigma$ in $\partial W_A$ and $\partial N(Y_0)$. 

\begin{figure}[ht]
\centering
\includegraphics[width=0.9\textwidth]{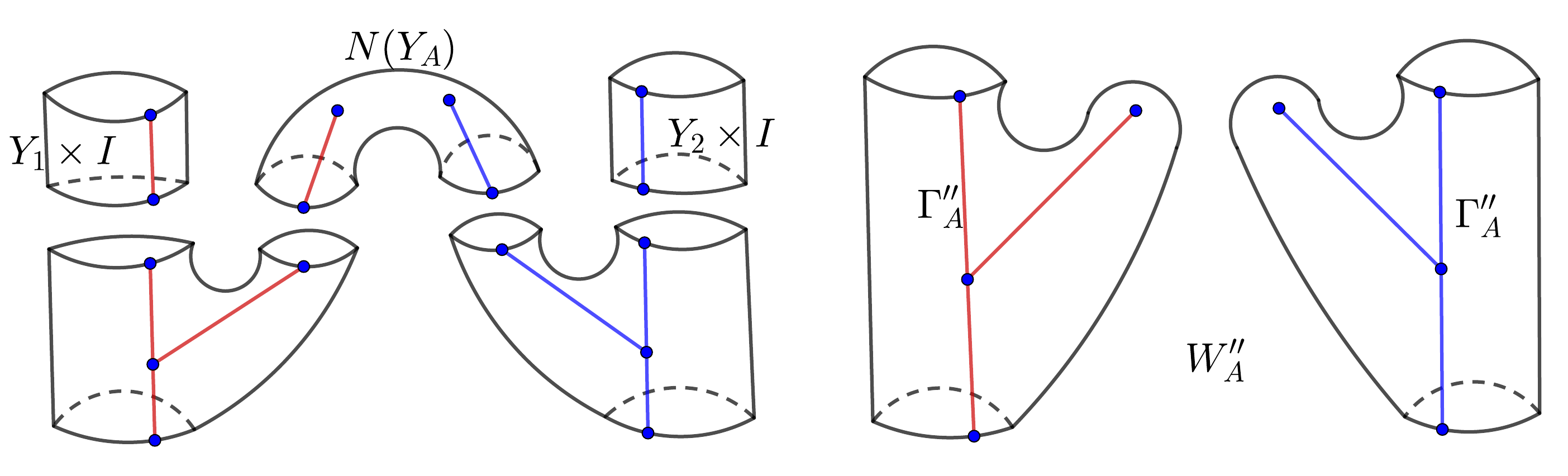}
\caption{Ribbon graph cobordisms $(W_A^\p,\Ga_A^\p)$ and $(W_A^\pp,\Ga_A^\pp)$.}
\label{excision52}
\end{figure}

By Proposition \ref{prop: composition}, we can decompose $(W_A^\p,\Ga_A^\p)$ into two parts as shown in the left subfigure of Figure \ref{excision52} and compute $HF(W_A,\Ga_A|G_A)$ by composition of two cobordism maps. The first part has three components corresponding to $Y_1\times I,N(Y_A)$, and $Y_2\times I$, respectively. By Lemma \ref{lem: two cobordism equal}, we can replace the component corresponding to $N(Y_A)$ by two components corresponding to $\Sigma\times D^2\sqcup (-\Sigma\times D^2)$ in the right subfigure of Figure \ref{excision11}. Then we know the cobordism map $HF(W_A,\Ga_A^\p|G_A)$ is the same as $HF(W_A^\pp,\Ga^\pp|G_A)$, where $(W_A^\pp,\Ga^\pp)$ is the ribbon graph cobordism in the right subfigure of Figure \ref{excision52}. By \cite[Section 11.2]{Zemke2019}, we can remove the arcs of $\Ga^\pp$ in the interior of the cobordism $W_A^\pp$. Then we know $HF(W_A^\pp,\Ga_A^\pp|G_A)$ is the identity map because $$(W_A^\pp,\Ga_A^\pp)\cong ((Y_1\sqcup Y_2)\times I,(w_1\sqcup w_2)\times I).$$Thus, the cobordism map $HF(W_A,\Ga|G_A)$ is the identity map.


{\bf Step 3}. For some restricted graph $\Ga_B$ and some surface $G_B$ in $W_B=W\cup_{Y}\bar{W}$, we show the cobordism map $$HF(W_B,\Ga_B|G_B)\deq HF^+(W_B,\Ga_B|G_B)=\hfm(W_B,\Ga_B|G_B)$$induces the identity map on $$HF(\widetilde{Y}|\widetilde{F})\deq HF^+(\widetilde{Y}|\widetilde{F})\cong \hfm(\widetilde{Y}|\widetilde{F}).$$

\begin{figure}[ht]
\centering
\includegraphics[width=0.95\textwidth]{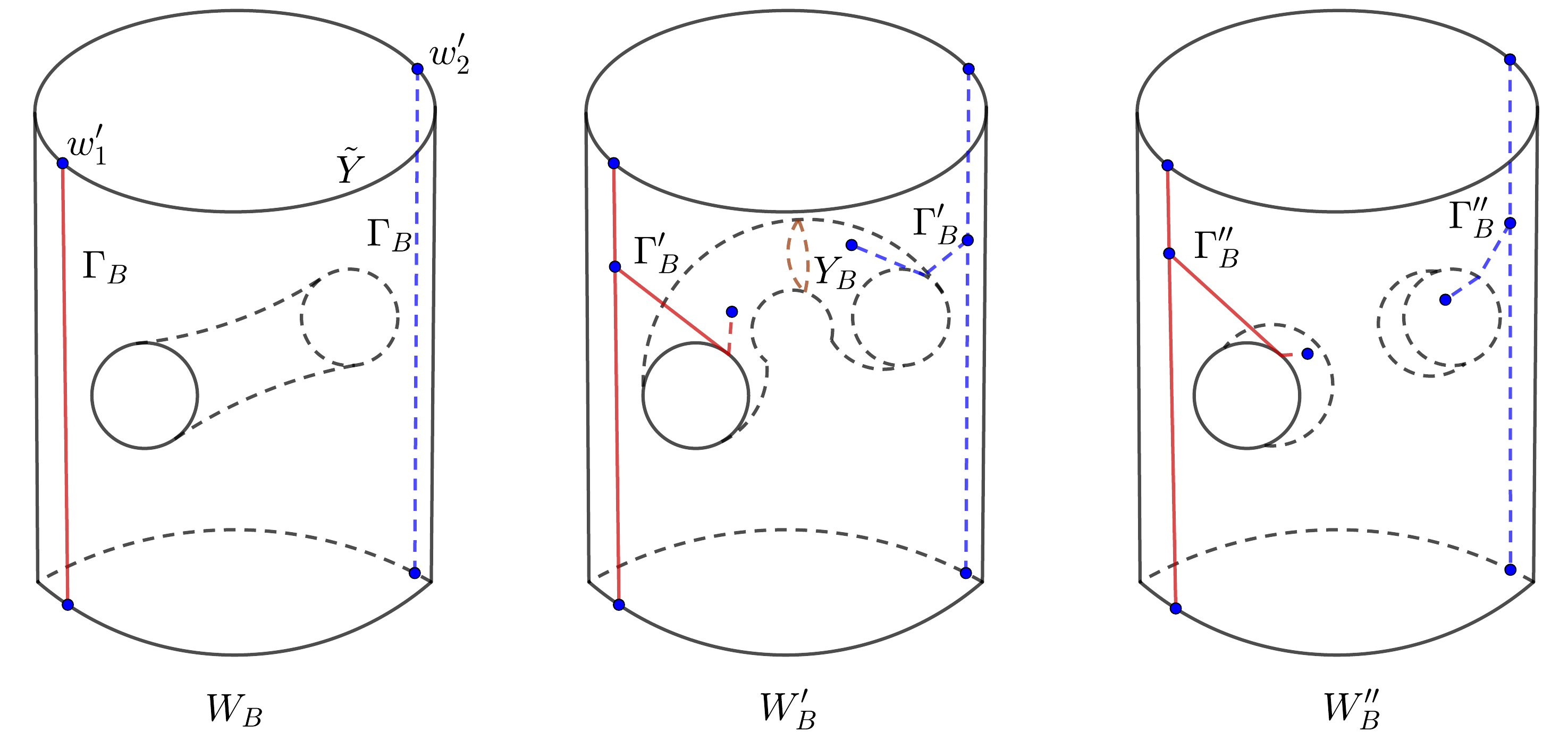}
\caption{Ribbon graph cobordisms $(W_B,\Ga_B),(W_B^\p,\Ga_B^\p)$, and $(W_B^\pp,\Ga_B^\pp)$.}
\label{excision61}
\end{figure}

We prove this for the case that $Y$ has two components $Y_1$ and $Y_2$. The proof for the case that $Y$ is connected is similar. The ribbon graph cobordism $(W_B,\Ga_B)$ is shown in the left subfigure of Figure \ref{excision61} and suppose endpoints of $\Ga_B$ correspond to $w_1^\p$ and $w_2^\p$ in $\widetilde{Y}$. The proof is essentially the same as that in Step 2. We first change the position of $W_B$ and add two arcs to $\Ga_B$ to obtain $(W_B^\p,\Ga_B^\p)$, as shown in the middle subfigure of Figure \ref{excision61}. Second, we choose $Y_B$ in the neck of $W_B^\p$ and set $G_B$ to be the images of $\Sigma$ in $\partial W_B^\p$ and $\partial N(Y_B)$. Third, we replace $N(Y_B)$ by $\Sigma\times D^2\sqcup (-\Sigma\times D^2)$ via Lemma \ref{lem: two cobordism equal} to obtain $(W_B^\pp,\Ga_B^\pp)$, as shown in the right subfigure of Figure \ref{excision61}. Finally we remove arcs in the interior of the cobordism and show it is the identity map because$$(W_B^\pp,\Ga_B^\pp)\cong (\widetilde{Y}\times I,(w_1^\p\sqcup w_2^\p)\times I).$$

Finally, we know Step 2 and Step 3 imply $$HF(Y|F)\cong HF(\widetilde{Y}|\widetilde{F})$$via cobordism maps associated to ribbon graph cobordisms $$(W,\Ga_A\cap W)\cong (W,\Ga_B\cap W)\aand(\bar{W},\Ga_A\cap \bar{W})\cong (\bar{W},\Ga_B\cap \bar{W}).$$Note that those ribbon graph cobordisms are restricted in the sense of Definition \ref{defn: restrict}.

\epf
\subsection{Sutured Heegaard Floer homology}\quad

In this subsection, we introduce two equivalent definitions of sutured Heegaard Floer homology. The first one is due to Juh\'{a}sz \cite{juhasz2006holomorphic}, based on balanced diagrams of balanced sutured manifolds. The other follows from the construction in Section \ref{sec: Floer theory for sutured manifold}, which is essentially due to Kronheimer and Mrowka \cite{kronheimer2010knots}. These definitions are denoted by $SFH$ and $\shf$, respectively. The equivalence of these definitions was shown by Lekili \cite{lekili2013heegaard} and Baldwin and Sivek \cite{Baldwin2020}. We will focus on the equality for graded Euler characteristics of two homologies.

\bdefn[{\cite[Section 2]{juhasz2006holomorphic}}]
A \textbf{balanced diagram} $\mch=(\Sigma,\al,\be)$ is a tuple  satisfying the following.
\begin{enumerate}[(1)]
    \item $\Sigma$ is a compact, oriented surface with boundary.
    \item $\al = \{ \al_1, \dots, \al_n \}$ and $\be = \{ \be_1,\dots, \be_n \}$ are two sets of pairwise disjoint simple closed curves in the interior of $\Sigma$.
    \item The maps $\pi_0(\partial \Sigma)\to \pi_0(\Sigma\backslash \al)$ and $\pi_0(\partial \Sigma)\to \pi_0(\Sigma\backslash\be)$ are surjective.
\end{enumerate}

For such triple, let $N$ be the 3-manifold obtained from $\Sigma\times [-1,1]$ by attaching 3–dimensional 2–handles along $\al_i \times \{-1\}$ and $\be_i \times \{1\}$ for $i=1,\dots,n$ and let $\nu=\partial \Sigma\times \{0\}$. A balanced diagram $(\Sigma,\al,\be)$ is called \textbf{compatible} with a balanced sutured manifold $(M,\ga)$ if the balanced sutured manifold $(N,\nu)$ is diffeomorphic to $(M,\ga)$. 

\edefn
Suppose $\mch=(\Sigma,\al,\be)$ is a balanced diagram with $g=g(\Sigma)$ and $n=|\al|=|\be|$. Suppose $\mch$ satisfies the admissible condition in \cite[Section 3]{juhasz2006holomorphic}. Consider two tori$$\mathbb{T}_\al\deq \al_1\times\cdots\times\al_n~{\rm and}~\mathbb{T}_\be\deq \be_1\times\cdots\times\be_n$$in the symmetric product $${\rm Sym}^n\Sigma\deq (\prod_{i=1}^n\Sigma)/S_n.$$

The chain complex $SFC(\mch)$ is a free $\ft$-module generated by intersection points $\bs{x}\in\mathbb{T}_\al\cap \mathbb{T}_\be$. Similar to the construction of $CF^-$, for a generic path of almost complex structures $J_s$ on $\sym^n\Sigma$, define the differential on $SFC(\mch)$ by$$\partial_{J_s}(\bs{x})=\sum_{\bs{y}\in \mathbb{T}_\al\cap \mathbb{T}_\be}\sum_{\substack{\phi\in\pi_2(\bs{x},\bs{y})\\\mu(\phi)=1}}\#\widehat{\mathcal{M}}_{J_s}(\phi)\cdot \bs{y}.$$



\bthm[{\cite{juhasz2006holomorphic,Juhasz2012}}]\label{thm: SFH}
Suppose $(M,\ga)$ is a balanced sutured manifold. Then there is an admissible balanced diagram $\mch$ compatible with $(M,\ga)$. The vector spaces $H(SFC(\mch),\partial_{J_s})$ for different choices of $\mch$ and $J_s$, together with some canonical maps, form a transitive system over $\mathbb{F}_2$. Let $SFH(M,\ga)$ denote this transitive system and also the associated actual group. Moreover, there is a decomposition$$SFH(M,\ga)=\bigoplus_{\mathfrak{s}\in \spin(M,\partial M)}SFH(M,\ga,\mathfrak{s}).$$
\ethm
\brem\label{rem: SFH}
The group $SFH(M,\ga)$ generalizes Heegaard Floer homology \cite{ozsvath2004holomorphic} and knot Flor homology \cite{ozsvath2004holomorphicknot,Rasmussen2003}. Suppose $Y$ is a closed 3-manifold and $K\subset Y$ is a knot. Let $Y(1)$ be obtained from $Y$ by removing a 3-ball and let $\delta$ be a simple closed curve on $\partial Y(1)$. Let $\ga$ consist of two meridians of $K$. Then there are isomorphisms \[SFH(Y(1),\delta)\cong \hf(Y) ~{\rm and}~ SFH(Y(K),\ga)\cong \hfk(Y,K).\]
\erem
\bdefn[]\label{defn: chi(sfh)}
For a balanced sutured manifold $(M,\ga)$, let the $\mathbb{Z}_2$-grading of $SFH(M,\ga)$ be induced by the sign of intersection points of $\mathbb{T}_\al$ and $\mathbb{T}_\be$ for some compatible diagram $\mch=(\Sigma,\al,\be)$ (\textit{c.f.} \cite[Section 3.4]{Friedl2009}). Suppose $H=H_1(M)/{\rm Tors}$ and choose any $\mathfrak{s}_0\in Spin^c(M,\ga)$. The \textbf{graded Euler characteristic} of $SFH(M,\ga)$ is $$\chi_{\rm gr}(SFH(M,\ga)) \colonequals \sum_{\substack{\mathfrak{s}\in \spin(M,\ga)\\\mathfrak{s}-\mathfrak{s}_0=h\in H^2(M,\partial M)}} \chi_{\rm gr}(SFH(M,\ga,\mathfrak{s})) \cdot p\circ {\rm PD}(h)\in \mathbb{Z}[H]/\pm H,$$where ${\rm PD}: H^2(M,\partial M)\to H_1(M)$ is the Poincar\'{e} duality map and $p: H_1(M)\to H_1(M)/{\rm Tors}$ is the projection map.
\edefn
\bthm[\cite{Friedl2009}]Suppose $(M,\ga)$ is a balanced sutured manifold. Then$$\chi_{\rm gr}(SFH(M,\ga))=p_*(\tau (M,\ga))\in\mathbb{Z}[H]/\pm H,$$where $\tau(M,\ga)$ is a (Turaev-type) torsion element computed from the map $$\pi_1(R_-(\ga))\to \pi_1(M)$$ by Fox calculus and $p_*$ is induced by $p: H_1(M)\to H_1(M)/{\rm Tors}=H$.
\ethm

Then we define the second version of sutured Heegaard Floer homology.
\bdefn
Suppose $(M,\ga)$ is a balanced sutured manifold and $(Y,R)$ is a closure of $(M,\ga)$ as in Definition \ref{defn: closure}. Define $$SHF(M,\ga)\deq HF(Y|R)= \bigoplus_{\substack{\mathfrak{s}\in \spin(Y|R)}}HF^+(Y,\mathfrak{s}).$$
\edefn
\brem
By work of Kutluhan, Lee, and Taubes \cite{kutluhan2010hf}, for any $\mathfrak{s}\in \spin(Y)$, there is an isomorphism$$HF^+(Y,\mathfrak{s})\cong \widecheck{HM}_*(Y,\mathfrak{s})=\widecheck{HM}_\bullet(Y,\mathfrak{s}).$$The last group is used to define $SHM$ in \cite{kronheimer2010knots}. 

Following the discussion in Section \ref{sec: Floer theory for sutured manifold}, we can prove the naturality of $SHF(M,\ga)$ based on Floer's excision theorem. Let $\shf(M,\ga)$ be the transitive system corresponding to $SHF(M,\ga)$.
\erem

\bthm[{\cite[Theorem 24]{lekili2013heegaard}, see also \cite[Theorem 3.26]{Baldwin2020}}]\label{thm: SFH equal SHF}
Suppose $(M,\ga)$ is a balanced sutured manifold and $(Y,R)$ is a closure of $(M,\ga)$. Then there exists a balanced diagram $\mch=(\Sigma,\al,\be)$ compatible with $(M,\ga)$ and a singly-pointed Heegaard diagram $\mch^\p=(\Sigma^\p,\al^\p,\be^\p,z)$ of $Y$ so that the following holds.
\benu
\item $\Sigma$ is a submanifold of $\Sigma^\p$.
\item $\al$ and $\be$ are subsets of $\al^\p$ and $\be^\p$, respectively.
\item Suppose $\al^\p=\al\cup \al^\pp$ and $\be^\p=\be\cup \be^\pp$. There exists an intersection point $\bs{x}_1\in\mathbb{T}_{\al^\pp}\cap \mathbb{T}_{\be^\pp}$ so that the map \begin{equation*}
    \begin{aligned}
    f:SFC(\mch)&\to CF^+(\mch^\p|R)\\
    \bs{c}\quad&\mapsto \bs{c}\times \bs{x}_1\\
\end{aligned}
\end{equation*}is a quasi-isomorphism, where $CF^+(\mch^\p|R)$ is the chain complex of $HF^+(Y|R)$ associated to $\mch^\p$.
\eenu
\ethm
\bcor[Proposition \ref{prop: SFH=SHF}]\label{cor: SFH=SHF}

Suppose $(M,\ga)$ is a balanced sutured manifold and $H=H_1(M)/{\rm Tors}\cong H^2(M,\partial M)/{\rm Tors}$. We have
\[SFH(M,\ga)\cong \shf(M,\ga)\]with respect to the grading associated to $H$ and the $\mathbb{Z}_2$ grading, up to a global grading shift.

In particular, we have$$\chi_{\rm gr}(SFH(M,\ga))=\chi_{\rm gr}(\shf(M,\ga))\in \mathbb{Z}[H]/\pm H,$$where $\chi_{\rm gr}(\shf(M,\ga))$ is defined as in Definition \ref{defn: chi(shg)}.
\ecor
\bpf
It suffices to show the quasi-isomorphism in Theorem \ref{thm: SFH equal SHF} respects spin$^c$ structures and $\mathbb{Z}_2$-gradings.

Consider the $\mathbb{Z}_2$-gradings at first. Suppose $\bs{c}_1$ and $\bs{c}_2$ are two generators of $SFC(\mch)$. Note that the $\mathbb{Z}_2$-grading of $\bs{c}_i$ is defined by the sign of the corresponding intersection point in $\mathbb{T}_\al\cap \mathbb{T}_\be$ for $i=1,2$. For $\bs{c}_i\times \bs{x}_1$, the $\mathbb{Z}_2$-grading is defined by mod 2 Maslov grading, which coincides with the sign of the corresponding intersection point in $\mathbb{T}_{\al^\p}\cap \mathbb{T}_{\be^\p}$. Thus, we have$${\rm gr}_2(\bs{c}_1)-{\rm gr}_2(\bs{c}_2)={\rm gr}_2(\bs{c}_1\times \bs{x}_1)-{\rm gr}_2(\bs{c}_2\times \bs{x}_1),$$where ${\rm gr}_2$ is the $\mathbb{Z}_2$-grading.  

Then we consider spin$^c$ structures. Consider $\bs{c}_i$ for $i=1,2$ as above. From \cite[Lemma 4.7]{juhasz2006holomorphic}, there is a one chain $\ga_{\bs{c}_1}-\ga_{\bs{c}_2}$ such that $$\mathfrak{s}(\bs{c}_1)-\mathfrak{s}(\bs{c}_2)={\rm PD}([\ga_{\bs{c}_1}-\ga_{\bs{c}_2}]),$$where $\mathfrak{s}(\cdot):\mathbb{T}_\al\cap \mathbb{T}_\be\to \spin(M,\partial M)$ is defined in \cite[Definition 4.5]{juhasz2006holomorphic}, and ${\rm PD}:H_1(M)\to H^2(M,\partial M)$ is the Poincar\'{e} duality map. 

From \cite[Lemma 2.19]{ozsvath2004holomorphic}, we have $$\mathfrak{s}_z(\bs{c}_1\times \bs{x}_1)-\mathfrak{s}_z(\bs{c}_2\times \bs{x}_1)={\rm PD}^\p(i_*([\ga_{\bs{c}_1}-\ga_{\bs{c}_2}])),$$where $\mathfrak{s}_z(\cdot):\mathbb{T}_{\al^\p}\cap \mathbb{T}_{\be^\p}\to \spin(Y)$ is defined in \cite[Section 2.6]{ozsvath2004holomorphic} and ${\rm PD}^\p:H_1(Y)\to H^2(Y)$ is the Poincar\'{e} duality map, and $i_*:H_1(M)\to H_1(Y)$ is the map induced by inclusion $i:M\to Y$. 

Hence we have $$c_1(\mathfrak{s}_z(\bs{c}_1\times \bs{x}_1))-c_1(\mathfrak{s}_z(\bs{c}_2\times \bs{x}_1))=2{\rm PD}^\p(i_*([\ga_{\bs{c}_1}-\ga_{\bs{c}_2}])).$$

Finally, the argument about graded Euler characteristics follows from definitions.
\epf

\section{The graded Euler characteristic of formal sutured homology}
\label{sec: Equivalence of graded Euler characteristics}

In this section, we prove the graded Euler characteristic of formal sutured homology is independent of the choice of the Floer-type theory. Throughout this section, we assume that $\hg$ is a Floer-type theory, \textit{i.e.}, it satisfies all three Axioms (A1), (A2), and (A3). For simplicity, we say `a property is independent of $\hg$' if a property is independent of the choice of the Floer-type theory. Suppose $(M,\ga)$ is a balanced sutured manifold. If the admissible surfaces and the closure of $(M,\ga)$ are fixed, then the graded Euler characteristic $\chi_{\rm gr}(\sh(M,\ga))$ in Definition \ref{defn: chi(shg)} is considered as a well-defined element in $\mathbb{Z}[H_1(M)/{\rm Tors}]$, rather than $\mathbb{Z}[H_1(M)/{\rm Tors}]/\pm (H_1(M)/{\rm Tors})$; see Remark \ref{rem: fix closure}. Note that in this subsection, we avoid using $H$ to denote $H_1(M)/{\rm Tors}$ and the symbol $H$ usually denotes a handlebody.

\subsection{Balanced sutured handlebodies}\quad

In this subsection, we deal with $\mathbb{Z}^n$-gradings for a balanced sutured handlebody. We start with the following lemma about the sign ambiguity.

\blem\label{lem: from one closure to another}
Suppose $(M,\ga)$ is a balanced sutured manifold, $S\subset (M,\ga)$ is an admissible surface. Suppose $(Y_1,R_1)$ and $(Y_2,R_2)$ are two closures of $(M,\ga)$ of the same genus so that $S$ extends to closed surfaces $\bar{S}_1$ and $\bar{S}_2$ as in Subsection \ref{subsec: gradings on SH}. If $\gr(\hg(Y_1|R_1))$ is already determined without the sign ambiguity, then $\gr(\hg(Y_2|R_2))$ is determined without the sign ambiguity from $\gr(\hg(Y_1|R_1))$ and the topological data of $(Y_1,R_1)$ and $(Y_2,R_2)$.
\elem
\bpf
In Subsection \ref{sec: Floer theory for sutured manifold}, we construct a canonical map
$$\Phi_{12}:\hg(Y_1|R_1)\ra \hg(Y_2|R_2).$$
From the proof of Theorem \ref{thm: grading is well defined, appendix}, the canonical map $\Phi_{12}$ necessarily preserves the grading induced by $S$. From the construction of $\Phi_{12}$ in Subsection \ref{sec: Floer theory for sutured manifold}, the canonical map is a composition of a few cobordism maps (or the inverse). Then the $\mathbb{Z}_2$-grading shift follows from Axiom (A3-3) .
\epf

Next, we consider gradings associated to admissible surfaces. To fix the ambiguity of $H_1(M)/{\rm Tors}$, we will fix the choices of admissible surfaces. For sutured handlebodies, we start with embedded disks.
\bprop\label{prop: Euler char for sutured handlebody with meridian disk}
Suppose $H$ is a genus $g>0$ handlebody and $\gamma\subset \partial H$ is a closed oriented 1-submanifold so that $(H,\gamma)$ is a balanced sutured manifold. Pick $D_1,\dots,D_g$ a set of pairwise disjoint meridian disks in $H$ so that
$[D_1],\dots,[D_g]$ generate $H_2(H,\partial H)$. Then for any fixed multi-grading $\bs{i}=(i_1,\dots,i_g)\in\intg^g$ associated to $D_1,\dots,D_g$, the Euler characteristic
$$\chi(\shg(-H,-\ga,\bs{i}))\in\intg\slash\{\pm1\}$$
depends only on $(H,\ga)$, $D_1,\dots,D_g$ and $\bs{i}\in\intg^g$, and is independent of $\hg$. Furthermore, if a particular closure of $(-H,-\ga)$ is fixed, then the sign ambiguity can be removed.
\eprop

\bpf
We fix the handlebody $H$ and the set of disks $D_1,\dots,D_g\subset H$. For any suture $\ga$ on $\partial H$, define
$$I(\ga)=\min_{\ga'\text{ is isotopic to }\ga}\sum_{j=1}^g|D_j\cap\ga'|,$$
where $|\cdot|$ denotes the number of points. We prove the proposition by induction on $I(\ga)$. Since $[\ga]=0\in H_1(\partial H)$, we know $|D_j\cap\ga|$ is always even for $j=1,\dots,g$.

Note that an isotopy of $\gamma$ can be understood as combinations of positive and negative stabilizations in the sense of Definition \ref{defn_2: stabilization of surfaces}, and the grading shifting behavior under such isotopies (stabilizations) is described by Theorem \ref{thm_2: grading shifting property}, which is determined purely by topological data and is independent of $\hg$. Hence we can assume that the suture $\ga$ has already realized $I(\ga)$.

First, if $I(\ga)<2g$, then there exists a meridian disk $D_j$ with $D_j\cap\ga=\emptyset.$ Then it follows from Proposition \ref{prop: SHI detects tautness} that $\shg(-H,-\ga)=0$ since $-H$ is irreducible while $(-H,-\ga)$ is not taut. Hence for any multi-grading $\bs{i}\in\intg^g$, we have $\chi(\shg(-H,-\ga,\bs{i}))=0.$

If $I(\ga)=2g$, then either there exists some integer $j$ so that $D_j\cap\ga=\emptyset$ or for $j=1,\dots,g$, we have $|D_j\cap\ga|=2.$ In the former case, we know that $\shg(-H,-\ga)=0$ and hence $\chi(\shg(-H,-\ga,\bs{i}))=0$
for any multi-grading $\bs{i}\in\intg^g$. In the later case, we know that $(-H,-\ga)$ is a product sutured manifold. It follows from Proposition \ref{prop: product} and Proposition \ref{thm_2: grading in SHG} that
$$\shg(-H,-\ga)=\shg(-H,-\ga,\bs{0})\cong\mathbb{F}.$$
Hence
$$\chi(\shg(-H,-\ga,\bs{i}))=\begin{cases}\pm1&\bs{i}=\bs{0}=(0,\dots,0)\\0& \bs{i}\in\intg^g\backslash\{\bs{0}\}\end{cases}$$
Note that the ambiguity $\pm1$ comes from the choice of the closure. If we choose a particular closure $Y$ of $(-H,-\ga)$, then the Euler characteristic has no sign ambiguity. Since $(H,\ga)$ is a product sutured manifold, there is a `standard' closure $(S^1\times\Sigma,\{1\}\times\Sigma)$ as in \cite{kronheimer2010knots}. By Axiom (A3-2), we have $$\chi(\hg(S^1\times\Sigma|\{1\}\times\Sigma))=-1.$$
Then for any other closure $(Y,R)$, by Lemma \ref{lem: from one closure to another} $\gr(\shg(Y|R))$ has no sign ambiguity.

Now suppose we have proved that, for all $\ga$ so that
$I(\ga)<2n$, the Euler characteristic of $\shg(-H,-\ga,\bs{i})$, viewed as an element in $\intg\slash\{\pm1\}$, is independent of $\hg$, and that when we choose any fixed closure of $(-H,-\ga)$, the sign ambiguity can be removed. Next we deal with the case when $I(\ga)=2n$. 

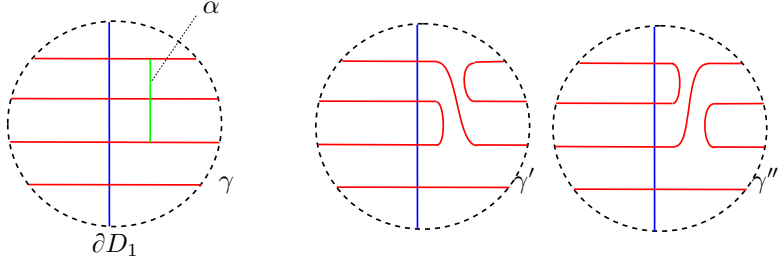
\begin{figure}[ht]
\centering
\begin{overpic}[width=0.7\textwidth]{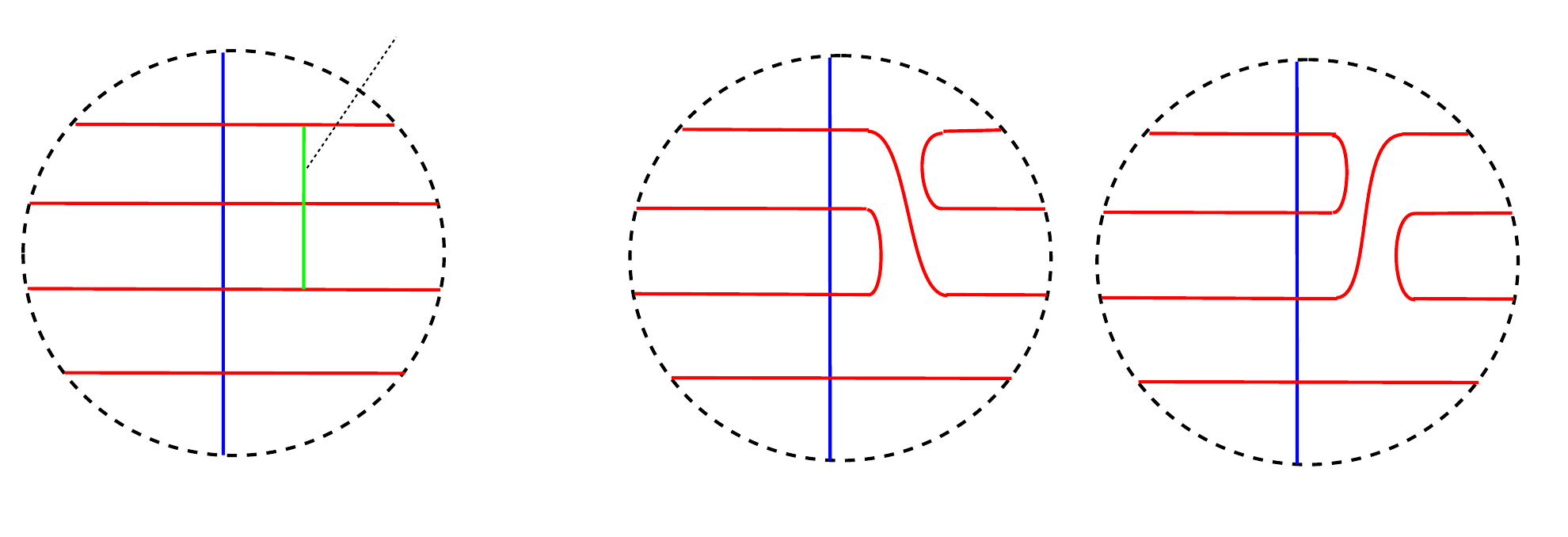}
    \put(26,32){$\alpha$}
    \put(12,2){$\partial D_1$}
    \put(28,10){$\ga$}
    \put(65,10){$\ga'$}
    \put(95,10){$\ga''$}
\end{overpic}
\vspace{-0.05in}
\caption{The bypass arc $\alpha$ that reduces the intersection function $I$.}\label{fig: reduce intersections}
\end{figure}

Note that we have dealt with the base case $I(\ga)\leq2g$, so we can assume that $n\geq g+1$. Hence, without loss of generality, we can assume that $|D_1\cap\ga|\geq 4$. Within a neighborhood of $\partial D_1$, the suture $\ga$ can be depicted as in Figure \ref{fig: reduce intersections}. We can pick the bypass arc $\alpha$ as shown in the same figure. From Proposition \ref{prop: bypass maps preserves gradings}, for any multi-grading $\bs{i}\in\intg^g$, we have an exact triangle
\begin{equation}\label{eq: by exact triangle for sutured handlebody}
\xymatrix{
&\shg(-H,-\ga,\bs{i})\ar[dr]&\\
\shg(-H,-\ga'',\bs{i})\ar[ur]&&\shg(-H,-\ga',\bs{i})\ar[ll]
}	
\end{equation}

Note that the suture $\ga'$ and $\ga''$ are determined by the original suture $\ga$ and the bypass arc $\alpha$, which are all topological data. From Figure \ref{fig: reduce intersections}, it is clear that
$$I(\ga')\leq I(\ga)-2~{\rm and}~I(\ga'')\leq I(\ga)-2.$$
Hence the inductive hypothesis applies, and we know that the Euler characteristics of $\shg(-H,-\ga'',\bs{i})$ and $\shg(-H,-\ga',\bs{i})$ can be fixed independently of $\hg$. Note that the maps in the bypass exact triangle (\ref{eq: by exact triangle for sutured handlebody}) are described by Proposition \ref{prop_2: description of bypass maps}. Hence we conclude that the Euler characteristic of $\shg(-H,-\ga,\bs{i})$ is also independent of $\hg$. Thus, we finish the proof by induction.
\epf

Next, we deal with gradings associated to general admissible surfaces.

\bprop\label{prop: Euler char for sutured handlebody with arbitrary surface}Suppose $H$ is a genus $g$ handlebody, and $S$ is a properly embedded surface in $H$. Suppose $\ga\subset \partial H$ is a suture so that $(H,\ga)$ is a balanced sutured manifold and $S$ is an admissible surface. Then the Euler characteristic
$$\chi(\shg(-H,-\ga,S,j))\in\intg\slash\{\pm1\}$$
depends only on $(H,\ga)$, $S$, and $j\in\intg$ and is independent of $\hg$. Furthermore, if we fix a particular closure of $(-H,-\ga)$, then the sign ambiguity can also be removed.
\eprop
Before proving the proposition, we need the following lemma.
\blem\label{lem: null homotopic boundary}
Suppose $(M,\ga)$ is a balanced sutured manifold and $S\subset (M,\ga)$ is a properly embedded admissible surface. Suppose $\alpha$ is a boundary component of $S$ so that $\alpha$ bounds a disk $D\subset\partial M$ and $|\al\cap\ga|=2$. Let $S'$ be the surface obtained by taking the union $S\cup D$ and then push $D$ into the interior of $M$. Then for any $i\in\intg$, we have
$$\shg(M,\ga,S,i)=\shg(M,\ga,S',i).$$
\elem
\bpf
Push the interior of $D$ into the interior of $M$ and make $D\cap S'=\emptyset$. It is clear that
$$[S]=[S'\cup D]\in H_2(M,\partial M){\rm~and~}\partial S=\partial(S'\cup D).$$
In Subsection \ref{subsec: gradings on SH}, when constructing the grading associated to $S'\cup D$, we can pick a closure $(Y,R)$ of $(M,\ga)$, so that $S'$ and $D$ extend to closed surfaces $\bar{S}'$ and $\bar{D}$ in $Y$, respectively. Since $|\partial D\cap\ga|=2$, we know that $\bar{D}$ is a torus. Since $\partial S=\partial(S'\cup D)$, we know that $S$ also extends to a closed surface $\bar{S}$ and from the fact that $[S]=[S'\cup D]$ we know that
$$[\bar{S}]=[\bar{S}'\cup\bar{D}]=[\bar{S}']+[\bar{D}].$$
Since $\bar{D}$ is a torus, from Axioms (A1-4) and (A1-6), we know that the decompositions of $\hg(Y|R)$ with respect to $\bar{S}$ and $\bar{S}'$ are the same. Thus it follows that
$$\shg(M,\ga,S,i)=\shg(M,\ga,S',i).$$
\epf

\bpf[Proof of Proposition \ref{prop: Euler char for sutured handlebody with arbitrary surface}.]
It is a basic fact that the map
$$\partial_*:H_2(H,\partial H)\ra H_1(\partial H)$$
is injective, and $H_2(H,\partial H)$ is generated by $g$ meridian disks, which we fix as $D_1,\dots,D_g$. Hence we assume that
$$[S]=a_1[D_1]+\dots+a_g[D_g]\in H_2(H,\partial H).$$

\textbf{Case 1}. $\partial S$ consists of only $\partial D_i$, \textit{i.e.},
$$\partial S=\bigcup_{i=1}^g(\cup_{a_i}\partial D_i),$$where $\cup_{a_i}\partial D_i$ means the union of $a_i$ parallel copies of $\partial D_i$. 

Then it follows immediately from the construction of the grading and Axiom (A1-6) that
\beq
\shg(-H,-\ga,S,j)&=\shg(-H,-\ga,\bigcup_{i=1}^g(\cup_{a_i}D_i),j)\\
&=\bigoplus_{j_1+\dots+j_g=j}\shg(-H,-\ga,(D_1,\dots,D_g),(j_1,\dots,j_g)).
\eeq
Hence this case follows from Proposition \ref{prop: Euler char for sutured handlebody with meridian disk}.

\textbf{Case 2}. $\partial S$ contains some component that is not parallel to $\partial D_i$ for $j=1,\dots,g$. 

\textbf{Step 1}. We modify $S$ and show that it suffices to deal with the case when $S\cap D_j=\emptyset$ for $j=1,\dots,g$.

Note that ${\rm im}(\partial_*)\subset H_1(\partial H)$ is generated by $[\partial D_1],\dots,[\partial D_g]$, so we have $\partial S\cdot \partial D_i=0$ for $j=1,\dots,g$. Here $\cdot$ denotes the algebraic intersection number of two oriented curves on $\partial H$. This means that for $j=1,\dots,g$, the intersection points of $\partial D_i$ with $\partial S$ can be divided into pairs. Suppose two intersection points of $\partial D_1$ with $\partial S$ of opposite signs are adjacent to each other on $\partial D_1$, as depicted in Figure \ref{fig: double curve surgery, 1}. We can perform a cut and paste surgery along $D_1$ and $S$ to obtain a new surface $S_1$. From the same figure, it is clear that after isotopy, we can make 
\begin{equation*}
|\partial D_1\cap\partial S_1|\leq |\partial D_1\cap \partial S|-2.
\end{equation*}

\begin{figure}[ht]
\centering
\begin{overpic}[width=0.7\textwidth]{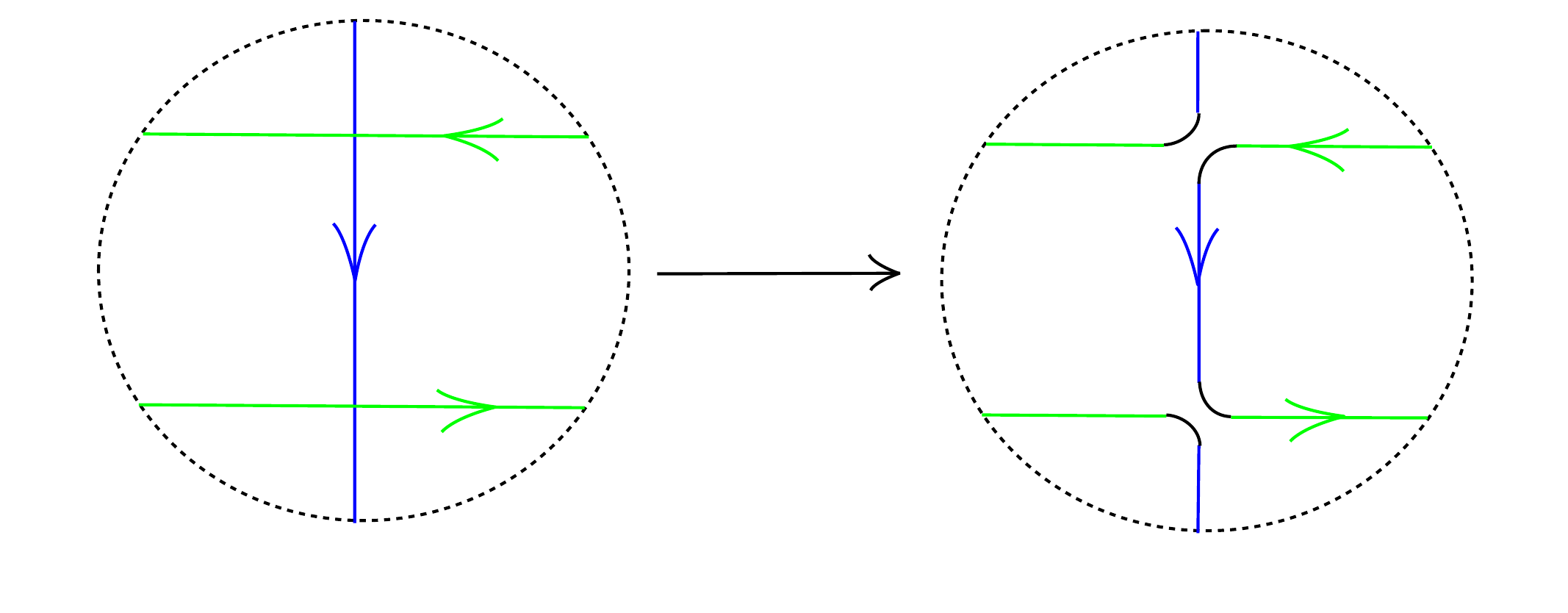}
    \put(21,2){$D_1$}
    \put(28,10){$S$}
    \put(42,23){cut and}
    \put(44,18){paste}
    \put(85,10){$S_1$}
\end{overpic}
\vspace{-0.05in}
\caption{The cut and paste surgery on $D_1$ and $S$.}\label{fig: double curve surgery, 1}
\end{figure}

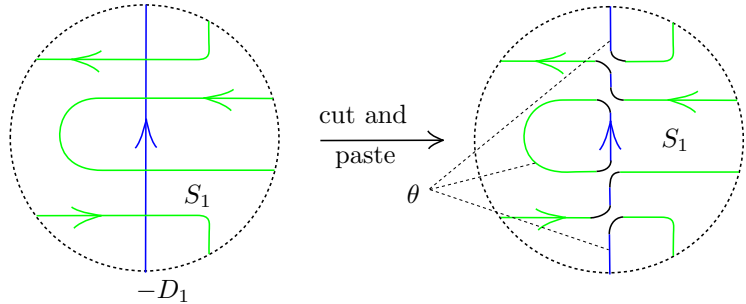
\begin{figure}[ht]
\centering
\begin{overpic}[width=0.7\textwidth]{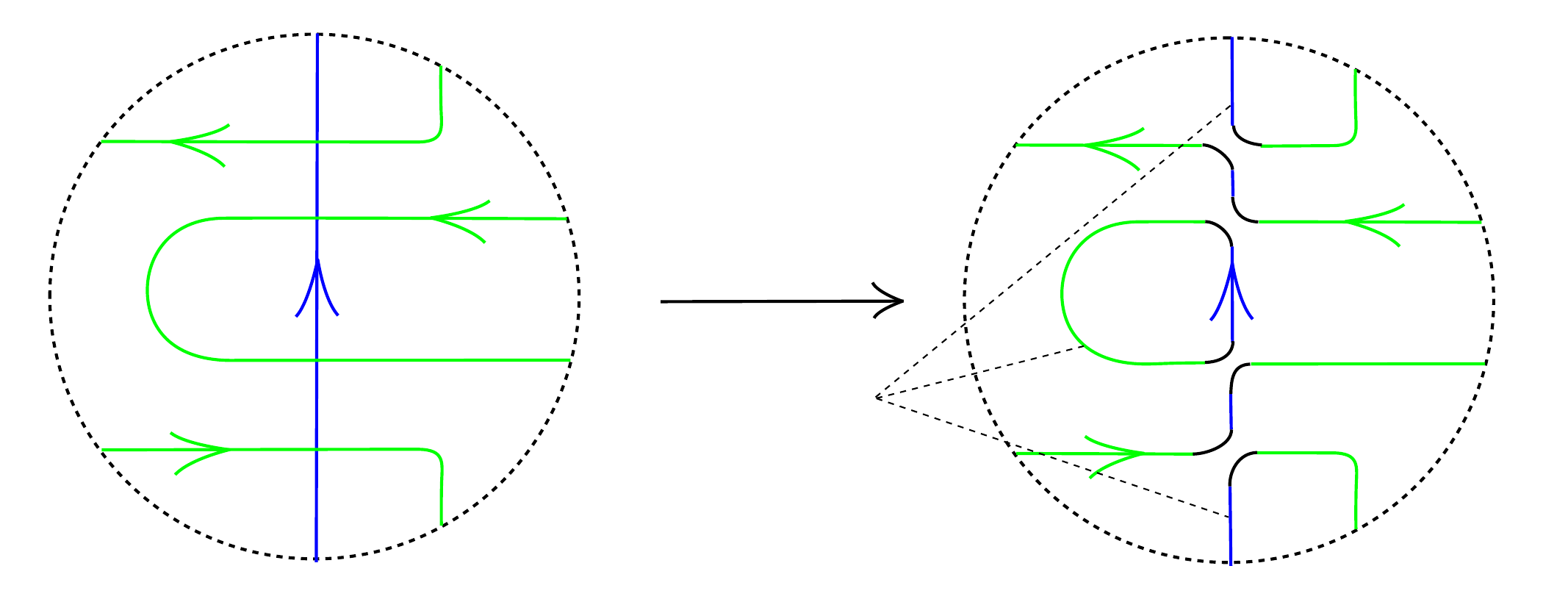}
    \put(19,0){$-D_1$}
    \put(25,12){$S_1$}
    \put(42,22){cut and}
    \put(44,17){paste}
    \put(85,19){$S_1$}
    \put(53,12){$\theta$}
\end{overpic}
\vspace{-0.05in}
\caption{The cut and paste surgery on $-D_1$ and $S_1$.}\label{fig: double curve surgery, 2}
\end{figure}

Note that if we perform a cut and paste surgery along $S_1$ and $-D_1$, we obtain another surface $S_2$. From Figure \ref{fig: double curve surgery, 2} it is clear that $\partial S_2=\partial S\cup\theta,$ where $\theta$ is the union of some null-homotopic closed curves on $\partial H$. We can isotope $S_2$ to make each component of $\theta$ intersects the suture twice. Let $S_3$ be the resulting surface of such an isotopy and $S_4$ be the surface obtained from $S_3$ by capping off every component of $\theta$. Then we have
$$[S]=[S_4]\in H_2(H,\partial H) {\rm~and~} \partial S=\partial S_4.$$
Hence from Lemma \ref{lem: null homotopic boundary} we know that
\beq
\shg(-H,-\ga,S,j)&=\shg(-H,\ga,S_4,j)\\
&=\shg(-H,\ga,S_3,j)\\
&=\shg(-H,-\ga,S_2,j+j(S_2,S_3))\\
&=\bigoplus_{j_1+j_2=j+j(S_2,S_3)}\shg(-H,-\ga,(D_1,S_1),(j_1,j_2))\\
\eeq
By Theorem \ref{thm_2: grading shifting property}, the shift $j(S_2,S_3)$ depends on the isotopy from $S_2$ to $S_3$, which is determined by the topological data and is independent of $\hg$. Hence we reduce the problem to understanding the Euler characteristic of $\shg(-H,-\ga)$ with multi-grading associated to $(D_1,S_1)$, with
$$|\partial D_1\cap\partial S_1|\leq |\partial D_1\cap \partial S|-2.$$
Repeating this argument, we finally reduce to the problem of understanding the Euler characteristic of $\shg(-H,-\ga)$ with multi-grading associated to $(D_1,\dots,D_g,S_g)$, with 
$$\partial D_i\cap \partial S_g=\emptyset\text{ for }j=1,\dots,g.$$

\textbf{Step 2}. We modify $S$ further to reduce to Case 1.

If every component of $\partial S_{g}$ is homotopically trivial, then we know that
$$[S_{g}]=0\in H_2(H,\partial H),$$since the map $H_2(H,\partial H)\to H_1(\partial H)$ is injective. We isotope each component of $\partial S_{g}$ by stabilization to make it intersect the suture $\ga$ twice and then cap it off by a disk. The resulting surface $S_{g+1}$ is a homologically trivial closed surface in $H$, so $\shg(-H,-\ga)$ is totally supported at grading $0$ with respect to $S_{g+1}$. The grading shift between $S_{g}$ and $S_{g+1}$ can then be understood by Theorem \ref{thm_2: grading shifting property}, and is independent of $\hg$.

Note that $\partial H\backslash(\partial D_1\cup\cdots\cup\partial D_g)$ is a $2g$-punctured sphere, so $\partial S$ is homotopically trivial when removing punctures on the sphere. If some component $C$ of $\partial S_{g}$ is not null-homotopic, then $C$ is obtained from some $\partial D_j$ by performing handle slides (or equivalently band sums) over $\partial D_1,\dots,\partial D_g$ for some times. 

If we isotope $C$ to make it intersect some $\partial D_i$ twice and then apply the cut and paste surgery, the resulting curve is isotopic to the one obtained by performing a handle slide over $\partial D_i$. Explicitly, in Figure \ref{fig: double curve surgery, 1}, suppose two right endpoints of arcs in $\partial S$ (the green arcs) are connected, then the right part of $\partial S_1$ is a trivial circle, and the left part of $\partial S_1$ is obtained from $\partial S$ by performing a handle slide over $\partial D_1$. Thus, we can apply the cut and paste surgery for many times, which is equivalent to performing handle slides over $\partial D_1,\dots,\partial D_g$ for some times. Finally, we reduce $C$ to the curve isotopic to $\partial D_j$. Then we reduce the problem to understanding the Euler characteristic of $\shg(-H,-\ga)$ with multi-grading associated to $(D_1,\dots,D_g,S_{g+2})$, where $S_{g+2}$ is a surface so that each component of $\partial S_{g+2}$ is parallel to $\pm\partial D_{i}$ for some $i$. Case 1 applies to $S_{g+2}$, and we finish the proof.
\epf

\bcor\label{cor: Euler char for sutured handlebody with arbitrary surface}
Suppose $H$ is a handlebody and $\ga$ is a suture on $\partial H$ so that $(H,\ga)$ is a balanced sutured manifold. Suppose $S_1,\dots,S_n$ are properly embedded admissible surfaces in $(H,\ga)$. Then the Euler characteristic
$$\chi(\shg(-H,-\ga,(S_1,\dots,S_n),(i_1,\dots,i_n)))\in\intg\slash\{\pm1\}$$
depends only on $(H,\ga)$, $S_1,\dots,S_n$, and $(i_1,\dots,i_n)\in\intg^n$, and is independent of $\hg$. Furthermore, if we fix a particular closure of $(-H,-\ga)$, then the sign ambiguity can also be removed.
\ecor
\bpf
The proof is similar to that for Proposition \ref{prop: Euler char for sutured handlebody with arbitrary surface}.
\epf

\subsection{Gradings about contact 2-handle attachments}\label{subsec: 2-handle attachment}\quad

In this subsection, we prove a technical proposition about the grading behavior for the map associated to contact 2-handle attachments.

Suppose $M$ is a compact oriented 3-manifold with boundary, and $S\subset M$ is a properly embedded surface. Suppose $\alpha\subset M$ is a properly embedded arc that intersects $S$ transversely and $\partial \al\cap\partial S=\emptyset$.
Let $N=M\backslash {\rm int}(N(\al))$, $S_N=S\cap N$, and $\mu\subset \partial N$ be a meridian of $\alpha$ that is disjoint from $S_N$. Let $\ga_N$ be a suture on $\partial N$ satisfies the following properties.
\benu
\item $(N,\ga_N)$ is balanced, $S$ is admissible, and $|\ga_N\cap\mu|=2$.
\item If we attach a contact 2-handle along $\mu$ in the sense of \cite[Section 4.2]{baldwin2016contact}, then we obtain a balanced sutured manifold $(M,\ga_M)$.
\eenu
From \cite[Section 4.2]{baldwin2016contact}, there is a map
$$C_{\mu}:\shg(-N,-\ga_N)\ra\shg(-M,-\ga_M)$$constructed as follows.

Push $\mu$ into the interior of $N$ to become $\mu'$. Suppose $(N_0,\ga_{N,0})$ is the manifold obtained from $(N,\ga_N)$ by a $0$-surgery along $\mu'$ with respect to the framing from $\partial N$. Equivalently, $(N_0,\ga_{N,0})$ can be obtained from $(M,\ga_M)$ by attaching a $1$-handle. Since $\mu'\subset {\rm int}(N)$, the construction of the closure of $(N,\ga_N)$ does not affect $\mu'$. Thus, we can construct a cobordism between closures of $(N,\ga_N)$ and $(N_0,\ga_{N,0})$ by attaching a $4$-dimensional $2$-handle associated to the surgery on $\mu'$. This cobordism induces a cobordism map
$$C_{\mu'}:\shg(-N,-\ga_N)\ra\shg(-N_0,-\ga_{N,0}).$$
It is shown in \cite[Section 4.2]{baldwin2016contact} (or also \cite[Section 6]{kronheimer2010knots}) that attaching a product 1-handle does not change the closure, so there is an identification
$$\iota:\shg(-M,-\ga_M)\xra{=}\shg(-N_0,-\ga_{N,0}).$$
Thus, we define $$C_{\mu}=\iota^{-1}\circ C_{\mu'}.$$

The main result of this subsection is the following proposition.
\bprop\label{prop: 2-handle map preserves the grading}Consider the setting as above. For any $i\in\intg$, we have
$$C_{\mu}(\shg(-N,-\ga_N,S_N,i))\subset \shg(-M,-\ga_M,S,i).$$
\eprop
\bpf
\textbf{Step 1}. We consider the grading behavior of the map $C_{\mu'}$ for gradings associated to $S_N$ and $S$. 

Since $\mu$ is disjoint from $S$, so we can also make $\mu'$ disjoint from $S_N=S\cap N$. As a result, the surface $S_N$ survives in $(N_0,\ga_{N,0})$. From Axiom (A1-7), the cobordism map associated to the $0$-surgery along $\mu'$ preserves the grading associated to $S_N$
$$C_{\mu'}(\shg(-N,-\ga_N,S_N,i))\subset \shg(-N_0,-\ga_{N,0},S_N,i).$$

\textbf{Step 2}. We show $\iota:\shg(-M,-\ga_M,S,i)\xra{=}\shg(-N_0,-\ga_{N,0},S,i).$

As discussed above, $(N_0,\ga_{N,0})$ is obtained from $(M,\ga_M)$ by a product 1-handle attachment. This product 1-handle can be described explicitly as follows. In $(N_0,\ga_{N,0})$, there is an annulus $A$ bounded by $\mu$ and its push-off $\mu'$. We can cap off $\mu'$ by the disk coming from the $0$-surgery, and hence obtain a disk $D$ with $\partial D=\mu$. By assumption, we know that $|\partial D\cap\ga_{N,0}|=|\mu\cap \ga_N|=2$. Hence $D$ is a compressing disk that intersects the suture twice. If we perform a sutured manifold decomposition on $(N_0,\ga_{N,0})$ along $D$, it is straightforward to check the resulting balanced sutured manifold is $(M,\ga_M)$. However, in \cite{juhasz2016cobordisms}, it is shown that decomposing along such a disk is the inverse operation of attaching a product 1-handle, and the disk is precisely the co-core of the product 1-handle. From this description, we can consider the product 1-handle attached to $(M,\ga_M)$ as along two endpoints of $\al$. Since $\partial \al\cap\partial S=\emptyset$, the surface $S$ naturally becomes a properly embedded surface in $(N_0,\ga_{N,0})$. From Axiom (A1-7), we know that the map $\iota$ preserves the gradings as claimed.

\textbf{Step 3}. We show $\shg(-N_0,-\ga_{N,0},S,i)=\shg(-N_0,-\ga_{N,0},S_N,i).$

If $S\cap \al=\emptyset$, then $S=S_N=S\cap N$ and the above argument is trivial. If $S\cap \al\neq\emptyset$, then $S_N$ is obtained from $S$ by removing disks containing intersection points in $\al\cap S$. Then $\partial S_N\backslash \partial S$ consists of a few copies of meridians of $\alpha$. For simplicity, we assume that there is only one copy of the meridian of $\alpha$, \textit{i.e.}, $\partial S_N\backslash \partial S=\mu$. The general case is similar to prove. 

After performing the $0$-surgery along $\mu'$, we know that the surface $S_N\subset N_0$ is compressible. Indeed, we can pick $\mu''\subset{\rm int}(S_N)$ parallel to $\mu\subset\partial S_N$. Then there is an annulus $A'$ bounded by $\mu''$ and $\mu'$, and we obtain a disk $D'$ by capping $\mu'$ off by the disk coming from the $0$-surgery. Performing a compression along the disk $D'$, we know that $S_N$ becomes the disjoint union of a disk $D''$ and the surface $S\subset N_0$. Note $\partial D''$ is parallel to the disk $D$ discussed above. Since 
$$\partial(D''\cup S)=\partial S_N~{\rm and}~[D''\cup S]=[S_N]\in H_2(N_0,\partial N_0),$$
From (A1-6), we know that
\beq
\shg(-N_0,-\ga_{N,0},S_N,i)&=\shg(-N_0,-\ga_{N,0},S\cup D'',i)\\
&=\sum_{i_1+i_2=i}\shg(-N_0,-\ga_{N,0},(S,D''),(i_1,i_2)).
\eeq
Since the disk $D''$ intersects $\ga_N'$ twice, from term (2) of Proposition \ref{thm_2: grading in SHG}, we know that 
$$\shg(-N_0,-\ga_{N,0})=\shg(-N_0,-\ga_{N,0},D'',0).$$
Hence we conclude that
\beq
\shg(-N_0,-\ga_{N,0},S_N,i)&=\sum_{i_1+i_2=i}\shg(-N_0,-\ga_{N,0},(S,D''),(i_1,i_2))\\
&=\shg(-N_0,-\ga_{N,0},S,i).
\eeq

\epf
\brem
Proposition \ref{prop: 2-handle map preserves the grading} is a generalization of \cite[Lemma 2.2]{BLY2020}, where $\al$ is a tangle and $S_N$ is an annulus.
\erem

\subsection{General balanced sutured manifolds}\quad

In this subsection, we prove the main theorem of this section, which is a restatement of the second part of Theorem \ref{thm: formal sutured Floer theory}. 
\bthm\label{thm: Euler characteristic of graded SHG}
Suppose $(M,\ga)$ is a balanced sutured manifold and $\{S_1,\dots,S_n\}$ is a collection of properly embedded admissible surfaces. Then the Euler characteristic
$$\chi(\shg(-M,-\ga,(S_1,\dots,S_n),(i_1,\dots,i_n)))$$
depends only on $(M,\ga)$, $S_1,\dots,S_n$, and $(i_1,\dots,i_n)\in \intg^n$, and is independent of $\hg$. 
\ethm

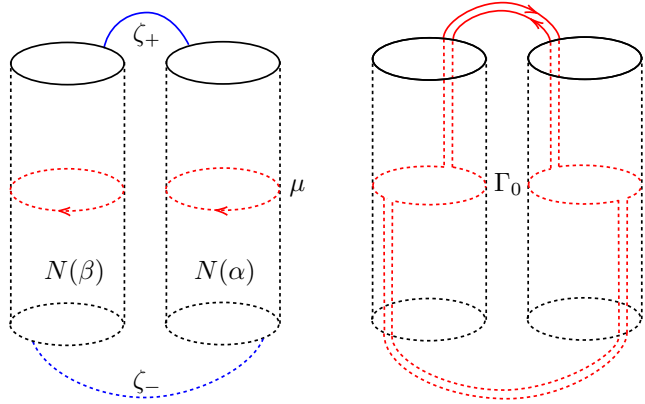
\begin{figure}[ht]
\centering
\begin{overpic}[width=0.6\textwidth]{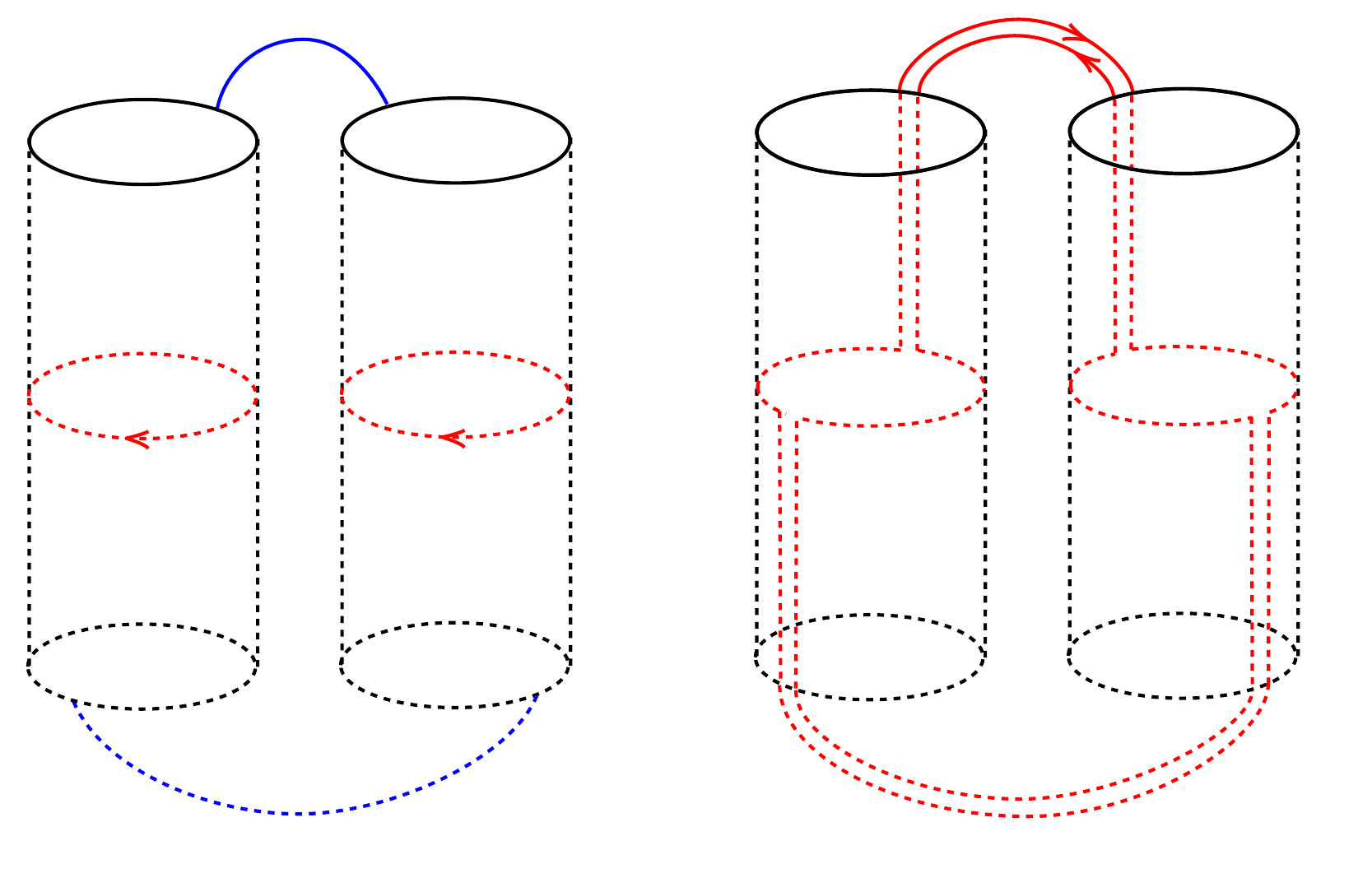}
    \put(20,6){$\zeta_-$}
    \put(20,57){$\zeta_+$}
    \put(7,22){$N(\be)$}
    \put(29,22){$N(\al)$}
    \put(43,35){$\mu$}
    \put(73,35){$\Ga_0$}
\end{overpic}
\vspace{-0.05in}
\caption{The suture $\Ga_{0}$.}\label{fig: the new Ga_0}
\end{figure}

\bcor
Suppose $(M,\ga)$ is a balanced sutured manifold and suppose $H=H_1(M)/{\rm Tors}$. Then the graded Euler characteristic $$\chi_{\rm gr}(\shg(M,\ga))=\chi_{\rm gr}(\shg^g(M,\ga))\in \mathbb{Z}[H]/\pm H$$is independent of the choice of the fixed genus $g$ of closures.
\ecor
\bpf
From Corollary \ref{cor: SFH=SHF} and Theorem \ref{thm: Euler characteristic of graded SHG}, we know$$\chi_{\rm gr}(\shg^g(M,\ga))=\chi_{\rm gr}(SFH(M,\ga))\in \mathbb{Z}[H]/\pm H,$$where the right hand side is independent of the choice of the fixed genus $g$ of closures.
\epf
\bpf[Proof of Theorem \ref{thm: Euler characteristic of graded SHG}]
First we can attach product $1$-handles disjoint from $S_1,\dots,S_n$. From \cite[Section 4.2]{baldwin2016contact}, attaching a product 1-handle does not change the closure and hence does not make any difference to the multi-grading associated to $(S_1,\dots,S_n)$. Hence we can assume that $\ga$ is connected from now on. From \cite[Section 3.1]{LY2020}, we can pick a disjoint union of properly embedded arcs
$$\al=\al_1\cup\cdots\cup\al_m$$
so that
\benu
\item for $k=1,\dots,m$, we have $\partial \al_k\cap R_{+}(\ga)\neq\emptyset {\rm~and~}\partial \al_k\cap R_{-}(\ga)\neq\emptyset,$
\item $M\backslash {\rm int}(N(\al))$ is a handlebody.
\eenu
Then we apply the arguments involved in \cite[Section 3.2]{LY2020}: since $\ga$ is connected, we can pick pairwise disjoint arcs $\zeta_1,\dots,\zeta_m$ so that for any $k=1,\dots,m,$ we have
$$\partial\zeta_k=\partial \al_k{\rm~and~}|\zeta_k\cap \ga|=1.$$
For any $k=1,\dots,g$, let $\be_k\subset \zeta_{k}$ be a neighborhood of the intersection point $\zeta_k\cap\ga$ and let
$$\zeta_{k}\backslash\beta_k=\zeta_{k,+}\cup\zeta_{k,-},$$
where $\zeta_{k,\pm}\subset R_{\pm}(\ga)$.
Push the interior of $\be_k$ into the interior of $M$ to make it a properly embedded arc, which we still call $\be_k$. Let
$$\beta=\beta_1\cup\cdots\cup\beta_m.$$
Let $N=M\backslash {\rm int}(N(\be))$, and let $\ga_N$ be the disjoint union of $\ga$ and a meridian for each component of $\be$. It is explained in \cite[Section 3.2]{LY2020} that $(N,\ga_N)$ can be obtained from $(M,\ga)$ by attaching product 1-handles disjoint from $S_1$,\dots,$S_m$, so there is a canonical identification
$$\shg(-M,-\ga,(S_1,\dots,S_n),(i_1,\dots,i_n))=\shg(-N,-\ga_N,(S_1,\dots,S_n),(i_1,\dots,i_n))$$
Let $H=M\backslash {\rm int}(N(\al\cup\be))$. It is straightforward to check that $H$ is a handlebody. Let $\Ga_{\mu}$ be the disjoint union of $\ga$ and a meridian for each component of $\al\cup\be$. Let the suture $\Ga_0$ be obtained from $\Ga_{\mu}$ by performing band sums along $\zeta_{k,+}$ and $\zeta_{k,-}$ for $k=1,\dots,m$. See Figure \ref{fig: the new Ga_0}. It is straightforward to check that $(N,\ga_N)$ can be obtained from $(H,\Ga_0)$ by attaching contact 2-handles along the meridians of all components of $\al$.

We prove the theorem in the case when $m=1$, while the general case follows from a straightforward induction. If $m=1$, then $\alpha$ is connected. Suppose $\mu$ is the meridian of $\alpha$. As explained in Subsection \ref{subsec: 2-handle attachment}, attaching a contact 2-handle along $\mu$ is the same as performing a $0$-surgery along a push-off $\mu'$ of $\mu$. There is an exact triangle associated to the surgeries along $\mu'$ that is discussed in \cite[Section 3.2]{LY2020} (see also \cite[Section 3.1]{li2019tau}):
\begin{equation}\label{eq: exact triangle, 1}
\xymatrix{
&\shg(-N,-\ga_N)\ar[dr]&\\
\shg(-H,-\Ga_0)\ar[ur]^{C_{\mu}}&&\shg(-H,-\Ga_1)\ar[ll]
}
\end{equation}
The map $C_{\mu}$ is the map associated to the contact 2-handle attachment as discussed in Subsection \ref{subsec: 2-handle attachment}. The suture $\Ga_1$ is obtained from $\Ga_0$ by twisting along $(-\mu)$ once. for $j=1,\dots,n$, let $S_{j,H}=S_j\cap H$. Since $\mu$ is disjoint from $S_{j,H}$ for $j=1,\dots,n$, the proof of Proposition \ref{prop: 2-handle map preserves the grading} implies there is a graded version of the exact triangle (\ref{eq: exact triangle, 1}):
\begin{equation}\label{eq: exact triangle, 2}
\xymatrix{
&\shg(-N,-\ga_N,(S_1,\dots,S_n),(i_1,\dots,i_n))\ar[dd]\\
\shg(-H,-\Ga_0,(S_{1,H},\dots,S_{n,H}),((i_1,\dots,i_n)))\ar[ur]^{C_{\mu}}&\\
&\shg(-H,-\Ga_1,(S_{1,H},\dots,S_{n,H}),((i_1,\dots,i_n)))\ar[lu]
}
\end{equation}
Then Theorem \ref{thm: Euler characteristic of graded SHG} follows from Proposition \ref{prop_2: parity of maps in an exact triangle} and Corollary \ref{cor: Euler char for sutured handlebody with arbitrary surface}.
\epf

\section{The canonical mod 2 grading}\label{sec: mod 2}

Throughout this section, we focus on special cases of balanced sutured manifolds obtained from connected closed 3-manifolds and knots in them (\textit{c.f.} Remark \ref{rem: SFH}). 

\bdefn\label{defn: khg}
Suppose that $Y$ is a closed 3-manifold and $z\in Y$ is a basepoint. Let $Y(1)$ be obtained from $Y$ by removing a 3-ball containing $z$ and let $\delta$ be a simple closed curve on $\partial Y(1)\cong S^2$. Suppose that $K\subset Y$ is a knot and $w$ is a basepoint on $K$. Let $Y(K)$ be the knot complement of $K$ and let $\ga=m\cup (-m)$ consist of two meridians with opposite orientations of $K$ near $w$. Then $(Y(1),\delta)$ and $(Y(K),\ga)$ are balanced sutured manifolds. Define

$$\hhat(Y,z)\deq \sh(Y(1),\delta) ~{\rm and}~\khg(Y,K,w)\deq \sh(Y(K),\ga).$$
\edefn
\begin{conv}
Different choices of the basepoints give isomorphism vector spaces. Since we only care about the isomorphism class of the vector spaces, we omit the basepoints and simply write $\hhat(Y)$ and $\khg(Y,K)$ instead.
\end{conv}
To be more specific and consistent with \cite{LY2020}, in this section, we focus on instanton theory. Based on the discussion in Subsection \ref{subsec: axioms for SH}, we specify the Floer homology $\hg(Y)$ and the cobordism map $\hg(W)$ to be $I^\omega(Y)$ and $I(W,\nu)$. For a connected closed 3-manifold, the framed instanton Floer homology $I^\sharp(Y)$ defined in \cite{kronheimer2011khovanov} is isomorphic to $\hhat(Y)$ when $\hg$ is instanton theory. Hence we replace $\hhat(Y)$ by $I^\sharp(Y)$ throughout this section. Also we replace $\shg$ and $\khg$ by $\shib$ and $\khib$, respectively. Recall that the definitions of $\shib$ and $\khib$ \textit{a priori} depend on the choice of a fixed and large genus $g$ of closures. We write $$\shib^g~{\rm and}~\khib^g$$explicitly in this section. However, for instanton theory, closures of different genus induce isomorphic groups and we can use closures of genus one to define sutured instanton homology (\textit{c.f.} \cite[Section 7]{kronheimer2010knots}).

In this section, we discuss the canonical $\mathbb{Z}_2$-grading on $\khib^g$ and the decomposition of $I^\sharp$ in Theorem \ref{thm: torsion spin c decomposition}.


\label{sec: mod 2 grading of KHI}
\subsection{The case of an unknot}\label{subsec: mod 2 grading for unknot}\quad

In this subsection, we study the model case: the unknot $U$ in $S^3$. Suppose $\mu_U$ and $\lambda_U$ are the meridian and the longitude of $U$, respectively. The knot complement is identified with a solid torus $S^1\times D^2$:
\begin{equation}\label{eq_5: complement of the unknot}
\rho:S^3(U)\xra{\cong}S^1\times D^2,
\end{equation}
where $\rho(\mu_U)=S^1\times \{1\}$, and $\rho(\lambda_U)=\{1\}\times \partial D^2$. For co-prime integers $x$ and $y$, let
$$\ga_{(x,y)}=\ga_{x\lambda_U+y\mu_U}\subset \partial S^3(U)$$
be the suture consisting of two disjoint simple closed curves representing $\pm(x\lambda_U+y\mu_U)$. 
\begin{conv}
Note $\ga_{(x,y)}=\ga_{(-x,-y)}$. From term (4) in Proposition \ref{thm_2: grading in SHG}, the orientation of the suture does not influence the isomorphism type of formal sutured homology. Hence we do not care about the orientation of the suture, and we always assume $y\ge 0$.
\end{conv} 

We describe a closure of the balanced sutured manifold $(S^3(U),\ga_{(x,y)})$ as follows.

Let $\Sigma$ be a connected closed surface of genus $g\geq1$. Suppose $Y_{\Sigma}=S^1\times \Sigma$ and $\Sigma=\{1\}\times \Sigma$. Pick a non-separating simple closed curve $\al\subset \Sigma$ and suppose its complement is $Y_{\Sigma}(\al)=Y_{\Sigma}\backslash{\rm int}(N(\al)).$ There is a framing on $\partial Y_{\Sigma}(\al)$ induced by the surface $\Sigma$. Let $\mu_{\al}$ and $\lambda_{\al}$ be the corresponding meridian and longitude, respectively. Also, suppose $p\in\Sigma$ is a point disjoint from $\al$. According to the discussion in Section \ref{sec: Axioms of sutured Floer homology}, we can form a closure $(\bar{Y},R,\omega)$ of $(S^3(U),\ga_{(x,y)})$ as follows:
\begin{equation}\label{eq_5: closure for unknot}
\bar{Y}=S^{3}(U)\mathop{\cup}_{\phi}Y_{\Sigma}(\al), R=\Sigma, {\rm~and~}\omega=S^1\times\{{\rm pt}\},
\end{equation}
where 
$\phi: \partial S^3(U)\xra{\cong} \partial Y(\al)$
is an orientation reversing diffeomorphism such that 
\begin{equation}\label{phi}
    \phi(x\lambda_{U}+y\mu_U)=\lambda_{\al}.
\end{equation} 
Note that different choices of the preimage of $\mu_{\al}$ lead to different closures of $(S^3(U),\ga_{(x,y)})$. From (\ref{phi}), we know that 
$$\phi(\lambda_U)=z\lambda_{\al}+y\mu_{\al},$$where $z=x$ if $y=0$, $z=y$ is arbitrary if $y=1$, and $zx \equiv 1\pmod y$ in other cases. Again different choices of $z$ lead to different closures. From now on, we fix the value of $z$ as follows: $z=x$ of $y=0$, $z=0$ if $y=1$, and $z$ is the minimal positive integer so that $y|(xz-1)$. Now, composing $\phi$ with the inverse of the map $\rho$ in (\ref{eq_5: complement of the unknot}), suppose
$$\widehat{Y}=Y(\al)\mathop{\cup}_{\phi\circ\rho^{-1}}S^1\times D^2,$$
where
$\phi\circ\rho^{-1}:\partial(S^1\times D^2)\ra\partial Y(\al)$ is a diffeomorphism such that
$$\phi\circ\rho^{-1}(\{1\}\times \partial D^2)=z\lambda_{\al}+y\mu_{\al}.$$
Hence, $\bar{Y}$ is obtained from $Y$ by performing a $y/z$ surgery and we also write
$\bar{Y}=\widehat{Y}_{y/z}.$

\blem\label{lem_5: Euler char of SHI for unknots}
For any suture $\ga_{(x,y)}$ on $\partial S^3(U)$, we have
$$\chi(\shib^g(S^3(U),\ga_{(x,y)}))=\pm y.$$
\elem
\bpf
First, we can focus on the closure $(\bar{Y}=\widehat{Y}_{y/z},R=\Sigma,\omega)$ as in (\ref{eq_5: closure for unknot}). We need to compute the Euler characteristic of
$$\shib^g(S^3(U),\ga_{(x,y)}):= I^{\omega}(\widehat{Y}_{y/z}|\Sigma).$$

If $y=0$, then $x=\pm1$, but $(S^3(U),\ga_{(\pm1,0)})$ are both irreducible and non-taut. By Proposition \ref{prop: SHI detects tautness}, we know that
$$\shib^g(S^3(U),\ga_{(1,0)})=0.$$

If $y=1$, then $z=0$ and $\widehat{Y}_{1/0}=S^1\times \Sigma.$ By Axiom (A1-5), we have
$$\chi(I^{\omega}(\widehat{Y}_{1/0}|\Sigma))=-1.$$

If $y>1$, we have $y>z\geq 1$. If $z=1$, then we have an exact triangle from Axiom (A2)
\begin{equation*}
\xymatrix{
I^{\omega}(\widehat{Y}_{y-1}|\Sigma)\ar[rr]&&I^{\omega}(\widehat{Y}_{y}|\Sigma)\ar[dl]^{f}\\
&I^{\omega}(\widehat{Y}_{1/0}|\Sigma)\ar[ul]&
}
\end{equation*}
where the parity of the map $f$ is odd and those of the rest two are even by Proposition \ref{prop_2: parity of maps in an exact triangle}. Hence we conclude by induction that 
$$\chi(I^{\omega}(\widehat{Y}_{y}|\Sigma))=-y.$$

Finally, when $y>z>1$, suppose the continued fraction of $-y/z$ is
$$-\frac{y}{z}=[a_0,\dots,a_n]=a_0-\frac{1}{a_1-\frac{1}{\dots-\frac{1}{a_n}}},$$where $a_n\le -2$. Define
$$-\frac{y'}{z'}=[a_0,\dots, a_{n-1}] ~{\rm and}~-\frac{y''}{z''}=[a_0,\dots,a_{n-1}+1],$$where $y^\p,y^\pp\ge 0$. From a basic property of continued fraction, we have
$$y=y'+y''~{\rm and~}z=z'+z''.$$From Axiom (A2), there exists an exact triangle
\begin{equation*}
\xymatrix{
I^{\omega}(\widehat{Y}_{y''/z''}|\Sigma)\ar[rr]&&I^{\omega}(\widehat{Y}_{y/z}|\Sigma)\ar[dl]\\
&I^{\omega}(\widehat{Y}_{y'/z'}|\Sigma)\ar[ul]^{f}&
}
\end{equation*}
where the parity of the map $f$ is odd, and those of the rest two are even by Proposition \ref{prop_2: parity of maps in an exact triangle}. 
\epf

\brem\label{rem_5: normalization of the Euler char}
It is worth mentioning that different papers have different normalizations for the canonical $\intg_2$-grading. Our choice of normalization in Axiom (A3) is the same as in \cite{kronheimer2010instanton}. In Lidman, Pinz\'on-Caicedo, and Scaduto's setup \cite{lidman2020framed}, they adapted another normalization and proved $\chi(I^{\omega}(S^1\times \Sigma|\Sigma))=1$ for $\Sigma$ of any genus that is at least one.
\erem

\bcor\label{cor: arbitrary closures for unknot complements}
Suppose $(Y',R',\omega')$ is a closure of $(S^3(U),\ga_{(x,y)})$, then
$$\chi(I(Y|R))=\pm y.$$
\ecor
\bpf
This corollary follows directly from the fact that canonical maps from $I^{\omega}(\widehat{Y}_1|\Sigma)$ to $I^{\omega'}(Y'|R')$ is a composition of cobordism maps and hence is homogeneous. 
\epf

\subsection{Sutured knot complements}\label{subsec: mod grading on knot complement}\quad

Suppose $Y$ is a closed 3-manifold and $K\subset Y$ is a null-homologous knot. Any Seifert surface $S$ of $K$ gives rise to a framing on $\partial Y(K)$: the meridian $\mu$ can be picked as the meridian of the solid torus $N(K)$, and the longitude $\lambda$ can be picked as $S\cap \partial Y(K)$. The `half lives and half dies' fact for 3-manifolds implies that the following map has a 1-dimensional image:
$$\partial_*: H_2(Y(K),\partial Y(K);\mathbb{Q})\ra H_1(\partial Y(K);\mathbb{Q}).$$
Hence any two Seifert surfaces lead to the same framing on $\partial Y(K)$. 
\bdefn\label{defn_4: surgery on pairs}
The framing $(\lambda,-\mu)$ defined as above is called the \textbf{canonical framing} of $(Y,K)$. With this canonical framing, let $$\ga_{(x,y)}=\ga_{x\lambda+y\mu}\subset \partial Y(K)$$be the suture consisting of two disjoint simple closed curves representing $\pm(x\lambda+y\mu)$. 
\edefn

Our goal in this subsection is to define a canonical $\mathbb{Z}_2$-grading on $\shib^g(Y(K),\ga_{(x,y)})$ for any fixed large enough $g$. Recall $\shib^g(M,\ga)$ is the projectively transitive system formed by closures of $(M,\ga)$ of a fixed genus $g$. We first assign a $\intg_2$-grading for any closure of $(Y(K),\ga_{(x,y)})$.

Suppose $(\bar{Y},R,\omega)$ is a closure of $(Y(K),\ga_{(x,y)})$. Then we can form a closure $(\bar{Y}_U,R,\omega)$ of $(S^3(U),\ga_{(x,y)})$ by taking
\begin{equation}\label{eq_5: corresponding closure for unknot}
\bar{Y}_U=\bar{Y}\backslash({\rm int}(Y(K)))\mathop{\cup}_{\rm id}S^3(U).
\end{equation}
Here ${\rm id}$ is the diffeomorphism between toroidal boundaries, which respect the canonical framings on both boundaries.

\bdefn\label{defn_5: modified mod 2 grading}
The {\bf modified $\intg_2$-grading} on $I^{\omega}(\bar{Y}|R)$ is defined as follows.
\begin{enumerate}[(1)]
\item If $\chi(I^{\omega}(\bar{Y}_U|R))$ is negative, then the grading is defined by the canonical $\intg_2$-grading on $I^{\omega}(\bar{Y}|R)$.
    \item If $\chi(I^{\omega}(\bar{Y}_U|R))$ is positive, then the grading is defined by switching the odd and even parts of $I^{\omega}(\bar{Y}|R)$ with the canonical $\intg_2$-grading.
\end{enumerate}
\edefn
 Suppose $(\bar{Y},R,\omega)$ and $(\bar{Y}',R,\omega)$ are two closures of $(Y(K),\ga_{(x,y)})$ so that $\bar{Y}'$ is obtained from $\bar{Y}$ by a Dehn surgery along a curve $\be\subset\bar{Y}$, which is disjoint from ${\rm int}(M)$, $R$ and $\omega$. Then there is a map
$$F:I^{\omega}(\bar{Y}|R)\ra I^{\omega}(\bar{Y}'|R)$$
associated to the Dehn surgery along $\be\subset \bar{Y}$. Let $(\bar{Y}_U,R,\omega)$ and $(\bar{Y}'_U,R,\omega)$ be the closures of $(S^3(U),\ga_{(x,y)})$ constructed as in (\ref{eq_5: corresponding closure for unknot}). There is also a map
$$F_U:I^{\omega}(\bar{Y}_U|R)\ra I^{\omega}(\bar{Y}'_U|R)$$
associated to the same Dehn surgery along $\be\subset \bar{Y}_U$. Then we have the following.

\blem\label{lem_5: maps arising from surgery have the same parity}
The maps $F$ and $F_U$ have the same parity with respect to the canonical $\intg_2$-gradings on corresponding instanton Floer homologies.
\elem
\bpf
Note that $H_1(S^3(U);\mathbb{Q})\cong \mathbb{Q}\langle \mu_U\rangle$ and the map
$$i^U_*:H_1(\partial S^3(U);\mathbb{Q})\ra H_1(S^3(U);\mathbb{Q})$$
induced by the inclusion has a 1-dimensional kernel generated by $\lambda_U$. For a null-homologous knot $K\subset Y$, we know that the map
$$i_*:H_1(\partial Y(K);\mathbb{Q})\ra H_1(Y(K);\mathbb{Q})$$
induced by the inclusion has a 1-dimensional kernel generated by the longitude $\lambda$ of $K$ and has a 1-dimensional image generated by the meridian $\mu$ of $K$.
Hence, from the Mayer-Vietoris sequence, we know that there is an injective map
$$j:H_1(\bar{Y}_U;\mathbb{Q})\hookrightarrow H_1(\bar{Y};\mathbb{Q}),$$
that sends $[\mu_U]$ to $[\mu]$ and sends every homology class in $\bar{Y}_U\backslash S^3(U)=\bar{Y}\backslash Y(K)$ using the natural map
$$i'_*:H_1(\bar{Y}\backslash Y(K);\mathbb{Q})\ra H_1(\bar{Y};\mathbb{Q}).$$

Similarly, since $\be\cap {\rm int}(Y(K))=\emptyset$, we know that there is an injective map
$$j^{\be}:H_1(\bar{Y}_U(\be);\mathbb{Q})\hookrightarrow H_1(\bar{Y}(\be);\mathbb{Q}),$$
which fits into the following commutative diagram
\begin{equation*}
\xymatrix{
H_1(\partial \bar{Y}_U(\be);\mathbb{Q})\ar[rr]^{\iota^{U}_*}\ar[d]^{=}&&H_1(\bar{Y}_U(\be);\mathbb{Q})\ar[d]^{j^{\be}}\\
H_1(\partial \bar{Y}(\be);\mathbb{Q})\ar[rr]^{\iota_*}&&H_1(\bar{Y}(\be);\mathbb{Q})
}    
\end{equation*}
where $\iota_*$ and $\iota_*^{U}$ are induced by natural inclusions. Hence we conclude, under the identification 
$H_1(\partial \bar{Y}_U(\be);\mathbb{Q})=H_1(\partial \bar{Y}(\be);\mathbb{Q}),$ two kernels are also identified:
$${\rm ker}(\iota_*^{U})={\rm ker}(\iota_*).$$
Since $F$ and $F_U$ are associated to Dehn surgeries along $\be$ of the same slopes, we conclude from \ref{prop_2: parity of maps in an exact triangle} that their parity must be the same.
\epf

\blem\label{lem_5: canonical maps are even}
Suppose $(\bar{Y},R,\omega)$ and $(\bar{Y}',R',\omega')$ are two different closures of $(Y(K),\ga_{(x,y)})$ so that $g(R)=g(R')$. There is a canonical isomorphism
$$\Phi:I^{\omega}(\bar{Y}|R)\xra{\cong} I^{\omega'}(\bar{Y}'|R')$$
as in Definition \ref{defn: canonical maps}.
Then with respect to the modified $\intg_2$-grading as in Definition \ref{defn_5: modified mod 2 grading}, the canonical map $\Phi$ is grading preserving.
\elem

\bpf
From the definition of $\Phi$ and Axiom (A3-3), we know that $\Phi$ is always homogeneous. For two closures $(\bar{Y},R,\omega)$ and $(\bar{Y}',R',\omega')$ of $(Y(K),\ga_{(x,y)})$, we can form two corresponding closures $(\bar{Y}_U,R,\omega)$ and $(\bar{Y}'_U,R',\omega')$ of $(S^3(U),\ga_{(x,y)})$ as in (\ref{eq_5: closure for unknot}), respectively. There is a canonical map
$$\Phi_U:I^{\omega}(\bar{Y}_U|R)\xra{\cong} I^{\omega'}(\bar{Y}'_U|R').$$

From Definition \ref{defn_5: modified mod 2 grading}, the modified $\intg_2$-gradings on $I^{\omega}(\bar{Y}|R)$ and $I^{\omega}(\bar{Y}^\p|R^\p)$ coincide if and only if the canonical $\intg_2$-gradings on $I^{\omega}(\bar{Y}_U|R)$ and $I^{\omega}(\bar{Y}^\p_U|R^\p)$ coincide. The definition of the modified $\intg_2$-grading automatically makes $\Phi_U$ even under the modified $\intg_2$-grading. Hence to show that $\Phi$ is grading preserving, it suffices to show that $\Phi$ and $\Phi_U$ have the same parity under the modified $\intg_2$-grading.

From the construction of the canonical maps, there is a sequence of simple closed curves $\be_1,\dots,\be_n$ on $R$, such that the map $\Phi$ is the composition of cobordism maps induced by a diffeomorphism and the sequence of Dehn surgeries. Similarly, the map $\Phi_U$ is the composition of the maps induced another diffeomorphism and the same sequence of Dehn surgeries on $\bar{Y}_U$. Since the surgery curves are all on $R$ and disjoint from ${\rm int}(Y(K))$, cobordism maps induced by diffeomorphisms are always with even degrees, \textit{i.e.} preserving the $\mathbb{Z}_2$-grading. For Dehn surgeries along $\be_i$, we can apply Lemma \ref{lem_5: maps arising from surgery have the same parity} and then finish the proof.
\epf
\bdefn\label{defn: mod 2 grading on SHI^g}
Suppose $(Y(K),\ga_{(x,y)})$ is the balanced sutured manifold constructed as before and suppose $g$ is the fixed large enough genus of closures. By Lemma \ref{lem_5: canonical maps are even}, we can define the \textbf{canonical $\mathbb{Z}_2$-grading} on $\shib^g(Y(K),\ga_{(x,y)})$ by the modified $\mathbb{Z}_2$-grading on the closures. In particular, there is a canonical $\mathbb{Z}_2$-grading on $\khib^g(Y,K)$.
\edefn
\bprop\label{prop: mod 2 grading on SHI^g}
The canonical $\mathbb{Z}_2$-grading on $\shib^g(Y(K),\ga_{(x,y)})$ is independent of the large enough genus $g$.
\eprop
\bpf
We need to compare the $\intg_2$-grading for closures of different genera. First we deal with the case of unknot. As in Subsection \ref{subsec: mod 2 grading for unknot}, we can constructed a standard closure $(\bar{Y}_{y/z},\Sigma)$ for $(Y(K),\ga_{(x,y)})$. Here, the genus of $\Sigma$ can be arbitrary. To specify the genus, in the proof of this corollary, we temporarily write $(\bar{Y}_{y/z},\Sigma)$ as $(\bar{Y}_{y/z}^g,\Sigma_g)$. In \cite{baldwin2015naturality}, the canonical map
$$\Phi^{g,g+1}:I^{\omega}(\bar{Y}_{y/z}^{g}|\Sigma_{g})\ra I^{\omega}(\bar{Y}_{y/z}^{g+1}|\Sigma_{g+1})$$
is constructed as follows: recall $\al\subset \Sigma_g$ is a non-separating simple closed curve, and
$$\bar{Y}_{y/z}^{g}=Y(K)\cup (S^1\times\Sigma_g)\backslash N(\al).$$
Let $\be\subset\Sigma_g$ be another simple closed curve so that $\al\cap\be=\emptyset$. Take a curve $\theta\subset \Sigma_2$ be non-separating simple closed curve as well. Then we can form $(\bar{Y}_{y/z}^{g+1},\Sigma_{g+1})$ from $(\bar{Y}_{y/z}^{g},\Sigma_{g})$ and $(S^1\times \Sigma_2,\Sigma_2)$ by cutting them open along $S^1\times\be$ and $S^1\times \theta$ respectively, and then glue the two pieces together along toroidal boundaries by the identifying the $S^1$ factor and $\beta=\theta$. Then, as in \cite{kronheimer2010knots}, there is a cobordism $W^{g,g+1}$ from $(\bar{Y}_{y/z}^{g},\Sigma_{g})$ and $(S^1\times \Sigma_2,\Sigma_2)$ to $(\bar{Y}_{y/z}^{g+1},\Sigma_{g+1})$, known as the Floer excision cobordism. In the proof of \cite[Proposition 3.6]{lidman2020framed} and \cite[Lemma 3.3]{kronheimer2010instanton}, the degree of the cobordism induced by cobordisms constructed in the same way as $W^{g,g+1}$ has been computed explicitly. By a similar computation, we know that the canonical map $\Phi^{g,g+1}$ preserves the modified $\intg_2$-grading (Note by the above argument, the modified $\intg_2$-grading for $(\bar{Y}_{y/z}^{g},\Sigma_{g})$ is the same as the canonical grading). For any two closures of $(Y(K),\ga_{(x,y)})$, as in \cite{baldwin2015naturality}, the canonical maps are constructed by composing the maps induced by some $W^{g,g+1}$ and canonical maps for closures of the same genera. Since both types of maps are grading preserving, we conclude that any canonical map preserves the modified $\mathbb{Z}_2$-grading.
\epf
\brem
The proof of Proposition \ref{prop: mod 2 grading on SHI^g} does not apply to other Floer-type theories $\hg$ because we need to use Floer's excision theorem along a surface of genus one.
\erem
\subsection{Computations and applications}\label{subsec: computing Euler char}\quad

In this subsection, we do some calculations based on techniques introduced before.

\begin{figure}[ht]
\centering
\includegraphics[width=0.45\textwidth]{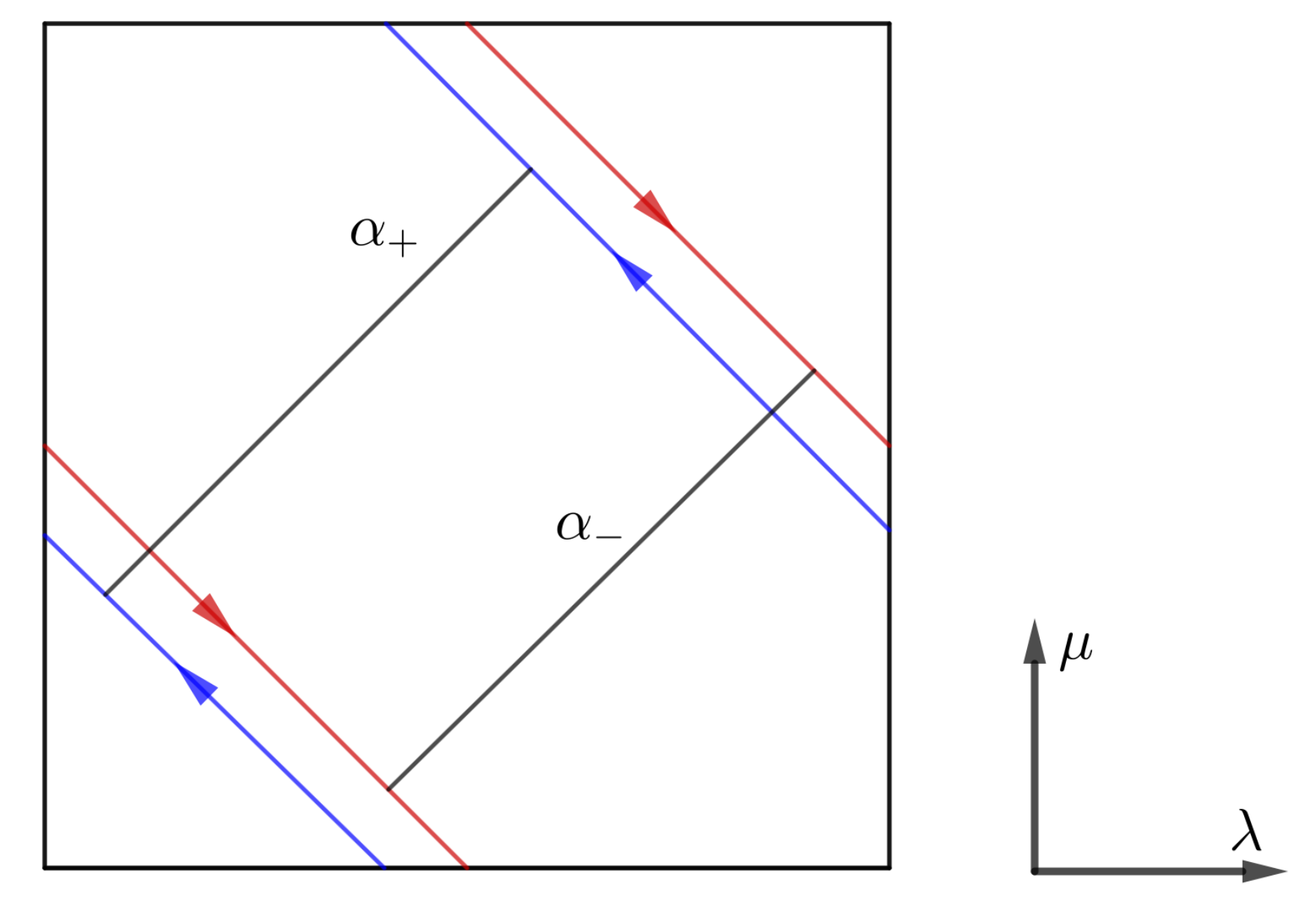}
\caption{Bypass arcs on $\ga_{(1,-1)}$.}
\label{bypass_arc}
\end{figure}

At first, we deal with bypass exact triangles. Suppose $y_3/x_3$ is a surgery slope with $y_3\ge 0$. According to Honda \cite{honda2000classification}, there are two basic bypasses on the balanced sutured manifold $(Y(K),\ga_{(x_3,y_3)})$, whose arcs are depicted as in Figure \ref{bypass_arc}. The sutures involved in the bypass triangles were described explicitly in Honda \cite{honda2000classification}. 
\begin{defn}\label{defn: xi and yi}
For a surgery slope $y_1/x_1$ with $y_1\ge 0$, suppose its continued fraction is
$$
\frac{y_1}{x_1}=[a_0,a_1,\dots,a_n]=a_0-\frac{1}{a_1-\frac{1}{\dots-\frac{1}{a_n}}},$$where $a_i\le -2$. If $y_1>-x_1> 0$, let
$$\frac{y_2}{x_2}=[a_0,\dots,a_{n-1}]~{\rm and~}\frac{y_3}{x_3}=[a_0,\dots,a_{n}+1].
$$
If $-x_1>y_1>0$, we do the same thing for $x_1/(-y_1)$. If $y_1>x_1>0$, we do the same thing for $y_1/(-x_1)$. If  $x_1>y_1>0$, we do the same thing for $x_1/(-y_1)$. If $y_1/x_1=1/0$, then $y_2/x_2=0/1$ and $y_3/x_3=1/(-1)$. If $y_1/x_1=0/1$, let $y_2/x_2=1/(-1)$ and $y_3/x_3=0/1$.
We always require that $y_2\ge 0$ and $y_3\ge 0$. 
\end{defn}
\brem
It is straightforward to use induction to verify that for $y_1>-x_1> 0$,
$$
x_1=x_2+x_3,~{\rm and}~y_1=y_2+y_3.$$
\erem
Then the bypass exact triangle in Theorem \ref{thm_2: bypass exact triangle on general sutured manifold} becomes the following.

\begin{prop}\label{prop_4: bypass exact triangle on knot complement}
Suppose $K\subset Y$ is a null-homologous knot, and suppose the surgery slopes $y_i/x_i$ for $i=1,2,3$ are defined as in Definition \ref{defn: xi and yi}. Suppose $\psi_{+,*}^*$ and $\psi_{-,*}^*$ are from two different bypasses, where $*$ means the corresponding slope. Then there are two exact triangles about $\psi_{+,*}^*$ and $\psi_{-,*}^*$, respectively.
\begin{equation*}
\xymatrix{
\shib^g(-Y(K),-\ga_{(x_2,y_2)})\ar[rr]^{\psi_{\pm,y_3/x_3}^{y_2/x_2}}&&\shib^g(-Y(K),-\ga_{(x_3,y_3)})\ar[dl]^{\psi_{\pm,y_1/x_1}^{y_3/x_3}}\\
&\shib^g(-Y(K),-\ga_{(x_1,y_1)})\ar[lu]^{\psi_{\pm,y_2/x_2}^{y_1/x_1}}&
}
\end{equation*}
\end{prop}
As stated in Proposition \ref{prop_2: description of bypass maps}, the bypass maps $\psi_1$, $\psi_2$, and $\psi_3$ are induced by some cobordism maps. Then we have the following.

\blem\label{lem_5: bypass and mod 2 grading}
Suppose $g$ is a large enough integer and suppose $y_i/x_i$ for $i=1,2,3$ is from Definition \ref{defn: xi and yi}. Suppose further $$x_1=x_2+x_3\aand y_1=y_2+y_3.$$With respect to the canonical $\intg_2$-grading on $\shib^g$ in Definition \ref{defn: mod 2 grading on SHI^g}, the parity of the map $\psi_2$ is odd and those of the rest two are even. As a consequence,
$$\chi(\shib^g(-Y(K),-\ga_{(x_1,y_1)}))=\chi(\shib^g(-Y(K),-\ga_{(x_3,y_3)}))+\chi(\shib^g(-Y(K),-\ga_{(x_2,y_2)})).$$
\elem

\bpf
As in Proposition \ref{prop_2: description of bypass maps}, we can fix a large enough $g$ so that for $i=1,2,3$, there are closures $(\bar{Y}_{i},R,\omega)$ for $(-Y(K),-\ga_{(x_i,y_i)})$ of genus $g$, and the bypass maps $\psi_1$, $\psi_2$, and $\psi_3$ have the same $\mathbb{Z}_2$ degree because the maps induced by Dehn surgeries along three curves $\zeta_1$, $\zeta_2$, and $\zeta_3$ in corresponding closures $\bar{Y}_i\backslash{\rm int}(Y(K))$. Since we only care about the $\mathbb{Z}_2$ degrees of maps, in a slight abuse of notation, we do not distinguish the bypass map and the map induced by Dehn surgery.

For $i=1,2,3$, we can form corresponding closures $(\bar{Y}_{i}^U,R,\omega)$ as in (\ref{eq_5: closure for unknot}) so that the curves $\zeta_1$, $\zeta_2$, and $\zeta_3$ still lie in closures. Moreover, suitable surgeries along these curves induces an exact triangle
\begin{equation*}
\xymatrix{
\shib^g(-S^3(U),-\ga_{(x_3,y_3)})\ar[rr]^{\psi^U_{3}}&&\shib^g(-S^3(U),-\ga_{(x_1,y_1)})\ar[dl]^{\psi^U_{1}}\\
&\shib^g(-S^3(U),-\ga_{(x_2,y_2)})\ar[lu]^{\psi^U_{2}}&
}
\end{equation*}

As in the proof of Lemma \ref{lem_5: canonical maps are even}, with the help of Lemma \ref{lem_5: maps arising from surgery have the same parity}, it suffices to check that the parities of maps $\psi_1^U$, $\psi_2^U$, and $\psi_3^U$ are odd or even as claimed, with respect to the canonical $\intg_2$-grading on $\shib^g$ from Definition \ref{defn: mod 2 grading on SHI^g}. 

For the case of the unknot, the argument becomes straightforward: from Definition \ref{defn_5: modified mod 2 grading} and Lemma \ref{lem_5: Euler char of SHI for unknots}, we know that for any $y>0$,
$$\chi(\shib^g(-S^3(U),-\ga_{(x,y)}))=-y$$
Then the equation $y_1=y_2+y_3$ implies that
$$\chi(\shib^g(-S^3(U),-\ga_{(x_1,y_1)}))=\chi(\shib^g(-S^3(U),-\ga_{(x_3,y_3)}))+\chi(\shib^g(-S^3(U),-\ga_{(x_2,y_2)})).$$

Note that the maps $\psi_i^U$ for $i=1,2,3$ are coming from a real surgery exact triangle as in Proposition \ref{prop_2: description of bypass maps}, while the $\mathbb{Z}_2$-gradings on $\shib^g$ could possibly be shifted due to the normalization in Definition \ref{defn_5: modified mod 2 grading} and the surgery along curves $\eta_1$ and $\eta_2$ as in Proposition \ref{prop_2: description of bypass maps}. Hence they still satisfies the hypothesis of Lemma \ref{lem_2: Floer's exact triangle after grading shift}. Thus, we conclude that the parity of $\psi_2^U$ is odd and those of the other two are even. Similarly, the parity of $\psi_2$ is odd and those of the other two are even, and we have
$$\chi(\shib^g(-Y(K),-\ga_{(x_1,y_1)}))=\chi(\shib^g(-Y(K),-\ga_{(x_3,y_3)}))+\chi(\shib^g(-Y(K),-\ga_{(x_2,y_2)})).$$
\epf

Let $Y$ be a closed 3-manifold and let $K\subset Y$ be a null-homologous knot. Suppose $S$ is a minimal genus Seifert surface of $K$. Its genus is always denoted by $g(S)$, which is distinguished with $g$, the fixed genus of closures. We refer \cite[Section 4]{LY2020} for the definitions of sutures $\Ga_n,\Ga_n(y/x)$, the admissible surface with stablization $S^{\tau}$, the bypass maps $\psi_{+,*}^*,\psi_{-,*}^*$, and numbers $i_{max}^*,i_{min}^*$. To simplify our notation, we write
\begin{equation}\label{eq_5: Euler char on knot complements, 1}
    \chi^{g}_{y/x}(-Y,K,i)=\chi(\shib^g(-Y(K),-\ga_{(x,y)},S^{\tau},i)),
\end{equation}
where the Euler characteristic is with respect to the canonical $\intg_2$-grading on $\shib^g$ as in Definition \ref{defn: mod 2 grading on SHI^g}.
We write
\begin{equation}\label{eq_5: Euler char on knot complements, 2}
    \chi^{g}_{y/x}(-Y,K)=\sum_{i\in\intg}\chi^{g}_{y/x}(-Y,K,i)
\end{equation}
When $|x|=1$, we write $y/x$ as an integer. 
Also, we write
$$\chi_{\mu}^g(-Y,K,i)=\chi_{1/0}^g(-Y,K,i)$$
to specify the meridional suture.

\blem\label{lem_5: Euler char for 0,1}
Suppose $Y$ is a closed oriented 3-manifold, and $K\subset Y$ is a null-homologous knot. For $g\in\intg$ large enough and any $i\in\intg$, we have
$$\chi_{1}^g(-Y,K,i)=\chi_{\mu}^g(-Y,K,i)\aand \chi_0^g(-Y,K,i)=0.$$
\elem
\bpf
From \cite[Proposition 4.14]{LY2020}, we have the following two bypass exact triangles:
\begin{equation*}
    \xymatrix{
    \shib^g(-Y(K),-\Ga_0,S^{\tau},i)\ar[rr]^{\psi_{-,1}^0}&&\shib^g(-Y(K),-\Ga_1,S^{\tau},i)\ar[dll]^{\psi_{-,\mu}^1}\\
    \shib^g(-Y(K),-\Ga_0,S^{\tau},i)\ar[u]^{\psi_{-,0}^{\mu}}&&
    }
\end{equation*}
and
\begin{equation*}
    \xymatrix{
    \shib^g(-Y(K),-\Ga_0,S^{\tau},i+1)\ar[rr]^{\psi_{+,1}^0}&&\shib^g(-Y(K),-\Ga_1,S^{\tau},i)\ar[dll]^{\psi_{+,\mu}^1}\\
    \shib^g(-Y(K),-\Ga_0,S^{\tau},i)\ar[u]^{\psi_{+,0}^{\mu}}&&
    }
\end{equation*}
Hence we obtain the following two equations from Lemma \ref{lem_5: bypass and mod 2 grading}
\begin{equation*}\label{eq_5: relating 0,1,mu sutures, -}
\chi_{1}^{g}(-Y,K,i)=\chi_{0}^{g}(-Y,K,i)+\chi_{\mu}^{g}(-Y,K,i)    
\end{equation*}
\begin{equation*}\label{eq_5: relating 0,1,mu sutures, +}
\chi_{1}^{g}(-Y,K,i)=\chi_{0}^{g}(-Y,K,i+1)+\chi_{\mu}^{g}(-Y,K,i)    
\end{equation*}
By Axiom (A1-4), for $i>g(S)$, we have
$$\shib^g(-Y(K),-\Ga_0,S^{\tau},i)=0.$$
Hence we conclude by (\ref{eq_5: Euler char on knot complements, 1}) and (\ref{eq_5: Euler char on knot complements, 2}) that
$$\chi_{1}^{g}(-Y,K,g(S))=\chi_{\mu}^{g}(-Y,K,g(S))\aand \chi_{0}^{g}(-Y,K,g(S))=0.$$The lemma follows from the induction on the grading $i$.
\epf

\blem\label{lem_5: Euler char for general (x,y)}
Suppose $Y$ is a closed oriented 3-manifold, and $K\subset Y$ is a null-homologous knot. For the suture $\ga_{(x,y)}$ with $y>0$, we know that
$$\chi_{y/x}^g(-Y,K,i)=\sum_{j=0}^{y-1}\chi_{\mu}^g(-Y,K,i-i^{y}_{max}+i^{\mu}_{max}+j).$$
\elem
\bpf
We only prove the case when $x<0$. The other case is similar. First, if $x=1$, then we have a bypass exact triangle (in this case we write $y=n$)
\begin{equation*}
\xymatrix{
\shib^g(-Y(K),-{\Ga}_{n},S^{\tau},i)\ar[rr]^{\quad\quad\psi^{n}_{-,n+1}}&&\shib^g(-Y(K),-{\Ga}_{n+1},S^{\tau},i)\ar[dll]^{\psi^{n+1}_{-,\mu}}\\
\shib^g(-Y(K),-{\Ga}_{\mu},S^{\tau},i+1)\ar[u]^{\psi^{\mu}_{-,n}}&&
}    
\end{equation*}
Hence, we can apply Lemma \ref{lem_5: bypass and mod 2 grading}, Lemma \ref{lem_5: Euler char for 0,1}, and the induction to conclude that
$$\chi_{n}^g(-Y,K,i)=\sum_{j=0}^{n-1}\chi_{\mu}^g(-Y,K,i-i^{n}_{max}+i^{\mu}_{max}+j).$$

If $x>1$, we can use the continued fraction description of $y/x$ and apply an induction in the same spirit as in the proof of Lemma \ref{lem_5: Euler char of SHI for unknots}.
\epf

\bcor\label{cor_5: Euler char for gradings in the middle}
Suppose $Y$ is a closed oriented 3-manifold, and $K\subset Y$ is a null-homologous knot. For the suture $\ga_{(x,y)}$ where $y>2g(S)$, and for any $i\in\intg$ so that
$$i^{y}_{max}-2g(S)\geq i\geq i^y_{min}+2g(S),$$
we know that
$$\chi^g_{y/x}(-Y,K,i)=\chi^g_{\mu}(-Y,K).$$
\ecor
\bpf
The corollary follows immediately from Lemma \ref{lem_5: Euler char for general (x,y)} and the fact that there are only $(2g(S)+1)$ gradings with nontrivial elements for $\shib^g(-Y(K),-\Ga_{\mu},S^{\tau})$.
\epf

\blem\label{lem_5: Euler char for mu and I sharp}
Suppose $Y$ is a closed oriented 3-manifold, and $K\subset Y$ is a null-homologous knot. Then we have
$$\chi_{\mu}^g(-Y,K)=\pm\chi(I^{\sharp}(-Y)).$$
\elem
\bpf
From Lemma \ref{lem_5: Euler char for general (x,y)}, we know that for any $n\in\intg_{\ge 0}$, we have
$$\chi_{n}^g(-Y,K)=n\cdot \chi_{\mu}^g(-Y,K).$$
From \cite[Lemma 4.9]{LY2020}, we know that there is an exact triangle
\begin{equation*}
\xymatrix{
\shib^g(-Y(K),-\Ga_{n})\ar[rr]&&\shib^g(-Y(K),-\Ga_{n+1})\ar[dl]\\
&I^{\sharp}(-Y)\ar[ul]&
}    
\end{equation*}
Hence, by Lemma \ref{lem_2: Floer's exact triangle after grading shift}, we know that there is a proper sign assignment for all $n$ so that
$$\pm\chi_{n}^g(-Y,K)\pm \chi_{n+1}^g(-Y,K)\pm\chi(I^{\sharp}(Y))=0.$$
Hence the only possibilities are 
$$\chi(I^{\sharp}(Y))=\pm\chi_{\mu}^g(-Y,K).$$
\epf

\bpf[Proof of Proposition \ref{prop: Euler Char of I sharp}]
It is an immediate corollary following Corollary \ref{cor_5: Euler char for gradings in the middle}, Lemma \ref{lem_5: Euler char for mu and I sharp} and the definition of the decomposition from \cite[Section 4.3]{LY2020}.
\epf

For a knot $K$ in $S^3$, we can actually fix the sign ambiguity coming from different choices of the fixed genus of the closures.
\blem\label{lem_5: Euler char for knots in S^3, mu}
For any knot $K\subset S^3$ and any positive integer $g$, we have
$$\chi^g_{\mu}(-S^3,K)=-1.$$
\elem
\bpf
Since we adapt the normalization from Kronheimer and Mrowka \cite[Section 2.6]{kronheimer2010instanton}, we can directly apply the results from them. In particular, for any knot $K$, in \cite[Section 2.4]{kronheimer2010instanton}, a preferred closure $(\bar{Y}_1,\Sigma_1,\omega_1)$ of $(-S^3(K),-\Ga_{\mu})$ with $g(\Sigma_1)=1$ is chosen. Then they proved that 
$$\chi(I^{\omega_1}(\bar{Y}_1|\Sigma_1))=-\Delta_K(1)=-1.$$
Note that this coincide with our choice of modified $\intg_2$-grading: when the Euler characteristic is negative, we do not perform any shift. 

The case $g=2$ has already been studied in the proof of \cite[Lemma 3.3]{kronheimer2010instanton}: another preferred closure $(\bar{Y}_2,\Sigma_2,\omega_2)$ with $g(\Sigma_2)=2$ for $(-S^3(K),-\Ga_{\mu})$ is constructed, and there is a cobordism $W_1$ from $\bar{Y}_1\cup S^1\times \Sigma_2$ to $\bar{Y}_2$ coming from Floer's exicison theorem. The canonical generator of $I^{\omega}(S^1\times \Sigma_2|\Sigma_2)$ is proved to be at the odd grading (\textit{c.f.} \cite[Lemma 3.8]{lidman2020framed}, though the normalization of the canonical $\intg_2$-grading is different). Moreover, the degree of the cobordism map $W_1$ is odd (\textit{c.f.} \cite[Proposition 3.6]{lidman2020framed}). 

For general $g$, it is straightforward to generalize the above construction for $\bar{Y}_1$ and $\bar{Y}_2$ to $\bar{Y}_g$ and $\bar{Y}_{g+1}$. There is a similar cobordism $W_g$ from $\bar{Y}_g\cup S^1\times \Sigma_2$ to $\bar{Y}_{g+1}$, the degree of which can be computed easily to be odd. Hence by induction we conclude that for all $g$,
$$\chi^g_{\mu}(-S^3,K)=-1.$$
\epf

By Lemma \ref{lem_5: Euler char for knots in S^3, mu}, we can identify $\chi^g_{\mu}(-S^3,K)$ for all large enough $g$, we simply write $\chi_{\mu}(-S^3,K)=\chi^g_{\mu}(-S^3,K)$ instead. Applying Lemma \ref{lem_5: Euler char for general (x,y)}, we know that for any $g$ large enough and $y>0$,
$$\chi^g_{y/x}(-S^3,K)=-y.$$
Similarly, we simply write $\chi_{y/x}(-S^3,K)$ instead.

Finally, we consider the projectively transitive system $\shi(M,\ga)$ for a balanced sutured manifold $(M,\ga)$ defined in \cite{baldwin2015naturality}, which is independent of the choices of the genus of the closures. The isomorphism class of $\shi(M,\ga)$ and $\shib^g(M,\ga)$ are the same. Similar to $\shib^g(M,\ga)$, it has a decomposition associated to an admissible surface $S\subset(M,\ga)$.

\bdefn\label{defn_3: shifting the grading}
Suppose $(M,\ga)$ is a balanced sutured manifold and $S$ is an admissible surface in $(M,\ga)$. For any $i,j\in\intg$, define
$$\shi(M,\ga,S,i)[j]=\shi(M,\ga,S,i-j).$$
\edefn

In \cite[Section 5]{li2019direct}, the first author constructed a minus version of the instanton knot homology via the direct system
\begin{equation}\label{eq_5: direct system}
\cdots\ra\shi(-S^3(K),\Ga_n,S^{\tau})[g(K)-i_{max}^n]\xra{\psi_{-,n+1}^n}\shi(-S^3(K),\Ga_{n+1},S^{\tau})[g(K)-i_{max}^{n+1}]\ra\cdots
\end{equation}
and define $\khi(-S^3,K)$ to be the direct limit of (\ref{eq_5: direct system}). All $\psi_{-,n+1}^n$ are grading preserving after shifting, so there is a well-defined $\intg$-grading (the Alexander grading) on $\khi(-S^3,K)$, which we write as
$$\khi(-S^3,K,i).$$ 

By \cite[Corollary 2.21]{li2019direct}, there is a commutative diagram
\begin{equation*}
\xymatrix{
\shi(-S^3(K),\Ga_n,S^{\tau})[g(K)-i_{max}^n]\ar[r]^{\psi_{-,n+1}^n}\ar[d]^{\psi_{+,n+1}^n}&\shi(-S^3(K),\Ga_{n+1},S^{\tau})[g(K)-i_{max}^{n+1}]\ar[d]^{\psi_{+,n+2}^{n+1}}\\
\shi(-S^3(K),\Ga_{n+1},S^{\tau})[g(K)-i_{max}^{n+1}]\ar[r]^{\psi_{-,n+2}^{n+1}}&\shi(-S^3(K),\Ga_{n+2},S^{\tau})[g(K)-i_{max}^{n+2}]
}    
\end{equation*}
Hence the maps $\{\psi_{+,n+1}^n\}$ induces an map $U$ on $\khi(-S^3,K)$. 
\begin{proof}[Proof of Proposition \ref{prop: Euler char of khi-minus}]
From Lemma \ref{lem_5: bypass and mod 2 grading}, the parity of maps $\psi_{-,n+1}^n$ are all even, hence there is a well-defined $\intg_2$-grading on $\khi(-S^3,K)$. 
Again from Lemma \ref{lem_5: bypass and mod 2 grading}, we know that the parity of the map $U$ is even, \textit{i.e.}, preserving the $\intg_2$-grading on $\khi(-S^3,K)$.
Finally, we can apply Lemma \ref{lem_5: Euler char for general (x,y)} and the fact$$\chi(\khib^g(-S^3,K))=-\Delta_K(t)$$ to conclude the desired formula. Note that by our normalization, the sign is negative.
\end{proof}

\bibliographystyle{alpha}

\end{document}